\documentclass[11pt]{article}
\usepackage[english]{babel}
\usepackage{amsmath,amsthm,amssymb}
\usepackage{mathrsfs}
\usepackage{bbm}
\usepackage{amsfonts}
\usepackage[margin=0.8in]{geometry}
\usepackage{graphicx}
\usepackage{caption}
\usepackage{subcaption}
\usepackage{wasysym}
\usepackage{enumerate}
\usepackage{algorithm2e}
\usepackage{tabularx}
\usepackage{color}
\usepackage{wrapfig}
\usepackage{todonotes}
\usepackage{mathrsfs}

\newtheorem{thm}{Theorem}[section]
\newtheorem{ass}[thm]{Assumption}
\newtheorem{cor}[thm]{Corollary}
\newtheorem{lem}[thm]{Lemma}
\newtheorem{prop}[thm]{Proposition}
\newtheorem{hyp}{Hypothesis}
\theoremstyle{definition}
\newtheorem{defn}[thm]{Definition}

\theoremstyle{remark}
\newtheorem{rem}[thm]{Remark}

\numberwithin{equation}{section}

\newcommand{\R}{\mathbb R}
\newcommand{\eps}{\varepsilon}

\newcommand{\boldD}{\mathbf D}
\newcommand{\bbD}{\mathbb D}
\newcommand{\bbF}{\mathbb F}
\newcommand{\bbT}{\mathbb T}

\newcommand{\bbN}{\mathbb N}

\newcommand{\mcA}{\mathcal{A}}

\newcommand{\mcB}{\mathcal{B}}
\newcommand{\mcC}{\mathcal C}
\newcommand{\mcD}{\mathcal D}
\newcommand{\mcE}{\mathcal E}

\newcommand{\mcF}{\mathcal F}

\newcommand{\msL}{\mathscr L}

\newcommand{\mcJ}{\mathcal J}
\newcommand{\mcT}{\mathcal T}
\newcommand{\mcP}{\mathcal P}

\newcommand{\mcH}{\mathcal H}

\newcommand{\mcU}{\mathcal U}

\newcommand{\mcR}{\mathcal R}
\newcommand{\mcS}{\mathcal S}

\newcommand{\mcV}{\mathcal V}

\newcommand{\E}{\mathbb{E}}
\newcommand{\Prob}{\mathbb{P}}

\newcommand{\vecb}{\mathbf{b}}

\newcommand{\vecv}{\mathbf{v}}
\newcommand{\esssup}{\mathop{\rm{ess}\,\sup}}
\newcommand{\essinf}{\mathop{\rm{ess}\,\inf}}

\newcommand{\ett}{\mathbbm{1}}

\newcommand{\cadlag}{c\`adl\`ag~}
\newcommand{\cadlagSTOP}{c\`adl\`ag}

\newcommand{\cadlagp}{c\`adl\`ag.~}

\newcommand{\Pred}{\mathcal{P}}
\newcommand{\Prog}{{\rm Prog}}

\newcommand{\ie}{\textit{i.e.\ }}
\newcommand{\eg}{\textit{e.g.\ }}
\newcommand{\etal}{\textit{et.~al.\ }}

\begin{document}

\title{Optimal Stopping of BSDEs with Constrained Jumps and Related Zero-Sum Games\footnote{This work was supported by the Swedish Energy Agency through grant number P2020-90032}}

\author{Magnus Perninge\footnote{M.\ Perninge is with the Department of Physics and Electrical Engineering, Linnaeus University, V\"axj\"o,
Sweden. e-mail: magnus.perninge@lnu.se.}} %
\maketitle
\begin{abstract}
In this paper, we introduce a non-linear Snell envelope which at each time represents the maximal value that can be achieved by stopping a BSDE with constrained jumps. We establish the existence of the Snell envelope by employing a penalization technique and the primary challenge we encounter is demonstrating the regularity of the limit for the scheme. Additionally, we relate the Snell envelope to a finite horizon, zero-sum stochastic differential game, where one player controls a path-dependent stochastic system by invoking impulses, while the opponent is given the opportunity to stop the game prematurely. Importantly, by developing new techniques within the realm of control randomization, we demonstrate that the value of the game exists and is precisely characterized by our non-linear Snell envelope.
\end{abstract}

\section{Introduction}
In recent decades, the optimal stopping problem has garnered considerable attention as one of the fundamental stochastic control problems. As a non-linear counterpart to the classical optimal stopping problem, El Karoui \etal introduced the notion of reflected backward stochastic differential equations (RBSDEs)~\cite{ElKaroui1}. Since their introduction, RBSDEs have found wide-ranging applications in the realm of stochastic control. These applications involve strategies of mixed type, seamlessly integrating stopping~\cite{Bayraktar11,ElAsri2022} (or more generally switching~\cite{HuTang,HamZhang} and impulse control~\cite{Perninge2022}) with classical control. Additionally, RBSDEs have also proven invaluable in addressing related challenges such as stochastic differential games (SDGs)~\cite{HamLep2000}. However, the application of BSDEs, including RBSDEs, is constrained by a notable limitation: their semi-linear nature allow us to relate the solution of a BSDE to a stochastic control problem only when the volatility is not immediately affected by the classical control.

Efforts to address this issue have led to the development of two distinct approaches. On one hand, there is the advancement of quasi-sure analysis~\cite{SonerQuasiSure} and the related concepts of second-order BSDEs (2BSDEs) (see \cite{2bsde07,STZ2bsde}) and $G$-nonlinear expectations~\cite{PengBook19}. On the other hand, there is the consideration of BSDEs driven by both a Brownian motion and an independent Poisson random measure, where the jumps are constrained to exceed a predefined barrier~\cite{Kharroubi2010}. The latter type of BSDEs were related to fully non-linear Hamilton-Jacobi-Bellman integro-partial differential equations (HJB-IPDEs) through a Feynman-Kac representation in \cite{KharroubiPham15}. This innovative approach to stochastic optimal control is commonly referred to as \emph{control randomization}. A significant breakthrough in this field was achieved with the seminal work of \cite{Fuhrman15}, which directly linked the value function of the randomized control problem to that of the original control problem. This eliminated the need for a Feynman-Kac representation, thereby expanding the theoretical framework to encompass stochastic systems with path-dependencies. Building upon this foundation, subsequent advancements extended their approach to the framework of partial information settings in \cite{Bandini18} and optimal switching problems in \cite{Fuhrman2020}.

Whereas approaches based on quasi-sure analysis and related techniques (notably that of \cite{NutzHandel}) have been successfully employed to solve various types of zero-sum stochastic differential games (see \eg \cite{PhamZhang14,NutzZhang15,Possamai20,ImpNL}), the extension of control randomization in this context appears to be constrained. It primarily manifests through a Feynman-Kac relation between RBSDEs, where the jumps are constrained to be non-positive, and fully non-linear variational inequalities that was established in \cite{Choukroun15}. Building upon a result presented in \cite{Bayraktar13}, the latter offers a probabilistic representation of the value function in Markovian controller-and-stopper games. It is worth noting that, the methodology employed in \cite{Choukroun15} to prove the existence of a solution to the RBSDE relies on a double penalization scheme. Therefore, similar to previous studies on doubly reflected BSDEs (see \eg \cite{CvitKar}), their approach assumes strong smoothness conditions on the reflecting barrier.

In the first part of the present work we take an altogether different approach and investigate the non-linear Snell envelope defined as $Y_t:=\esssup_{\tau\in\mcT_t}Y^\tau_t$, where for each $[0,T]$-valued stopping time $\tau$, the quadruple $(Y^\tau,Z^\tau,V^\tau,K^{\tau})$ is the maximal solution to the stopped BSDEs with constrained jumps
\begin{align}\label{ekv:bsde-c-jmp}
  \begin{cases}
    Y^\tau_t=S_\tau+\int_t^\tau f(s,Y^\tau_s,Z^\tau_s,V^\tau_s)ds-\int_t^\tau Z^\tau_s dW_s-\int_t^\tau\!\!\!\int_U V^\tau_s(e)\mu(ds,de)-(K^\tau_\tau-K^\tau_t),\, \forall t\in [0,\tau]\\
    V^\tau_t(e)\geq - \chi(t,Y^\tau_{t-},Z^\tau_t,e),\quad d\Prob\otimes dt\otimes \lambda(de)-a.e.,
  \end{cases}
\end{align}
in which $\chi:[0,T]\times \Omega\times \R^{d+1}\times U\to [0,\infty)$ provides a lower barrier for $V^\tau$. We study a general setting where the barrier, $S$, is only required to be \cadlag and quasi-left upper semi-continuous. Under this assumption, along with mild conditions on the data $f$ and $\chi$, we demonstrate the existence of a \cadlag process $Y$ that satisfies standard integrability assumptions and fulfills the aforementioned relation.


In the second part, we shift our focus to the scenario where the stopped BSDE takes on a linear form,
\begin{align}\label{ekv:bsde-c-jmp-lin}
  \begin{cases}
    Y^\tau_t=\Psi(\tau,X)+\int_t^\tau f(s,X)ds-\int_t^\tau Z^\tau_s dW_s-\int_t^\tau\!\!\!\int_U V^\tau_s(e)\mu(ds,de)-(K^\tau_\tau-K^\tau_t),\, \forall t\in [0,\tau]\\
    V^\tau_t(e)\geq - \chi(t-,X,e),\quad d\Prob\otimes dt\otimes \lambda(de)-a.e.,
  \end{cases}
\end{align}
where we have introduced a state-process, $X$, that solves the path-dependent SDE
\begin{align}\label{ekv:fwd-sde}
X_t&=x_0+\int_0^ta(s,X)ds+\int_0^t\sigma(s,X)dW_s+\int_0^t\!\!\!\int_U\gamma(s-,X,e)\mu(ds,de),\quad\forall t\in [0,T].
\end{align}
The primary contribution in this part lies in establishing a relationship between the non-linear Snell envelope, defined over solutions to \eqref{ekv:bsde-c-jmp-lin}, and a path-dependent SDG of impulse control versus stopping. Specifically, for any $t\in[0,T]$ and any given impulse control $u:=(\eta_j,\beta_j)_{j=1}^\infty\in\mcU_t$ (the set of impulse controls where the first intervention is made after $t$), we let $X^{t,u}$ solve the path-dependent SDE with impulses,
\begin{align}\nonumber
X^{t,u}_s&=x_0+\int_0^s a(r,X^{t,u})dr+\int_0^s\sigma(r,X^{t,u})dW_r+\int_0^{s\wedge t}\!\!\!\int_U\gamma(r-,X^{t,u},e)\mu(dr,de)
\\
&\quad+\sum_{j=1}^N\ett_{[\eta_j\leq s]}\gamma(\eta_j,X^{t,[u]_{j-1}},\beta_j),\quad\forall s\in [0,T],\label{ekv:fwd-sde-2}
\end{align}
where $[u]_k:=(\eta_j,\beta_j)_{j=1}^k$ and $N:=\sup\{j:\eta_j\leq T\}$. We then consider the game of impulse control versus stopping with lower (resp.~upper) value process defined as
\begin{align*}
  \underline Y_t:=\inf_{u^S\in\mcU^{S,W}_t}\sup_{\tau\in\mcT^W_t}J_t(u^{S}(\tau);\tau) \qquad (\text{resp.~}\,\bar Y_t:=\sup_{\tau^{S}\in\mcT^{S,W}_t}\inf_{u\in\mcU^W_t}J_t(u;\tau^{S}(u))).
\end{align*}
Within our problem formulation, the cost/reward functional takes on the form
\begin{align*}
  J_t(u;\tau):=\E\Big[\Psi(\tau,X^{t,u})+\int_t^{\tau}f(s,X^{t,u})ds + \sum_{j=1}^N\ett_{[\eta_j\leq\tau]}\chi(\eta_j,X^{t,[u]_{j-1}},\beta_j)\Big|\mcF_t\Big]
\end{align*}
and $\mcU^{S,W}_t$ (resp.~$\mcT^{S,W}_t$) is the set of non-anticipative maps from the set of stopping times with respect to $\bbF^{t,W}$ (the filtration generated by $\mu(\cdot\cap[0,t],\cdot)$ and $W$) valued in $[t,T]$, denoted $\mcT^W_t$, to the set of $\bbF^{t,W}$-adapted impulse controls with the first intervention after $t$, denoted $\mcU^W_t$ (resp.~$\mcU^W_t\to\mcT^W_t$).

In particular, we demonstrate that, under fairly general assumptions on the involved coefficients, the non-linear Snell envelope $Y_t:=\esssup_{\tau\in\mcT_t}Y^\tau_t$ serves as a representation of the game's value by satisfying $Y_t=\underline Y_t=\bar Y_t$. This finding extends the existing results on path-dependent impulse control in \cite{DjehiceImpulse, Joensson2023}, as well as the recent advancements in path-dependent SDGs involving impulse controls in \cite{Perninge2022,ImpNL}. Notably, our work expands this framework to incorporate scenarios where the opponent employs a stopping rule, while also providing opportunities for the development of more efficient numerical solution methods.
Of greater significance,  however, is that our work bridges a void in the literature on control randomization by extending the applicability of this methodology to incorporate zero-sum SDGs with path-dependencies.
Additionally, our assumptions are formulated in a way to allow the results to transfer to other types of SDGs. In particular, it should be fairly straightforward to adapted the developed methodology to handle controller-stopper games, thereby extending the results in \cite{Choukroun15} to the non-Markovian framework while allowing for a more general setting compared to \cite{BayraktarYao14,NutzZhang15}, as these works are based on a non-degeneracy assumption on the volatility.

The remainder of the article is structured as follows. In the next section, we establish all notation that will be used throughout the first part of the paper and recall some important results on reflected BSDEs with jumps. In Section~\ref{sec:Snell}, we show existence of a non-linear Snell envelope $Y_t=\esssup_{\tau\in\mcT_t}Y^\tau_t$, where $Y^\tau$ is the first component in the unique maximal solution to the general BSDE in~\eqref{ekv:bsde-c-jmp}. In Section~\ref{sec:zsg}, we meticulously outline the framework for our SDG. In addition, we provide preliminary estimates on the solution to the controlled SDE as defined in equation \eqref{ekv:fwd-sde-2} and introduce approximations of the value functions $\underline Y$ and $\bar Y$, based on truncation and discretization. Subsequently, in Section~\ref{sec:game-value} we demonstrate the existence of a value for the game by showing that the upper and lower value functions both coincide with the same non-linear Snell envelope, as defined in Section~\ref{sec:Snell}.

\section{Preliminaries}

\subsection{Probabilistic setup}
Let $(\Omega,\mcF,\Prob)$ be a complete probability space that supports a $d$-dimensional Brownian motion denoted by $W$, and an independent Poisson random measure $\mu$ defined on a compact set $U\subset \mathbb{R}^d$ with a finite compensator $\lambda$. We denote the $\Prob$-augmented natural filtration generated by $W$ and $\mu$ as $\bbF := (\mcF_t)_{t\geq 0}$. In this context, for each $E \in \mcB(U)$, the process $\tilde\mu([0,\cdot],E) := (\mu([0,t],E)-t\lambda(E):t\geq 0)$ is a martingale.

\subsection{Notations}

\noindent Throughout, we will use the following notations, where $T> 0$ is the fixed problem horizon:
\begin{itemize}
  \item For a measure space $(\tilde\Omega,\tilde\mcF)$ and a filtration $\tilde\bbF$ on $\tilde\mcF$ we let $\Prog(\tilde\bbF)$ (resp. $\Pred(\tilde\bbF)$) denote the $\sigma$-algebra of $\tilde\bbF$-progressively (resp. $\tilde\bbF$-predictably) measurable subsets of $\R_+\times \Omega$.
  \item For $p\geq 1$ and a measure space $(E,\mcE,m)$ we let $L^p(E,\mcE,m)$ denote the set of functions $\xi:E\to\R$ which are $\mcE$-measurable and such that $|\xi|^p$ is integrable under $m$. When $m=\Prob$ we often use the shorthand $L^p(E,\mcE)$ and when $(E,\mcE,m)=(\Omega,\mcF,\Prob)$ we sometimes write $L^p$.
  \item We let $L^p_\lambda:=L^p(U,\mcB(U),\lambda)$, \ie the set of all measurable functions $\ell:U\to \R$ such that $\|\ell\|_{L^p_\lambda}:=\big(\int_U|\ell(e)|^p\lambda(de)\big)^{1/p}<\infty$. The set $L^2_\lambda$ is a Hilbert space when equipped with the inner product $\langle \ell_1,\ell_2\rangle_{L^2_\lambda}:=\int_U\ell_1(e)\ell_2(e)\lambda(de)$.
  \item We let $\mcT$ be the set of all $[0,T]$-valued $\bbF$-stopping times and for each $\eta\in\mcT$, we let $\mcT_\eta$ be the corresponding subset of stopping times $\tau$ such that $\tau\geq \eta$, $\Prob$-a.s.
  \item For $p\geq 1$ and $\tau\in\mcT$, we let $\mcS^{p,\tau}$ be the set of all $\R$-valued, $\Prog(\bbF)$-measurable \cadlag processes $(Z_t: t\in[0,\tau])$ such that $\|Z\|_{\mcS^{p,\tau}}:=\E\Big[\sup_{t\in [0,\tau]} |Z_t|^p\Big]^{1/p}<\infty$.
      When $\tau=T$, we use the shorter notation $\mcS^{p}$.
  \item We let $\mcA^{p,\tau}$ be the subset of $\mcS^{p,\tau}$ with all $\mcP(\bbF)$-measurable processes, $Z$, that are non-decreasing and start at $Z_0=0$. Moreover, we let $\mcA^p:=\mcA^{p,T}$.
  \item We let $\mcH^{p,\tau}(W)$ denote the set of all $\R^d$-valued, $\mcP(\bbF)$-measurable processes $(Z_t: t\in[0,\tau])$ such that $\|Z\|_{\mcH^{p,\tau}(W)}:=\E\left[(\int_0^\tau |Z_t|^2 dt)^{p/2}\right]^{1/p}<\infty$. When $\tau=T$, we use the notation $\mcH^{p}(W)$.
  \item We let $\mcH^{p,\tau}(\mu)$ denote the set of all $\R$-valued, $\mcP(\bbF)\otimes \mcB(U)$-measurable mappings $(Z_t(e): t\in[0,\tau],e\in U)$ such that $\|Z\|_{\mcH^{p,\tau}(\mu)}:=\E\Big[\int_0^\tau |Z_t(e)|^p\lambda(de)dt\Big]^{1/p}<\infty$ and set $\mcH^{p}(\mu):=\mcH^{p,T}(\mu)$.
\end{itemize}

Unless otherwise specified, all inequalities involving random variables are assumed to hold $\Prob$-a.s.

\subsection{Reflected BSDEs with Jumps}

Our approach will heavily rely on the existing theory of reflected backward stochastic differential equations (RBSDEs) with jumps. Several studies have addressed the existence and uniqueness of such RBSDEs, with varying assumptions on the involved coefficients and the obstacle, as documented in works such as~\cite{HamOuk2003,Essaky08,HamOuk2016}. We recall the following important result:

\begin{thm}\label{thm:HamOuk-rbsde}\emph{(Hamad\'ene-Ouknine \cite{HamOuk2016})}
  Assume that
  \begin{enumerate}[a)]
  \item $\xi\in L^2(\Omega,\mcF_T)$.
  \item The barrier $S$ is real-valued, $\Prog(\bbF)$-measurable and \cadlag with $S^+\in\mcS^2$ and $S_T\leq\xi$.
  \item $f:[0,T]\times \Omega\times \R^{1+d}\times L^2_\lambda\to\R$ is $\Prog(\bbF)\otimes\mcB(\R^{1+d})\otimes\mcB(L^2_\lambda)$-measurable, such that $\|f(\cdot,0,0,0)\|_{\mcH^2(W)}<\infty$ and for some $k_f>0$ we have that $\Prob$-a.s., for all $(t,y,y',z,z')\in[0,T]\times\R^{2(1+d)}$ and $v,v'\in L^2_\lambda$,
  \begin{align*}
     |f(t,y',z',v')-f(t,y,z,v)|\leq k_f(|y'-y|+|z'-z|+\|v'-v\|_{L^2_\lambda}).
  \end{align*}
  \end{enumerate}
  Then, there exists a unique quadruple $(Y,Z,V,K^+)\in \mcS^2\times\mcH^2(W)\times \mcH^2(\mu)\times\mcA^2$, such that
\begin{align}\label{ekv:HamOuk-rbsde}
  \begin{cases}
    Y_t=\xi+\int_t^T f(s,Y_s,Z_s,V_s)ds-\int_t^T Z_s dW_s-\int_t^T\!\!\!\int_U V_s(e)\tilde\mu(ds,de)+K^+_T-K^+_t,\quad\forall t\in[0,T],\\
    Y_t\geq S_t,\, \forall t\in [0,T],\quad \int_0^T \left(Y_t-S_t\right)dK^c_t=0\mbox{ and }\Delta K^d_t=(S_{t-}-Y_t)^+\ett_{[Y_{t-}=S_{t-}]},
  \end{cases}
\end{align}
where $K^c$ is the continuous and $K^d$ the purely discontinuous part of $K$, respectively.
\end{thm}

The comparison principle for BSDEs with jumps is not as straightforward as for BSDEs driven solely by Brownian motion. Early work on this topic was presented in \cite{Barles97}, where a first result was obtained, and later expanded upon in \cite{Royer06} and \cite{QuenSul13}. In addition to the prerequisites for Theorem~\ref{thm:HamOuk-rbsde}, these studies assumed an integral constraint on the driver. In a related context, \cite{QuenSul14} employed a similar comparison result to establish a connection between the solution of a reflected BSDE with jumps and a stopping problem. In light of these findings, we recall the following:

\begin{thm}\label{thm:Quenez-Sulem-rbsde}\emph{(Quenez-Sulem \cite{QuenSul13,QuenSul14})}
Assume that $(\xi,S,f)$ satisfies the assumptions in Theorem~\ref{thm:HamOuk-rbsde} and that $d\Prob\otimes dt$-a.e.~for all $(y,z)\in \times\R^{1+d}$ and $v,v'\in L^2_\lambda$ we have
\begin{align*}
     f(t,y,z,v')-f(t,y,z,v)\geq \langle\theta_t^{y,z,v',v}, v'-v\rangle_{L^2_\lambda},
\end{align*}
where $\theta:[0,T]\times\Omega\times \R^{1+d}\times (L_\lambda^2)^2\to L^2_\lambda$ is $\mcP(\bbF)\otimes\mcB(\R^{d+1})\otimes\mcB((L_\lambda^2)^2)$-measurable, bounded and satisfies $d\Prob\otimes dt\otimes d\lambda$-a.e.
\begin{align*}
     \theta_t^{y,z,v',v}(e)\geq-1 \quad \text{and}\quad |\theta_t^{y,z,v',v}(e)|\leq \psi(e),
\end{align*}
with $\psi\in L^2_\lambda$. We then have:
\begin{enumerate}[i)]
  \item The unique solution to RBSDE \eqref{ekv:HamOuk-rbsde} satisfies $Y_t=\esssup_{\tau\in\mcT_t}Y^{\tau}_t$, where the triple $(Y^{\tau},Z^{\tau},V^\tau)\in \mcS^{2,\tau}\times\mcH^{2,\tau}(W)\times\mcH^{2,\tau}(\mu)$ is the unique solution to
\begin{equation}\label{ekv:Quenez-stopped-bsde}
      Y^\tau_t=\xi\ett_{[\tau=T]}+S_{\tau}\ett_{[\tau<T]}+\int_t^\tau f(s,Y^\tau_s,Z^\tau_s,V^\tau_s)ds-\int_t^\tau Z^\tau_s dW_s-\int_t^\tau\!\!\!\int_U V^\tau_s(e)\tilde\mu(ds,de),\quad\forall t\in [0,\tau].
\end{equation}
  \item If $S$ is quasi-left upper semi-continuous, then with $D_t:=\inf\{s\geq t: Y_s=S_s\}\wedge T$ we have the representation $Y_t=Y^{D_t}_t$ and get that $K^+$ is continuous and satisfies $K^+_{D_t}-K^+_t=0$, $\Prob$-a.s.
  \item Assume that the parameters of another RBSDE, $(\tilde \xi,\tilde S,\tilde f)$, satisfy the requirements of Theorem~\ref{thm:HamOuk-rbsde} in addition to $\tilde \xi\leq \xi$, $\Prob$-a.s., $\tilde S_t\leq S_t$, $\Prob$-a.s., for all $t\in[0,T]$ and that $\tilde f(t,y,z,v)\leq f(t,y,z,v)$, for all $(z,y,v)\in \R^{d+1}\times L^2_\lambda$, $dt\otimes d\Prob$-a.e., then the solution $(\tilde Y^\tau,\tilde Z^\tau,\tilde V^\tau)$ of \eqref{ekv:Quenez-stopped-bsde} with parameters $(\tilde \xi,\tilde S,\tilde f)$ satisfies $\tilde Y^\tau_t\leq Y^\tau_t$, $\Prob$-a.s.~for each $t\in[0,T]$ and $\tau\in\mcT_t$. In particular, if $\tilde Y$ is the first component in the solution to \eqref{ekv:HamOuk-rbsde} with parameters $(\tilde \xi,\tilde S,\tilde f)$ we get that $\tilde Y_t\leq Y_t$, $\Prob$-a.s.
\end{enumerate}
\end{thm}

\bigskip

\begin{rem}
Recall here the concept of quasi-left continuity: A \cadlag process $(X_t:t\geq 0)$ is quasi-left continuous if for each predictable stopping time $\theta$ and every announcing sequence of stopping times $\theta_k\nearrow\theta$ we have $X_{\theta -}:=\lim\limits_{k\to\infty}X_{\theta_k} = X_\theta$, $\Prob$-a.s. Similarly, $X$ is quasi-left upper semi-continuous if $X_{\theta -}\leq  X_\theta$, $\Prob$-a.s.
\end{rem}


\section{Optimal stopping of BSDEs with constrained jumps\label{sec:Snell}}
In this section, we consider optimal stopping of BSDEs with constrained jumps. For each $\tau\in\mcT$, we recall the definition of the quadruple $(Y^\tau, Z^\tau, V^\tau, K^{-,\tau})\in\mcS^{2,\tau}\times\mcH^{2,\tau}(W)\times\mcH^{2,\tau}(\mu)\times\mcA^{2,\tau}$ as the maximal\footnote{Maximal in the sense that $Y^{\tau}_t\geq \tilde Y^{\tau}_t$ whenever $\tilde Y^{\tau}$ is the first component of another solution} solution to the following equation:
\begin{align}\label{ekv:bsde-c-jmp-rec}
  \begin{cases}
    Y^{\tau}_t=S_\tau+\int_t^\tau f(s,Y^{\tau}_s,Z^{\tau}_s,V^\tau_s)ds-\int_t^\tau Z^{\tau}_s dW_s-\int_t^\tau\!\!\!\int_U V^{\tau}_s(e)\mu(ds,de)-(K^{-,\tau}_\tau-K^{-,\tau}_t),\quad\forall t\in[0,\tau]
    \\
    V^{\tau}_s(e)\geq-\chi(s,Y^{\tau}_{s-},Z^{\tau}_s,e),\,d\Prob\otimes\lambda(de)\otimes ds-a.e.,
  \end{cases}
\end{align}
where the data $(f,S,\chi)$ satisfies the assumptions in Theorem~\ref{thm:nl-Snell} below. Indeed, repeating the steps in~\cite{Kharroubi2010} we find that, under these assumptions, \eqref{ekv:bsde-c-jmp-rec} admits a unique maximal solution for each $\tau\in\mcT$. The main contribution of the present section is that we show the existence of an aggregator $Y\in\mcS^2$ satisfying $Y_\eta=\esssup_{\tau\in\mcT_\eta} Y^{\tau}_\eta$ for every $\eta\in\mcT$, in addition to a corresponding optimal stopping time. We summarize this result in the following theorem:
\begin{thm}\label{thm:nl-Snell}
Assume that,
\begin{itemize}
  \item $S\in\mcS^2$ is left upper semi-continuous at predictable stopping times, in particular $\xi:=S_T\in L^2(\Omega,\mcF_T)$ and $\lim_{t\to T}S_t\leq\xi$, $\Prob$-a.s.;
  \item $f:[0,T]\times\Omega\times\R^{1+d}\times L^2_\lambda\to\R$ is $\Prog(\bbF)\otimes\mcB(\R^{1+d})\otimes\mcB(L^2_\lambda)$-measurable, Lipschitz continuous in $(z,y,v)$, \ie
   \begin{align*}
     |f(t,y,z,v)-f(t,y',z',v')|\leq k_f(|y-y'|+|z-z'|+\|v-v'\|_{L^2_\lambda})
   \end{align*}
   for all $(t,y,y',z,z')\in [0,T]\times R^{2(1+d)}$, $\Prob$-a.s.,~and such that $\|f(\cdot,0,0,0)\|_{\mcH^2(W)}<\infty$. Moreover, $d\Prob\otimes dt$-a.e.~for all $(y,z)\in \times\R^{1+d}$ and $v,v'\in L^2_\lambda$ we have
\begin{align*}
     f(t,y,z,v')-f(t,y,z,v)\geq \langle \tilde\theta_t^{y,z,v',v}, v'-v\rangle_{L^2_\lambda}
\end{align*}
where $\tilde\theta:[0,T]\times\Omega\times \R^{1+d}\times (L_\lambda^2)^2\to L^2_\lambda$ is $\mcP(\bbF)\otimes\mcB(\R^{1+d})\otimes\mcB((L_\lambda^2)^2)$-measurable, bounded and satisfies $d\Prob\otimes dt\otimes d\lambda$-a.e.
\begin{align*}
     \tilde\theta_t^{y,z,v',v}(e)\geq 0 \quad \text{and}\quad |\tilde\theta_t^{y,z,v',v}(e)|\leq \psi(e),
\end{align*}
with $\psi\in L^2_\lambda$; and
  \item $\chi:[0,T]\times\Omega\times\R^{1+d}\times U\to[0,\infty)$ is $\Pred(\bbF)\otimes\mcB(\R^{1+d})\otimes\mcB(U)$-measurable, Lipschitz continuous in $(z,y)$, \ie
   \begin{align*}
     |\chi(t,y,z,e)-\chi(t,y',z',e)|\leq k_\chi(|y-y'|+|z-z'|)
   \end{align*}
   for all $(t,y,y',z,z')\in [0,T]\times R^{2(1+d)}$, $d\Prob\otimes \lambda(de)$-a.e.,~and such that $\chi(\cdot,0,0,\cdot)\in\mcH^2(\mu)$.
\end{itemize}
Then there exists a $Y\in\mcS^2$ such that for every $\eta\in\mcT$, we have $Y_\eta=\esssup_{\tau\in\mcT_\eta}Y^\tau_\eta=Y^{\tau^*}_\eta$, where $\tau^*:=\inf\{t\geq\eta:Y_t=S_t\}\in\mcT_\eta$.
\end{thm}

As noted in the introduction, a similar results was shown in \cite{Choukroun15}. In particular, assuming that $\chi$ is identically zero and introducing an additional condition of regularity on the barrier $S$, the work presented in \cite{Choukroun15} established the existence of a process $Y$ satisfying the conditions of Theorem~\ref{thm:nl-Snell}. Furthermore, it was shown that there are processes $(Z, V, K^{-},K^+)\in\mcH^{2}(W)\times\mcH^{2}(\mu)\times\mcA^{2}\times\mcA^{2}$ such that the quintet $(Y, Z, V, K^{-},K^+)$ is the unique maximal solution to the reflected BSDE:
\begin{align}\label{ekv:bsde-c-jmp-rec-2}
  \begin{cases}
    Y_t=S_\tau+\int_t^\tau f(s,Y_s,Z_s,V_s)ds-\int_t^\tau Z_s dW_s-\int_t^\tau\!\!\!\int_U V_s(e)\mu(ds,de)-(K^{-}_\tau-K^{-}_t)+K^{+}_\tau-K^{+}_t,
    \\
    Y_t\geq S_t,\, \forall t\in [0,T],\quad \int_{0}^T \big(Y_t-S_t\big)dK^{+}_t=0
    \\
    V_s(e)\geq 0,\,d\Prob\otimes\lambda(de)\otimes ds-a.e.
  \end{cases}
\end{align}
However, the specific regularity assumption on $S$ made in \cite{Choukroun15} and their choice of $\chi\equiv 0$ render their methodology unsuitable for addressing impulse control problems. Conversely, Theorem~\ref{thm:nl-Snell} is meticulously tailored to accommodate impulses and the corresponding intervention costs, making it a more appropriate framework for handling this particular application.

To prove Theorem~\ref{thm:nl-Snell} we apply an approximation routine based on penalization and for each $n\in\bbN$, we let $(Y^{n},Z^{n},V^{n},K^{+,n})\in\mcS^2\times\mcH^2(W)\times\mcH^2(\mu)\times\mcA^2$ be the unique solution to
\begin{align}\label{ekv:rbsde-pen-n}
  \begin{cases}
    Y^{n}_t=\xi+\int_t^T f(s,Y^{n}_s,Z^{n}_s,V^n_s)ds-\int_t^T Z^{n}_s dW_s-\int_t^T\!\!\!\int_U V^{n}_s(e)\mu(ds,de)+K^{+,n}_T-K^{+,n}_t
    \\
    \quad-n\int_t^T\!\!\!\int_U(V^{n}_s(e)+\chi(s,Y^{n}_{s},Z^{n}_s,e))^-\lambda(de)ds,\quad \forall t\in[0,T]\\
    Y^{n}_t\geq S_t,\, \forall t\in [0,T],\quad \int_{0}^T \big(Y^{n}_t-S_t\big)dK^{+,n}_t=0
  \end{cases}
\end{align}
and set $K^{-,n}_t:=n\int_0^t\!\!\int_U(V^{n}_s(e)+\chi(s,Y^{n}_{s},Z^{n}_s,e))^-\lambda(de)ds$.

\begin{rem}\label{rem:rbsde-pen-n-has-comparison}
Note that for any $(t,y,z,v,v')\in[0,T]\times\R^{1+d}\times (L^2_\lambda)^2$, we have
\begin{align*}
  -n\int_U(v(e)+\chi(s,y,z,e))^-\lambda(de)+n\int_U(v'(e)+\chi(s,y,z,e))^-\lambda(de)\geq n\int_U\ett_{[v(e)> v'(e)]}(v'(e)-v(e))\lambda(de).
\end{align*}
Since $(Y^{n},Z^{n},V^{n},K^{+,n})$ satisfies
\begin{align*}
  Y^{n}_t=\xi+\int_t^T \tilde f^n(s,Y^{n}_s,Z^{n}_s,V^n_s)ds-\int_t^T Z^{n}_s dW_s-\int_t^T\!\!\!\int_U V^{n}_s(e)\tilde\mu(ds,de)+K^{+,n}_T-K^{+,n}_t,\quad \forall t\in[0,T]
\end{align*}
with
\begin{align*}
  \tilde f^n(t,y,z,v):=f(t,y,z,v)-\int_Uv(e)\lambda(de)-n\int_U(v(e)+\chi(s,y,z,e))^-\lambda(de),
\end{align*}
it thus follows that \eqref{ekv:rbsde-pen-n} satisfies the conditions of Theorem~\ref{thm:Quenez-Sulem-rbsde} with $\theta^{y,z,v,v'}_t(e)=\tilde\theta^{y,z,v,v'}_t(e)-1+n\ett_{[v(e)> v'(e)]}$.
\end{rem}

\begin{lem}\label{lem:bound-Yn}
There is a $C>0$ such that
\begin{align*}
  \|Y^n\|_{\mcS^2}+\|Z^n\ett_{[\eta,\tau_n]}\|_{\mcH^2(W)}+\|V^n\ett_{[\eta,\tau_n]}\|_{\mcH^2(\mu)}+\|K^{-,n}_{\tau_n}-K^{-,n}_{\eta}\|_{L^2}\leq C
\end{align*}
for all $\eta\in\mcT$ and $n\in\bbN$, with $\tau_n:=\inf\{s\geq\eta:Y^n_s=S_s\}\wedge T\in\mcT_\eta$.
\end{lem}

\noindent\emph{Proof.} First note that $S\leq Y^n\leq Y^0$, ensuring the existence of a $C>0$ such that $\|Y^n\|_{\mcS^2}\leq C$ for all $n\in\bbN$. Next, Ramark~\ref{rem:rbsde-pen-n-has-comparison} and Theorem~\ref{thm:Quenez-Sulem-rbsde}.(ii) implies that for each $\eta\in\mcT$ and $n\in\bbN$, the corresponding stopping time $\tau_n$ is optimal for \eqref{ekv:rbsde-pen-n} in the sense that
\begin{align*}
    Y^{n}_t=S_{\tau_n}+\int_t^{\tau_n} f(s,Y^{n}_s,Z^{n}_s,V^n_s)ds-\int_t^{\tau_n} Z^{n}_s dW_s-\int_t^{\tau_n}\!\!\!\int_U V^{n}_s(e)\mu(ds,de)
    \\
    \quad-n\int_t^{\tau_n}\!\!\!\int_U(V^{n}_s(e)+\chi(s,Y^{n}_s,Z^{n}_s,e))^-\lambda(de)ds,\quad \forall t\in[\eta,\tau_n].
\end{align*}
Consequently, It{\^o}'s formula applied to $|Y^n|^2$ gives
\begin{align*}
 |S_{\tau_n}|^2&= |Y^n_t|^2 - 2\int_t^{\tau_n} Y^n_s f(s,Y^n_s,Z^n_s,V^n_s)ds + 2\int_t^{\tau_n}Y^n_sZ^n_sdW_s + 2\int_t^{\tau_n} Y^n_sdK^{-,n}_s + \int_t^{\tau_n}|Z^n_s|^2ds
  \\
  &\quad+\int_t^{\tau_n}\!\!\!\int_U (2Y^n_{s-}V^n_{s}(e)+|V^n_{s}(e)|^2)\mu(ds,de),\quad \forall t\in[\eta,\tau_n].
\end{align*}
Rearranging, taking expectation and using the inequality $2ab\leq a^2+b^2$ we find that\footnote{Throughout, $C$ will denote a generic positive constant that may change value from line to line.}
\begin{align*}
 &\E\Big[\int_\eta^{\tau_n}|Z^n_s|^2ds\Big]+\E\Big[\int_\eta^{\tau_n}\!\!\!\int_U |V^n_{s}(e)|^2\lambda(de)ds\Big]
 \\
 &\leq \E\Big[|S_{\tau_n}|^2 +2\int_\eta^{\tau_n} Y^n_s f(s,Y^n_s,Z^n_s,V^n_s)ds -  2\int_\eta^{\tau_n} Y_sdK^{-,n}_s +
\int_\eta^{\tau_n}\!\!\!\int_U 2Y^n_{s}V^n_s(e)\lambda(de)ds\Big]
 \\
 &\leq \E\Big[|S_{\tau_n}|^2+C(1+\frac{1}{\alpha})\sup_{s\in [\eta,\tau_n]}|Y^n_s|^2+\frac{1}{2}\int_\eta^{\tau_n} |f(s,0,0,0)|ds+\frac{1}{2}\int_\eta^{\tau_n} |Z^n_s|^2ds
 \\
 &\quad+\frac{1}{2}\int_\eta^{\tau_n}\!\!\!\int_U |V^n_s(e)|^2\lambda(de)ds+\frac{\alpha}{2} |K^{-,n}_{\tau_n}- K^{-,n}_{\eta}|^2 \Big]
\end{align*}
for any $\alpha>0$ and we conclude that
\begin{align*}
 & \E\Big[\int_\eta^{\tau_n}|Z^n_s|^2ds\Big]+\E\Big[\int_\eta^{\tau_n}\!\!\!\int_U |V^n_{s}(e)|^2\lambda(de)ds\Big]
\leq C(1+\frac{1}{\alpha})+ \alpha\E\big[ |K^{-,n}_{\tau_n}- K^{-,n}_{\eta}|^2 \big].
\end{align*}
On the other hand, squaring
\begin{align*}
    K^{-,n}_{\tau_n}- K^{-,n}_{\eta} = S_{\tau_n}-Y^n_\eta+\int_\eta^{\tau_n} f(s,Y^n_s,Z^n_s,V^n_s)ds-\int_\eta^{\tau_n} Z^n_s dW_s-\int_\eta^{\tau_n}\!\!\!\int_U V^n_{s-}(e)\mu(ds,de)
\end{align*}
gives
\begin{align*}
    \E\big[|K^{-,n}_{\tau_n}- K^{-,n}_{\eta}|^2\big]\leq C(1+\E\Big[\int_\eta^{\tau_n}|Z^n_s|^2ds+\int_\eta^{\tau_n}\!\!\!\int_U |V^n_{s}(e)|^2\lambda(de)ds\Big])
\end{align*}
and by choosing $\alpha>0$ sufficiently small the assertion follows.\qed\\

We introduce the set of admissible densities, denoted by $\mcV$, as the collection of all $\Pred(\bbF)\otimes\mcB(U)$-measurable, bounded maps $\nu=\nu_t(\omega,e):[0,T]\times\Omega\times U\to [0,\infty)$. For each $\nu\in\mcV$ and $\tau\in\mcT$, we let $(Y^{\tau;\nu},Z^{\tau;\nu},V^{\tau;\nu})\in\mcS^{2,\tau}\times\mcH^{2,\tau}(W)\times\mcH^{2,\tau}(\mu)$ be the unique solution to the BSDE
\begin{align}
    Y^{\tau;\nu}_t&=S_{\tau}+\int_t^\tau f^\nu(s,Y^{\tau;\nu}_s,Z^{\tau;\nu}_s,V^{\tau;\nu}_s)ds-\int_t^\tau Z^{\tau;\nu}_s dW_s
 -\int_t^\tau\!\!\!\int_U V^{\tau;\nu}_s(e)\mu(ds,de),\quad \forall t\in[0,\tau],\label{ekv:reward-BSDE}
\end{align}
where
\begin{align*}
f^\nu(t,y,z,v):=f(t,y,z,v)+\int_U(v(e)+\chi(t,y,z,e))\nu_t(e)\lambda(de)
\end{align*}

The following lemma gives a representation of $Y^n$ as the value of a zero-sum game:

\begin{lem}\label{lem:Yn-repr}
For any $\eta\in\mcT$, we have
\begin{align}\label{ekv:Yn-repr}
  Y^{n}_\eta=\esssup_{\tau\in\mcT_\eta}\essinf_{\nu\in\mcV^n}Y^{\tau;\nu}_\eta = \essinf_{\nu\in\mcV^{n}}\esssup_{\tau\in\mcT_\eta}Y^{\tau;\nu}_\eta,
\end{align}
where $\mcV^n:=\{\nu\in\mcV: \nu_t(\omega,e)\in [0,n],\,\forall(t,\omega,e)\in [0,T]\times\Omega\times U\}$.
\end{lem}

\noindent\emph{Proof.} Letting $\nu^{*,n}_s(e):=n\ett_{[V_s^{n}(e)<-\chi(s,Y^{n}_{s-},Z^{n}_s,e)]}$, we note that $\nu^{*,n}\in\mcV^n$ and \begin{align*}
f^{\nu^{*,n}}(t,Y^{n}_{t},Z^{n}_t,V^n_t) = f(t,Y^{n}_{t},Z^{n}_t,V^n_t)-n\int_U(V^n(e)+\chi(t,Y^{n}_{t},Z^{n}_t,e))^-\lambda(de),
\end{align*}
$dt\otimes d\Prob$-a.e., giving, by Theorem~\ref{thm:Quenez-Sulem-rbsde}.(i-ii), that $Y^{n}_\eta=\esssup_{\tau\in\mcT_\eta}Y^{\tau;\nu^{*,n}}_\eta=Y^{\tau_n;\nu^{*,n}}_\eta$. Now, pick $\nu\in\mcV^n$ and note that since
\begin{align*}
  -n\int_U(v(e)+\chi(t,y,z,e))^-\lambda(de)\leq \int_U(v(e)+\chi(t,y,z,e))\nu_t(e)\lambda(de)
\end{align*}
pointwisely, we can argue as in Remark~\ref{rem:rbsde-pen-n-has-comparison} to conclude that $Y^{\tau;\nu^{*,n}}_\eta\leq Y^{\tau;\nu}_\eta$ for any $\tau\in\mcT_\eta$. Combined, this gives that for all $(\tau,\nu)\in \mcT_\eta\times \mcV^n$, we have
\begin{align*}
  Y^{\tau;\nu^{*,n}}_\eta\leq Y^{\tau_n;\nu^{*,n}}_\eta\leq Y^{\tau_n;\nu}_\eta.
\end{align*}
Hence, there exists a saddle-point for the game and consequently the representation \eqref{ekv:Yn-repr} holds.\qed\\

An immediate consequence of Lemma~\ref{lem:Yn-repr} is that the sequence $(Y^n)_{n\in\bbN}$ is non-increasing. As it is, moreover, bounded from below by the barrier $S$, we find that there is a progressively measurable process $\tilde Y$ such that $Y^{n}\searrow \tilde Y$, pointwisely. Since each $\nu\in\mcV$ is bounded and thus belongs to $\mcV^n$ for some $n$, Lemma~\ref{lem:Yn-repr} gives that
\begin{align*}
  \tilde Y_\eta = \essinf_{\nu\in\mcV}\esssup_{\tau\in\mcT_\eta}Y^{\tau;\nu}_\eta
\end{align*}
for any stopping time $\eta\in\mcT$. Using this relation we are able to prove the following lemma:

\begin{lem}
The process $\tilde Y\in\mcS^2$. In particular, $\tilde Y$ is \cadlagp
\end{lem}

\noindent\emph{Proof.} Since $Y^{n}$ is a sequence of optional processes, $\tilde Y$ is optional. To prove right-continuity, we thus only need to show that $\tilde Y$ is right-continuous at each stopping time. Moreover, since $\tilde Y$ is the limit of a non-increasing sequence of \cadlag processes it is clearly right upper semi-continuous. We argue by contradiction and assume that there is an $\vartheta\in\mcT$ and an $\eps>0$ such that $\liminf_{t\searrow\vartheta}\tilde Y_{t}\leq \tilde Y_{\vartheta}-\eps$ on some set $B\in\mcF_\vartheta$ of positive measure. If this would be true, then the sequence of stopping times defined as $\vartheta_n:=\inf\{s>\vartheta: Y^{n}_{s}\leq \tilde Y_{\vartheta}-\eps/2\}\wedge T$ would tend to $\vartheta$ on $B$. Now,
\begin{align*}
\tilde Y_{\vartheta}-\tilde Y_{\vartheta_i}&=\essinf_{\nu\in\mcV}\esssup_{\tau\in\mcT_\vartheta}Y^{\tau;\nu}_\vartheta - \essinf_{\nu\in\mcV}\esssup_{\tau\in\mcT_{\vartheta_i}}Y^{\tau;\nu}_{\vartheta_i}
\\
&\leq \esssup_{\nu\in\bar\mcV_{\vartheta_i},\,\tau\in\mcT_{\vartheta}}(Y^{\tau;\nu}_{\vartheta}-Y^{\tau\vee\vartheta_i;\nu}_{\vartheta_i}),
\end{align*}
where
\begin{align*}
  \bar \mcV_t:=\{\nu\in\mcV:\nu_s\equiv 1,\,\forall s\in [0,t]\} \cap \{\nu\in\mcV:\esssup_{\tau\in\mcT_t}Y^{\nu,\tau}_t\leq Y^0_t\}.
\end{align*}
On the other hand, for each $\nu\in\bar\mcV_{\vartheta_i}$ and $\tau\in\mcT_{\vartheta}$, we have
\begin{align*}
  Y^{\tau;\nu}_{\vartheta} &= \ett_{[\tau\leq \vartheta_i]}S_\tau + \ett_{[\tau> \vartheta_i]}Y^{\tau;\nu}_{\vartheta_i}+\int_\vartheta^{\vartheta_i\wedge\tau} \big(f(s,Y^{\tau;\nu}_s,Z^{\tau;\nu}_s,V^{\tau;\nu}_s)+\chi(s,Y^{\tau;\nu}_s,Z^{\tau;\nu}_s)\big)ds
  \\
  &\quad-\int_\vartheta^{\vartheta_i\wedge\tau} Z^{\tau;\nu}_s dW_s -\int_{\vartheta+}^{\vartheta_i\wedge\tau}\!\!\int_U V^{\tau;\nu}_s(e)\tilde\mu(ds,de)
\end{align*}
and we find that
\begin{align*}
  \ett_B(Y^{\tau;\nu}_{\vartheta}- Y^{\tau\vee\vartheta_i;\nu}_{\vartheta_i}) &\leq \ett_B\Big(S_{\tau\wedge\vartheta_i}-S_{\vartheta_i}+ C\int_\vartheta^{\vartheta_i}\big(|f(s,0,0,0)|+\chi_s+|Y^{\tau;\nu}_s|+|Z^{\tau;\nu}_s|+\|V^{\tau;\nu}_s\|_{L^2_\lambda}\big)ds
  \\
  &\quad+\int_\vartheta^{\vartheta_i\wedge\tau} Z^{\tau;\nu}_s dW_s +\int_{\vartheta+}^{\vartheta_i\wedge\tau}\!\!\int_U V^{\tau;\nu}_s(e)\tilde\mu(ds,de)\Big),
\end{align*}
where $\chi:=\int_U\chi(\cdot,0,0,e)\lambda(de)$. Taking the expectation on both sides and utilizing the fact that $B\in\mcF_\vartheta$, we can deduce that
\begin{align}\nonumber
  &\E\big[\ett_B\esssup_{\nu\in\bar\mcV_{\vartheta_i},\,\tau\in\mcT_{\vartheta}}(Y^{\tau;\nu}_{\vartheta}-Y^{\tau\vee\vartheta_i;\nu}_{\vartheta_i})\big]
  \\
  &\leq \sup_{\nu\in\bar\mcV_{\vartheta_i},\,\tau\in\mcT_{\vartheta}}\E\Big[\ett_B\big(S_{\tau\wedge\vartheta_i}-S_{\vartheta_i}+ C\int_\vartheta^{\vartheta_i}(|f(s,0,0,0)|+\chi_s+|Y^{\tau;\nu}_s|+|Z^{\tau;\nu}_s|+\|V^{\tau;\nu}_s\|_{L^2_\lambda})ds\Big].\label{ekv:disc-bound}
\end{align}
Concerning the first term, right-continuity of $S$ and dominated convergence implies that $\ett_B\big(S_{\tau\wedge\vartheta_i}-S_{\vartheta_i})$ tends to 0 in $L^1$ as $i\to\infty$, uniformly over all $\tau\in\mcT_\vartheta$. Moreover, $|Y^{\tau;\nu}|\leq |Y^0| + |S|$, where the right-hand side belongs to $\mcS^2$, and a simple application of Jensen's inequality together with the fact that $f(\cdot,0,0,0)$, $\chi$ both have finite $\mcH^2(W)$-norms, and since there is a $C>0$ such that
\begin{align*}
  &\sup_{\nu\in\bar\mcV_{\vartheta_i},\,\tau\in\mcT_{\vartheta}}\E\Big[\int_\vartheta^{\vartheta_i}|Z^{\tau;\nu}_s|^2ds + \int_\vartheta^{\vartheta_i}\!\!\!\int_U|V^{\tau;\nu}_s(e)|^2\lambda(de)ds\Big]\leq C
\end{align*}
for all $i\in\bbN$, we find that the right-hand side of \eqref{ekv:disc-bound} tends to 0 as $i\to\infty$. This is a contradiction, since our construction above implies that
\begin{align*}
  \liminf_{i\to\infty}\E\big[\ett_B(\tilde Y_{\vartheta}-\tilde Y_{\vartheta_i})\big]\geq \Prob[B]\eps/2>0.
\end{align*}
In particular, we conclude that $\tilde Y$ is right-continuous.

We turn to the existence of left limits and consider the number of downcrossings of the interval $[a,b]$ made by $\tilde Y$ for any two constants $a,b\in\R$ with $a<b$. We thus define $\vartheta^b_i:=\inf\{s>\vartheta^a_{i-1}:\tilde Y_s\geq b\}\wedge T$, with $\vartheta^a_0:=0$ and $\vartheta^a_i:=\inf\{s>\vartheta^b_{i}:\tilde Y_s\leq a\}\wedge T$ for $i=1,\ldots$. Then, $\vartheta^b_1,\vartheta^a_1,\vartheta^b_2,\ldots$ is a non-decreasing sequence of stopping times (strictly until it hits $T$, since $\tilde Y$ is right-continuous) and $\vartheta^a_i\nearrow \vartheta^{a,b}$ (and then also $\vartheta^b_i\nearrow \vartheta^{a,b}$) for some predictable stopping time $\vartheta^{a,b}\in\mcT$. Moreover, on the set $B_i:=\{\omega:\vartheta_i^a<T\}$ we have $\tilde Y_{\vartheta^b_i}-\tilde Y_{\vartheta^a_i}\geq b-a$. On the other hand, arguing as above we find that
\begin{align*}
  \E\big[\tilde Y_{\vartheta^b_i}- \tilde Y_{\vartheta^a_i}\big] &\leq
  \sup_{\nu\in\bar\mcV_{\vartheta^a_i},\,\tau\in\mcT_{\vartheta^b_i}}\E\Big[S_{\tau\wedge \vartheta^a_i}-S_{\vartheta^a_i}+ C\int_{\vartheta^b_i}^{\vartheta^a_i}(|f(s,0,0,0)|+\chi_s+|Y^{\tau;\nu}_s|+|Z^{\tau;\nu}_s|+\|V^{\tau;\nu}_s\|_{L^2_\lambda})ds\Big],
\end{align*}
where the right-hand side tends to 0 since $S$ has left limits and $\vartheta^a_i-\vartheta^b_i$, tends to 0, $\Prob$-a.s.  Now, by construction
\begin{align*}
  \tilde Y_{\vartheta^b_i}- \tilde Y_{\vartheta^a_i}\geq \ett_{B_i}(b-a)+\ett_{[\vartheta^b_i<T=\vartheta^a_i]}(\tilde Y_{\vartheta^b_i}-\xi)
\end{align*}
and thus
\begin{align}\nonumber
  \ett_{B_i}(b-a) &\leq \tilde Y_{\vartheta^b_i}- \tilde Y_{\vartheta^a_i} +\ett_{[\vartheta^b_i<T=\vartheta^a_i]}(\xi-\tilde Y_{\vartheta^b_i})
  \\
  &\leq \tilde Y_{\vartheta^b_i}- \tilde Y_{\vartheta^a_i} +\ett_{B_{i-1}\setminus B_i}(|\xi|+|\tilde Y_{\vartheta^b_i}|).\label{ekv:Bi-bound}
\end{align}
Since $B_i\subset B_{i-1}$, letting $B:=\cap_{i=1}^\infty B_i$ gives that $B_{i-1}\setminus B_{i}\subset B_{i-1}\setminus B$ and thus $\limsup_{i\to\infty}\Prob[B_{i-1}\setminus B_{i}]=0$ and by taking the expectation of both sides in \eqref{ekv:Bi-bound}, we conclude that $\Prob[B]=0$. As the number of downcrossings of the interval $[a,b]$, denoted $D([a,b])$, is finite on the set $\Omega\setminus B$ and since $a<b$ were arbitrary, we conclude that $D([a,b])$ is $\Prob$-a.s.~finite for any $a,b\in\R$ with $a<b$. By countable additivity this holds simultaneously for all rational pairs $a<b$ (outside of a $\Prob$-null set) and thus also simultaneously for all real pairs $a<b$ proving that $\tilde Y$ has left limits everywhere.

Combining Lemma \ref{lem:bound-Yn} and Fatou's lemma gives that $\|\tilde Y\|_{\mcS^2}<\infty$ and we conclude that $\tilde Y\in\mcS^2$.\qed\\

\begin{rem}
By construction, it also follows that the limit $\tilde Y$ is quasi-left upper semi-continuous. In particular, if $(\vartheta_i)_{i\in\bbN}$ is an increasing sequence in $\mcT$ with $\vartheta_i\nearrow\vartheta$, then $\lim_{i\to\infty}\tilde Y_{\vartheta_i}\leq \lim_{i\to\infty} Y^n_{\vartheta_i}\leq Y^n_{\vartheta}$ for any $n\in\bbN$, whereas $Y^n_{\vartheta}\searrow \tilde Y_{\vartheta}$ as $n\to\infty$ by definition.
\end{rem}

As $Y^{n}$ is non-increasing in $n$, the sequence of stopping times $(\tau_n)_{n\in\bbN}$ is non-increasing and by right-continuity of the filtration, we conclude that there is a $\tau^\diamond\in\mcT_\eta$ such that $\tau_n\searrow \tau^\diamond$. Moreover, right-continuity of $\tilde Y$ and $S$ implies that
\begin{align*}
  \tilde Y_{\tau^\diamond}=\lim_{n\to\infty}\tilde Y_{\tau_n}\leq \limsup_{n\to\infty}\tilde Y^n_{\tau_n}=\lim_{n\to\infty}S_{\tau_n}=S_{\tau^\diamond}.
\end{align*}
On the other hand, $Y^n_{\tau_\diamond}\geq S_{\tau^\diamond}$ for each $n\in\bbN$ from which it follows that $\tilde Y_{\tau^\diamond}=\lim_{n\to\infty}\tilde Y^n_{\tau^\diamond}\geq S_{\tau^\diamond}$ and we conclude that $\tilde Y_{\tau^\diamond}=S_{\tau^\diamond}$. In particular, the stopping time $\tilde\tau:=\inf\{t\geq\eta:\tilde Y_t=S_t\}$ satisfies $\eta\leq \tilde\tau\leq \tau^\diamond$, $\Prob$-a.s.

\begin{lem}
Let $\eta$ and $\tilde\tau$ be as defined above, then $\tilde Y_\eta=Y^{\tilde\tau}_\eta$.
\end{lem}

\noindent\emph{Proof.} First, since $Y^n_\eta\geq Y^{\tilde\tau}_\eta$, we get that $\tilde Y_\eta\geq Y^{\tilde\tau}_\eta$. To establish that the lemma holds we show that there exists a triple of processes $(\tilde Z,\tilde V,\tilde K)\in \mcH^{2,\tilde\tau}(W)\times  \mcH^{2,\tilde\tau}(\mu)\times\mcA^{2,\tilde\tau}$ such that the quadruple $(\tilde Y,\tilde Z,\tilde V,\tilde K)$ is a solution to \eqref{ekv:bsde-c-jmp-rec-2}. Our proof follows along the lines of the proofs detailed in Section 3.2 of \cite{Kharroubi2010}. By dominated convergence we get that
\begin{align}\label{ekv:Yn-conv}
  \E\big[\int_0^T|\tilde Y_t- Y^n_t|^2dt+|\tilde Y_{\tilde\tau}- Y^n_{\tilde\tau}|^2|\big]\to 0
\end{align}
as $n\to\infty$. Letting $f^n_s:=f(s,Y^n_s,Z^n_s,V^n_s)$, Lemma~\ref{lem:bound-Yn} and the fact that $\tilde\tau\leq \tau_n$, $\Prob$-a.s., for each $n\in\bbN$ implies that there is a triple of processes $(\phi,\tilde Z,\tilde V)$, where $\phi$ is $\R$-valued, $\Prog(\bbF)$-measurable and $dt\otimes d\Prob$-square integrable and $(\tilde Z,\tilde V)\in  \mcH^{2,\tilde\tau}(W)\times  \mcH^{2,\tilde\tau}(\mu)$, such that (up to a subsequence) for any stopping time $\tau\in\mcT_{[\eta,\tilde\tau]}$, the following weak convergence holds in $L^2(\Omega,\mcF_\tau)$
\begin{align*}
  &\int_\eta^\tau f^n_sds\rightharpoonup\int_\eta^\tau \phi_s ds,\qquad \int_\eta^\tau Z^n_sdW_s \rightharpoonup \int_\eta^\tau \tilde Z_sdW_s
\\
&\int_\eta^\tau\!\!\!\int_U V^n_s(e) \mu(ds,de)\rightharpoonup \int_\eta^\tau\!\!\!\int_U \tilde V_s(e) \mu(ds,de).
\end{align*}
Moreover, since $K^{+,n}_{\tilde\tau}-K^{+,n}_\eta=0$ for all $n\in\bbN$ we have that
\begin{align*}
  K^{-,n}_{\tau}-K^{-,n}_{\eta}\rightharpoonup \tilde K^-_{\tau}:=\tilde Y_{\tau}-\tilde Y_\eta+\int_\eta^\tau \phi_s ds-\int_\eta^\tau \tilde Z_s dW_s-\int_\eta^\tau\!\!\!\int_U \tilde V_s(e)\mu(ds,de)
\end{align*}
and letting $\tilde K^-_{t}\equiv 0$ on $[0,\eta]$ gives a non-decreasing, predictable, \cadlag process that thus belongs to $\mcA^{2,\tilde\tau}$. In particular, we find that
\begin{align*}
  \tilde Y_{t}=\tilde Y_\eta-\int_\eta^t \phi_s ds+\tilde K_t+\int_\eta^t \tilde Z_s dW_s+\int_\eta^t\!\!\!\int_U \tilde V_s(e)\mu(ds,de),\quad \forall t\in [\eta,\tilde\tau]
\end{align*}
In order to show that we in fact have $\phi_s=f(s,\tilde Y_s,\tilde Z_s,\tilde V_s)$, $ds\otimes d\Prob$-a.e.~and that $\ett_{[\eta,\tilde\tau]}V_s(e)\geq -\ett_{[\eta,\tilde\tau]}\chi(s,\tilde Y_{s-},\tilde Z_s,e)$, $dt\otimes d\Prob\otimes \lambda(de)$-a.e.~we show that $(\ett_{[\eta,\tilde\tau]}Z^n,\ett_{[\eta,\tilde\tau]}V^n)$ converges to $(\ett_{[\eta,\tilde\tau]}\tilde Z,\ett_{[\eta,\tilde\tau]}\tilde V)$ in $\mcH^{p,\tilde\tau}(W)\times \mcH^{p,\tilde\tau}(\mu)$ for some $p\in (1,2)$. Applying It\^o's formula to $|Y^n-\tilde Y|^2$ gives
\begin{align*}
 |Y^n_{\tilde\tau}- \tilde Y_{\tilde\tau}|^2&= |Y^n_{\eta}- \tilde Y_{\eta}|^2 - 2\int_\eta^{\tilde\tau} (Y^n_{s}- \tilde Y_{s}) (f^n_s-\phi_s)ds + 2\int_\eta^{\tilde\tau}(Y^n_{s}- \tilde Y_{s})(Z^n_s-\tilde Z_s)dW_s
\\
&\quad + 2\int_\eta^{\tilde\tau} (Y^n_{s}- \tilde Y_{s})dK^{-,n}_s - 2\int_\eta^{\tilde\tau} (Y^n_{s-}- \tilde Y_{s-}+\Delta \tilde K^-_s)d\tilde K^-_s-\sum_{s\in[\eta,\tilde\tau]}|\Delta\tilde K^-_s|^2 + \int_\eta^{\tilde\tau}|Z^n_s-\tilde Z_s|^2 ds
  \\
  &\quad+\int_\eta^{\tilde\tau}\!\!\!\int_U (2(Y^n_{s-}-\tilde Y_{s-})(V^n_{s-}(e)-\tilde V_{s-}(e))+|V^n_{s-}(e)-\tilde V_{s-}(e)|^2)\mu(ds,de).
\end{align*}
Rearranging, taking expectation and using the inequality $2ab\leq a^2+b^2$ together with the fact that $(Y^n_{s}- \tilde Y_{s})dK^{-,n}_s\geq 0$, we find that
\begin{align*}
 &\E\Big[\int_\eta^{\tilde\tau}|Z^n_s-\tilde Z_s|^2ds\Big]+\frac{1}{2}\E\Big[\int_\eta^{\tilde\tau}\!\!\!\int_U |V^n_{s-}(e)-\tilde V_{s-}(e)|^2\lambda(de)ds\Big]
 \\
 &\leq \E\Big[|Y^n_{\tilde\tau}- \tilde Y_{\tilde\tau}|^2 +C\int_\eta^{\tilde\tau} |Y^n_{s}- \tilde Y_{s}|^2ds +2\int_\eta^{\tilde\tau} (Y^n_{s}- \tilde Y_{s}) (f^n_s-\phi_s)ds
\\
&\quad+  2\int_\eta^{\tilde\tau} |Y^n_{s-}- \tilde Y_{s-}+\Delta \tilde K^-_s|d\tilde K^-_s + \sum_{s\in[\eta,\tilde\tau]}|\Delta\tilde K^-_s|^2\Big].
\end{align*}
By \eqref{ekv:Yn-conv} and uniform $L^2([0,T]\times\Omega,\mcB([0,T])\otimes\mcF_T,dt\otimes d\Prob)$-boundedness of $\ett_{[\eta,\tilde\tau]}(f^n-\phi)$ we conclude that the first three terms tend to zero as $n\to\infty$. Then, by applying the approximation routine of \cite{PengMLT} as was done in the proof of Theorem 3.1 in \cite{Kharroubi2010} we find that $(\ett_{[\eta,\tilde\tau]}Z^n,\ett_{[\eta,\tilde\tau]}V^n)$ converges to $(\ett_{[\eta,\tilde\tau]}\tilde Z,\ett_{[\eta,\tilde\tau]}\tilde V)$ in measure which by uniform integrability gives strong converges in $\mcH^{p}(W)\times \mcH^{p}(\mu)$ for $p\in [1,2)$.

It thus holds that
\begin{align*}
    \tilde Y_t=S_{\tilde\tau}+\int_t^{\tilde\tau} f(s,\tilde Y_s,\tilde Z_s,\tilde V_s)ds-\int_t^{\tilde\tau} \tilde Z_s dW_s-\int_t^{\tilde\tau}\!\!\!\int_U \tilde V_s(e)\mu(ds,de)-(\tilde K^{-}_{\tilde\tau}-\tilde K^{-}_t),\quad\forall t\in[\eta,\tilde\tau].
\end{align*}
Moreover, since $\chi$ is Lipschitz,
\begin{align*}
  \E\Big[\int_\eta^{\tilde\tau}\!\!\!\int_U(V^{n}_s(e)+\chi(s,Y^{n}_s,Z^{n}_s,e))^-\lambda(de)ds\Big]\to \E\Big[\int_\eta^{\tilde\tau}\!\!\!\int_U(\tilde V_s(e)+\chi(s,\tilde Y_s,\tilde Z_s,e))^-\lambda(de)ds\Big]
\end{align*}
as $n\to\infty$ and since $K^{n}_{\tilde\tau}-K^{n}_{\eta}$ is bounded in $L^2(\Omega,\mcF_T)$, uniformly in $n\in\bbN$, we conclude that
\begin{align*}
  \E\Big[\int_\eta^{\tilde\tau}\!\!\!\int_U(\tilde V_s(e)+\chi(s,\tilde Y_s,\tilde Z_s,e))^-\lambda(de)ds\Big]=0
\end{align*}
from which it follows that $\ett_{[\eta,\tilde\tau]}\tilde V_s(e)\geq-\ett_{[\eta,\tilde\tau]}\chi(s,\tilde Y_{s-},\tilde Z_s,e)$, $d\Prob\otimes\lambda(de)\otimes ds$-a.e.

We thus conclude that $(\tilde Y,\tilde Z,\tilde V,\tilde K^-)$ solves \eqref{ekv:bsde-c-jmp-rec-2} for $\tau=\tilde\tau$ on the interval $[\eta,\tilde\tau]$. On the other hand, $(Y^{\tilde\tau},Z^{\tilde\tau},V^{\tilde\tau},K^{\tilde\tau,-})$ is the unique maximal solution to this BSDE and since $\tilde Y\geq Y^{\tilde\tau}$, we conclude that $(\ett_{[\eta,\tilde\tau]}\tilde Y,\ett_{[\eta,\tilde\tau]}\tilde Z,\ett_{[\eta,\tilde\tau]}\tilde V,\tilde K^-)=(\ett_{[\eta,\tilde\tau]}Y^{\tilde\tau},\ett_{[\eta,\tilde\tau]}Z^{\tilde\tau},\ett_{[\eta,\tilde\tau]} V^{\tilde\tau}, K^{\tilde\tau,-}_{\cdot\vee \eta}-K^{\tilde\tau,-}_\eta)$ in $\mcS^{2,\tilde\tau}\times\mcH^{2,\tilde\tau}(W)\times\mcH^{2,\tilde\tau}(\mu)\times\mcS^{2,\tilde\tau}$.\qed\\

\emph{Proof of Theorem~\ref{thm:nl-Snell}.} It remains to demonstrate that $\tilde Y_\eta\geq Y^{\tau}_\eta$ for any $\tau\in\mcT_\eta$.  We observe that $Y^n_\eta\geq Y^{n,\tau}_\eta$ for each $n\in\bbN$ and $\tau\in\mcT_\eta$, where the quadruple $(Y^{\tau,n},Z^{\tau,n},V^{\tau,n},K^{+,\tau,n})\in\mcS^{2,\tau}\times\mcH^{2,\tau}(W)\times\mcH^{2,\tau}(\mu)\times\mcA^{2,\tau}$ is the unique solution to
\begin{align}\nonumber
    Y^{\tau,n}_t&=S_{\tau}+\int_t^\tau f(s,Y^{\tau,n}_s,Z^{\tau,n}_s,V^{\tau,n}_s)ds-\int_t^\tau Z^{\tau,n}_s dW_s-\int_t^\tau\!\!\!\int_U V^{\tau,n}_s(e)\mu(ds,de)
    \\
    &\quad-n\int_t^\tau\!\!\!\int_U(V^{\tau,n}_s(e)+\chi(s,Y^{\tau,n}_s,Z^{\tau,n}_s,e))^-\lambda(de)ds,\quad \forall t\in[0,\tau].\label{ekv:bsde-pen-n}
\end{align}
On the other hand, repeating the steps above it easily follows that $Y^{\tau,n}\searrow Y^\tau$ pointwisely and we conclude that
\begin{align*}
  \tilde Y_\eta=\lim_{n\to\infty}Y^n_\eta\geq \lim_{n\to\infty}Y^{\tau,n}_\eta=Y^\tau_\eta.
\end{align*}
In particular, $(\tilde Y,\tilde\tau)$ fulfills the condition regarding $(Y,\tau^*)$ as stated in the theorem.\qed\\

\begin{rem}
An additional point of significance is that the sequence of reflected BSDEs, represented by \eqref{ekv:rbsde-pen-n}, presents an efficient means for numerically approximating the non-linear Snell envelope, $Y$. Alternatively, one can first approximate $\tau^*$ using the sequence $\tau_n$ and then employ a numerical scheme for BSDEs with constrained jumps (see \eg \cite{Kharroubi14}) to find an efficient numerical approximation of $Y$.
\end{rem}


\section{A zero-sum game of impulse control versus stopping\label{sec:zsg}}
We delve into formulating the game of impulse control versus stopping, which is closely connected to the problem of optimal stopping of BSDEs with constrained jumps treated above and, in particular, to the corresponding non-linear Snell envelope. In this section, we will precisely define the problem, state the main result and provide preliminary estimates for the involved processes. Subsequently, in the following section, we will prove that the game has a value by establish a relationship between the upper and lower value functions and the aforementioned non-linear Snell envelope.

\subsection{Additional notations and definitions}
We introduce the following additional notations:
\begin{itemize}
  \item We let $\mcT^W$ be the set of all $\bbF^W$-stopping times valued in $[0,T]$ and for each $t\in [0,T]$, we let $\mcT^W_t$ be the set of $\bbF^{t,W}$-stopping times $\tau$ such that $\tau\geq t$, $\Prob$-a.s.
    \item We let $\boldD$ be the set of \cadlag paths $[0,T]\to\R^d$ with a finite number of jumps and for each $t\in[0,T]$, we introduce the semi-norm $\|x\|_t:=\sup_{s\in[0,t]}|x_s|$ and the semi-metric on $[0,T]\times \boldD$,
      \begin{align*}
        \mathbf d[(t,x),(t',x')]&:=|t-t'|+\|\mcC(x_{\cdot\wedge t})-\mcC(x'_{\cdot\wedge t'})\|_T+\sum_{j=1}^{\infty}(|\mcJ_i(x)\ett_{[\mcJ_i(x)\leq t]}-\mcJ_i(x')\ett_{[\mcJ_i(x')\leq t']}|
        \\
        &\quad+|\Delta x_{\mcJ_i(x)}\ett_{[\mcJ_i(x)\leq t]}-\Delta x'_{\mcJ_i(x')}\ett_{[\mcJ_i(x')\leq t']}|),
      \end{align*}
      where $\mcC(x)$ is the continuous part of $x$, $\mcJ_i(x)$ is the time of the $i^{\rm th}$ jump of $x$ and $\Delta x_t:=x_t-x_{t-}$.
  \item We introduce the filtration $\bbD:=(\mcD_t)_{t\in[0,T]}$, where $\mcD_t$ is the $\sigma$-algebra generated by the canonical maps $\boldD\to\R^d:x\mapsto x(s)$ for $s\in [0,t]$.
  \item We let $\mcU$ be the set of all $u=(\eta_j,\beta_j)_{j=1}^N$, where $(\eta_j)_{j=1}^\infty$ is a non-decreasing sequence of $\bbF$-stopping times, $\beta_j$ is a $U$-valued, $\mcF_{\eta_j}$-measurable random variable and $N=N^u_T:=\sup\{j\geq 0:\eta_j\leq T\}$ is $\Prob$-a.s.~finite. Moreover, for $t\in[0,T]$, we let $\mcU_t$ be the subset of $\mcU$ with $\eta_1\geq t$, $\Prob$-a.s.
  \item We let $\mcU^W_t$ be the subset of $\mcU_t$, where the $\eta_j$ are $\bbF^{t,W}$-stopping times and $\beta_j$ is $\mcF^{t,W}_{\eta_j}$-measurable.
  \item For $k\geq 0$ and $t\in[0,T]$, we let $\mcU^{k}_t$ (resp.~$\mcU^{W,k}_t$) be the set of all $u:=(\eta_j,\beta_j)_{j=1}^N$ in $\mcU_t$ (resp.~$\mcU^W_t$) for which $N\leq k$, $\Prob$-a.s.
  \item For any $\tau\in\mcT$ and $u=(\eta_j,\beta_j)_{j=1}^N\in\mcU$ we let $N^u_\tau:=\max\{j\geq 0:\eta_j\leq\tau\}$ and set $u_\tau:=[u]_{N^u_\tau}=(\eta_j,\beta_j)_{j=1}^{N^u_\tau}$. Similarly, $u_{\tau-}:=[u]_{N^u_{\tau-}}$, with $N^u_{\tau-}:=\max\{j\geq 0:\eta_j<\tau\}$.
  \item We let $\mathbf V_t:=\{(t_i,b_i)_{i=1}^\infty\in\mathbf V:t_1\geq t\}$ and define the concatenation operator $\otimes_t$ as
  \begin{align*}
    (t_i,b_i)_{i=1}^\infty\otimes_t (\tilde t_i,\tilde b_i)_{i=1}^\infty:=((t_1,b_1),\ldots,(t_{n_t},b_{n_t}),(\tilde t_1,\tilde b_1),(\tilde t_2,\tilde b_2),\ldots)
  \end{align*}
  for all $(t_i,b_i)_{i=1}^\infty\in\mathbf V$ and $(\tilde t_i,\tilde b_i)_{i=1}^\infty\in\mathbf V_t$.
\end{itemize}
Note that an alternative, and indeed a more general, metric on the space of \cadlag paths that we could have used as a basis for our definition of $\mathbf d$ is the Skorokhod metric. However, with the type of temporal distortions that we expect, $\mathbf d$ proves to be more convenient and practical to utilize.

\begin{defn}
For $t\in [0,T]$, the set of non-anticipative stopping strategies, denoted $\mcT^{S,W}_t$, is defined as the set of maps $\tau^S:\mcU^W_t\to\mcT^W_t$ such that for any $\eta\in\mcT^W_t$ and $u,\tilde u\in \mcU^W_t$, the difference between the sets $\{\omega:\tau^S(u)\leq \eta, u_\eta=\tilde u_\eta\}$ and $\{\omega:\tau^S(\tilde u)\leq \eta, u_\eta=\tilde u_\eta\}$ is a $\Prob$-null set.
\end{defn}

\begin{defn}
For $t\in [0,T]$, the set of non-anticipative impulse control strategies $\mcU^{S,W}_t$ (resp.~$\mcU^{S,W,k}_t$) is defined as the set of maps $u^S:\mcT^W_t\to\mcU^W_t$ (resp.~$u^S:\mcT^W_t\to\mcU^{W,k}_t$) such that for any $\tau,\tau'\in\mcT^W_t$, there is a $\Prob$-null set $E$ for which $u^S(\tau)_{(\tau\wedge\tau')-}=u^S(\tau')_{(\tau\wedge\tau')-}$ outside of $E$ and $u^S(\tau)=u^S(\tau')$ on $\{\tau=\tau'\}\setminus E$.
\end{defn}

\subsection{Problem formulation and main result}
Let us recall the definition of the cost/reward process. For $t \in [0, T]$, $u \in \mathcal{U}_t$, and $\tau \in \mathcal{T}_t$, the random variable $J_t(u, \tau)$ was defined as:
\begin{align*}
  J_t(u,\tau)=\E\Big[\Psi(\tau,X^{t,u})+\int_t^{\tau} f(s,X^{t,u})ds+\sum_{j=1}^{N}\ett_{[\eta_ j\leq\tau]}\chi(\eta_j,X^{t,[u]_{j-1}},\beta_j)\Big|\mcF_t\Big],
\end{align*}
where we recall that \eqref{ekv:fwd-sde-2} defined $X^{t,u}$ as the soluton to
\begin{align*}
X^{t,u}_s&=x_0+\int_0^sa(r,X^{t,u})dr+\int_0^s\sigma(r,X^{t,u})dW_r+\int_0^{t\wedge s}\!\!\!\int_U\gamma(r,X^{t,u},e)\mu(dr,de)
\\
&\quad+\sum_{j=1}^N \ett_{[\eta_j\leq s]}\gamma(\eta_j,X^{t,[u]_{j-1}},\beta_j),\quad\forall s\in [0,T].
\end{align*}
We then define the lower value process
\begin{align}\label{ekv:lvf-def}
\underline Y_t:=\essinf_{u^S\in\mcU^{S,W}_t}\esssup_{\tau\in\mcT^W_t}J_t(u^S(\tau),\tau)
\end{align}
and the upper value process
\begin{align}\label{ekv:uvf-def}
\bar Y_t:=\esssup_{\tau^S\in\mcT^{S,W}_t}\essinf_{u\in\mcU^W_t}J_t(u,\tau^S(u))
\end{align}
for each $t\in[0,T]$.

The main result of the second part of the paper is the following:

\begin{thm}\label{thm:sdg}
Under the assumptions detailed in Section~\ref{sub-sec:SDG-ass} and (H.\ref{hyp:SFDEflow})-(H.\ref{hyp:unif-conv}) below, there is a process $Y\in\mcS^2$ such that $Y_\eta=\esssup_{\tau\in\mcT_\eta}Y^\tau_\eta$ for any $\eta\in\mcT$, where the triple $(Y^\tau,Z^\tau,V^\tau,K^{-,\tau})\in\mcS^2\times\mcH^2(W)\times\mcH^2(\mu)\times\mcA^2$ is the maximal solution to \eqref{ekv:bsde-c-jmp-lin}. Moreover, the upper and lower value processes satisfy $\underline Y_t=\bar Y_t=Y_t$, $\Prob$-a.s., for each $t\in [0,T]$.
\end{thm}

A comment regarding the use of filtrations seems appropriate in this context. For $t=0$, as $\mathcal{F}_0$ is the trivial filtration, we encounter a SDG of impulse control versus stopping in a Brownian filtration. Our motivation for investigating the game formulation in a conditional setting, allowing $t\in (0,T]$, goes beyond providing a representation for the non-linear Snell envelope discussed in Section~\ref{sec:Snell}. It also serves as a foundational element for considering more general scenarios. Subsequent research can build upon this framework, similar to the extension of reflected BSDEs, known to be associated with classical control versus stopping games~\cite{HamLep2000,Bayraktar11}, to encompass control problems involving a combination of switching and classical control~\cite{HuTang}, as well as impulse control versus classical control games~\cite{Perninge2022}.  Furthermore, the conditional framework facilitates the efficient use of the penalization routine described in Section~\ref{sec:Snell} to approximate the value function over a wide range of historical information.

\subsection{Assumptions\label{sub-sec:SDG-ass}}
To be able to represent an impulse control $u\in\mcU^W_t$ by a randomized control, defined in terms of the random measure $\mu$, we need the following assumption:
\begin{ass}\label{ass:lambda}
$\lambda$ is a finite positive measure on $(U,\mcB(U))$ with full topological support.
\end{ass}
We assume that the coefficients of the SDE satisfy the following:
\begin{ass}\label{ass:onSFDE}
There is a $C>0$ such that for any $t,t'\geq 0$, $b,b'\in U$ and $x,x'\in\boldD$, we have:
\begin{enumerate}[i)]
  \item\label{ass:onSFDE-Gamma} The function $\gamma:[0,T]\times\boldD\times U\to\R^d$ is $\Prog(\bbD)\otimes\mcB(U)$-measurable, such that $\chi(\cdot,x,b)$ is \cadlag and satisfies the growth condition
    \begin{align*}
      |x_t+\gamma(t,x,b)|\leq C\vee |x_t|.
    \end{align*}
  \item\label{ass:onSFDE-a-sigma} The coefficients $a:[0,T]\times\boldD\to\R^{d}$ and $\sigma:[0,T]\times\boldD\to\R^{d\times d}$ are $\Prog(\bbD)$-measurable and satisfy the growth condition
    \begin{align*}
      |a(t,x)|+|\sigma(t,x)|&\leq C(1+\|x\|_t)
    \end{align*}
    and the Lipschitz continuity
    \begin{align*}
      |a(t,x)-a(t,x')|+|\sigma(t,x)-\sigma(t,x')|&\leq C\|x'-x\|_t.
    \end{align*}
\end{enumerate}
\end{ass}

The coefficients of the reward/cost satisfy the following assumptions:

\begin{ass}\label{ass:oncoeff}
There are constants $C>0$ and $q>0$ in addition to a family of modulus of continuity functions $(\varpi_K)_{K\geq 0}$ such that for any $t,t'\geq 0$, $b,b'\in U$ and $x,x'\in\boldD$, we have:
\begin{enumerate}[i)]
  \item\label{ass:oncoeff-f} The running cost $f:[0,T]\times \boldD\to\R$ is $\Prog(\bbD)$-measurable and satisfies the growth condition
  \begin{align*}
    |f(t,x)|\leq C(1+\|x\|_t^q)
  \end{align*}
  and for any $K>0$
  \begin{align*}
    |f(t,x)-f(t,x')|\leq \varpi_K(\mathbf d[(t,x),(t,x')])
  \end{align*}
  whenever $\|x\|_t\vee \|x'\|_t\leq K$.
  \item\label{ass:oncoeff-ell} The intervention cost $\chi:[0,T]\times \boldD\times U\to [0,\infty)$ is $\Prog(\bbD)\otimes\mcB(U)$-measurable with $\chi(\cdot,x,b)$ \cadlagSTOP, of polynomial growth, \ie
  \begin{align*}
    |\chi(t,x,b)|\leq C(1+\|x\|_t^q),
  \end{align*}
  and satisfies
  \begin{align*}
    \chi(t',x',b')-\chi(t,x,b)\leq \varpi_K(\mathbf d[(t,x),(t',x')]+|b-b'|)
  \end{align*}
  whenever $\|x\|_t\vee \|x'\|_{t'}\leq K$ and $t\leq t'$.
  \item\label{ass:@stopping} The barrier function $\Psi:[0,T]\times \boldD\to \R$ is $\Prog(\bbD)$-measurable with $\Psi(\cdot,x)$ \cadlagSTOP, of polynomial growth, \ie $|\Psi(t,x)|\leq C(1+\|x\|_t^q)$ and satisfies the regularity property
      \begin{align*}
        \Psi(t,x)-\Psi(t',x')\leq \varpi_K(\mathbf d[(t,x),(t',x')])
      \end{align*}
      whenever $\|x\|_t\vee \|x'\|_{t'}\leq K$ and $t\leq t'$. In addition, we assume that
      \begin{align}\label{ekv:ass@stopping}
        \Psi(t,x)\leq \Psi(t,x_\cdot+\ett_{[\cdot\geq t]}\gamma(t,x,b))+\chi(t,x,b).
      \end{align}
\end{enumerate}
\end{ass}

In impulse control it is generally assumed that intervening on the system at the end of the horizon is suboptimal. Inequality \eqref{ekv:ass@stopping} extends this assumption and ascertain that it is never optimal for the impulse controller to intervene on the system at the time that the game ends. Note also that we may (and will), without loss of generality assume that there is a $C>0$ such that $\varpi_K$ is bounded by $C(1+K^q)$ for any $K\geq 0$.

\subsection{Preliminary estimates}

\begin{lem}\label{lem:SFDEmoment}
Under Assumption~\ref{ass:onSFDE}, the path-dependent SDE \eqref{ekv:fwd-sde-2} admits a unique solution for each $t\in[0,T]$ and $u\in\mcU_t$. Furthermore, the solution has moments of all orders, in particular we have for $p\geq 0$, that
\begin{align}\label{ekv:SFDEmoment}
\E\Big[\|X^{u}_s\|_T^{p}\Big|\mcF_t\Big]\leq C(1+\|X^u\|_t^p),
\end{align}
where $C>0$ does not depend on $u$ and $\lambda$ and $X^{u}:=X^{0,u}$.
\end{lem}

\noindent\emph{Proof.} The result was proved for a Brownian filtration in \cite{Perninge2022} (see Proposition 4.2) and the method of that proof readily extends to cover our framework.\qed\\

\begin{rem}\label{rem:SFDEmoment}
Since $(\sigma_j,\zeta_j)_{j\geq 1}\in\mcU$, the above proposition immediately gives that there is a $C>0$ such that $\esssup_{u\in\mcU_t}\E\big[\|X^{t,u}\|^p_T\big|\mcF_t\big]\leq C(1+\|X\|^p_t)$ for all $t\in [0,T]$.
\end{rem}

\begin{cor}
There is a $C>0$ such that
\begin{align}\label{ekv:vfs-bound}
|\underline Y_t| + |\bar Y_t| \leq C(1+\|X\|_t^q)
\end{align}
for all $t\in[0,T]$.
\end{cor}

\noindent\emph{Proof.} We have
\begin{align}
\Psi(t,X)\leq \,\underline Y_t,\,\bar Y_t\,\leq \E\Big[ \|\Psi(\cdot,X^{t,\emptyset})\|_T+\int_t^{T} |f(r,X^{t,\emptyset})|dr\Big|\mcF_t\Big]
\end{align}
and the assertion follows by the polynomial growth assumptions on $\Psi$ and $f$ and Remark~\ref{rem:SFDEmoment}.\qed\\

In addition to the assumptions stated in Section~\ref{sub-sec:SDG-ass}, we introduce two hypotheses that are formulated in a more implicit manner to allow for broader applicability of the results presented in the following section. These hypotheses can be demonstrated to hold under relatively mild additional assumptions. For instance, the first hypothesis is shown to hold in Lemma 4.3 of \cite{Perninge2022} by imposing a continuity requirement on $\gamma$ and $L^1$ (resp. $L^2$) continuity on $a$ (resp. $\sigma$). Furthermore, Proposition~\ref{prop:lvf-approx} below establishes that the second hypothesis holds when $\chi$ is bounded from below by a positive constant.

Moreover, it is worth noting that both hypotheses can be easily verified for a wide range of classical control versus stopping games if we initially approximate the classical control $\alpha$ using a piece-wise constant process $\hat\alpha_t := a_0\ett_{[0,\eta_1]}(t) + \sum_{j=1}^{N}\beta_j\ett_{(\eta_j,\eta_{j+1}]}(t)$, as done in \cite{Fuhrman15}.

\begin{hyp}[H.\ref{hyp:SFDEflow}]\label{hyp:SFDEflow}
For each $k\in\bbN$, there is a modulus of continuity $\rho_k$ such that
\begin{align*}
  \E\big[(\mathbf d[(T,X^{t,\tilde u}),(T,X^{t,u})])^2\big]\leq \rho_k(\eps),
\end{align*}
for all $t\in[0,T]$ and $u,\tilde u\in\mcU^k_t$ that satisfy $|u-\tilde u|\leq\eps$, $\Prob$-a.s.
\end{hyp}

Under (H.\ref{hyp:SFDEflow}) we have the following:

\begin{lem}\label{lem:SFDEflow}
Suppose that for some $t\in[0,T]$ and $k\in \bbN$, we pick sequences $\eps_l\searrow 0$ and $(u^l)_{l\in\bbN},(\tilde u^l)_{l\in\bbN}\subset\mcU^k_t$ so that $|u^l-\tilde u^l|\leq \eps_l$ and let $\tau_l\in\mcT_{\eta^l_k}$ and $\tilde\tau_l\in\mcT_{\tilde\eta^l_k}$ be such that $|\tau_l-\tilde\tau_l|\leq \eps_l$. Then, by possibly going to a subsequence, it holds that
\begin{align*}
  \mathbf d[(\tau_l,X^{t,u^l})-(\tilde\tau_l,X^{t,\tilde u^l})]\to 0,
\end{align*}
$\Prob$-a.s., as $l\to\infty$.
\end{lem}

\noindent\emph{Proof.} By (H.\ref{hyp:SFDEflow}) and Lemma~\ref{lem:SFDEmoment} together with the Burkholder-Davis-Gundy inequality it follows that there is a modulus of continuity $\rho$ and a constant $C_1>0$ such that
\begin{align*}
  \E\big[(\mathbf d[(\tau_l,X^{t,\tilde u^l}),(\tilde\tau_l,X^{t,u^l})])^2\big]\leq C_1\rho(\eps_l),\quad\forall l\in\bbN.
\end{align*}
Let $A_l:=\{(\mathbf d[(\tau_l,X^{t,\tilde u}),(\tilde\tau_l,X^{t,u})])^2>(C_1\rho(\eps_l))^{1/2}\}$ and note that $\Prob[A_l]\leq (C_1\rho(\eps_l))^{1/2}$. Hence, there is a subsequence $(l_j)_{j\in\bbN}$ such that $\sum_{j=1}^\infty \Prob[A_{l_j}]<\infty$ implying that $\limsup_j A_{l_j}$ is $\Prob$-negligible. In particular, $\mathbf d[(\tau_{l_j},X^{t, u^{l_j}}),(\tilde\tau_{l_j},X^{t,\tilde u^{l_j}})]\to 0$, $\Prob$-a.s.~as $j\to\infty$.\qed\\

The above lemma allows us to deduce the following important continuity result:
\begin{cor}\label{cor:rew-cont}
For each $k\in\bbN$ and $t\in [0,T]$, we have
\begin{align*}
  \limsup_{\eps\to 0}&\sup_{u,\tilde u\in\mcU^k_t,\,|u-\tilde u|\leq\eps}\sup_{\tau\in\mcT_{\eta_k}}\sup_{\tilde\tau\in\mcT_{[\tilde\eta_k\vee\tau,\tau+\eps]}}\E\Big[(\Psi(\tau,X^{t,u}) - \Psi(\tilde\tau,X^{t,\tilde u}))^+
  \\
  &+ \int_t^T|f(s,X^{t,u})-f(s,X^{t,\tilde u})|ds+\sup_{b,\tilde b\in U,\,|b-\tilde b|\leq\eps}(\chi(\tau,X^{t,u},b)-\chi(\tilde \tau,X^{t,\tilde u},\tilde b))^-\Big]=0.
\end{align*}
\end{cor}

\noindent\emph{Proof.} We treat only the term containing $\Psi$ as identical arguments can be used to for the other terms. We seek a contradiction and assume that there is a $\varrho>0$, a sequence $\eps_l$ such that $\eps_l\to 0$ as $l\to\infty$ and sequences $(u^l)_{l\in\bbN},(\tilde u^l)_{l\in\bbN}\subset \mcU^{k}_t$ with $|u^l-\tilde u^l|\leq \eps_l$ and $(\tau_l)_{l\in\bbN},(\tilde\tau_l)_{l\in\bbN}\subset \mcT$  with $\eta_k^l\leq \tau_l\leq \tilde \tau_l\leq \tau_l+\eps_l$ and $\tilde\tau_l\geq\tilde\eta^l_k$ such that
\begin{align}\label{ekv:geq-rho}
\E\big[(\Psi(\tau_l,X^{t,u^l})-\Psi(\tilde\tau_l,X^{t,\tilde u^l}))^+\big]\geq \varrho
\end{align}
for every $l\in\bbN$. By Assumption~\ref{ass:oncoeff}.(\ref{ass:@stopping}) we find that for any $K\geq 0$, it holds that
\begin{align*}
  \E\big[(\Psi(\tau_l,X^{t,u^l})-\Psi(\tilde\tau_l,X^{t,\tilde u^l}))^+\big]&\leq\E\big[\varpi_K(\mathbf d[(\tau_l,X^{t,u^l}),(\tilde\tau_l,X^{t,\tilde u^l})])\big]
  \\
  &\quad+\E\big[\ett_{[\|X^{t,u^l}\|_T\vee\|X^{t,\tilde u^l}\|_T> K]}|\Psi(\tau_l,X^{t,u^l})-\Psi(\tilde\tau_l,X^{t,\tilde u^l})|\big]
  \\
  &\leq  \E\big[\varpi_K(\mathbf d[(\tau_l,X^{t,u^l}),(\tilde\tau_l,X^{t,\tilde u^l})])\big]+\frac{C}{K}.
\end{align*}
Now, by possibly going to a subsequence we have by Lemma~\ref{lem:SFDEflow} that
\begin{align*}
  \varpi_K(\mathbf d[(\tau_l,X^{t,u^l},(\tilde\tau_l,X^{t,\tilde u^l})])\to 0
\end{align*}
$\Prob$-a.s. Since $\varpi_K$ is uniformly bounded, we can use dominated convergence to conclude that there is a subsequence such that
\begin{align*}
&\lim_{l\to\infty}\E\big[(\Psi(\tau_l,X^{t,u^l})-\Psi(\tilde\tau_l,X^{t,\tilde u^l}))^+\big]\leq \frac{C}{K}
\end{align*}
contradicting \eqref{ekv:geq-rho} as $K>0$ was arbitrary.\qed\\

We note that (H.\ref{hyp:SFDEflow}) presumes that the impulse controls have a limited number of interventions. The hypothesis becomes particularly significant when combined with the upper and lower value functions under a truncation on the maximal number of interventions presented in the next subsection.

\subsection{Truncation of the impulse control set}
Next, we develop useful approximations of the value processes by their counterparts with a truncated number of interventions. For this, we introduce
\begin{align}\label{ekv:lvf-k-def}
\underline Y^k_t:=\essinf_{u^S\in\mcU^{S,W,k}_t}\esssup_{\tau\in\mcT^W_t}J_t(u^S(\tau),\tau)
\end{align}
and
\begin{align}\label{ekv:uvf-k-def}
\bar Y^k_t:=\esssup_{\tau^S\in\mcT^{S,W}_t}\essinf_{u\in\mcU^{W,k}_t}J_t(u,\tau^S(u)).
\end{align}

\begin{hyp}[H.\ref{hyp:unif-conv}]\label{hyp:unif-conv}
For every $t\in[0,T]$, we have $\underline Y^k_t\searrow \underline Y_t$ and $\bar Y^k_t\searrow \bar Y_t$, $\Prob$-a.s.
\end{hyp}

As the following proposition shows, a vital example where (H.\ref{hyp:unif-conv}) is satisfied is when the intervention costs are strictly positive.

\begin{prop}\label{prop:lvf-approx}
Assume that $\chi$ is bounded from below by a positive constant, \ie $\chi(t,x,b)\geq\delta >0$, then there is a $C>0$ such that
\begin{align}\label{ekv:lvf-approx}
\underline Y^k_t \leq \underline Y_t+\frac{C}{\sqrt{k}}(1+\|X\|_t^{2q})
\end{align}
and similarly
\begin{align}\label{ekv:uvf-approx}
\bar Y^k_t \leq \bar Y_t+\frac{C}{\sqrt{k}}(1+\|X\|_t^{2q}),
\end{align}
$\Prob$-a.s., for all $k\in\bbN$ and $t\in [0,T]$.
\end{prop}

\noindent\emph{Proof.} We prove \eqref{ekv:uvf-approx} as this proof is slightly more involved. Let $\tau^S\in\mcT^{S,W}_t$ and note that for any  $u\in\mcU^W_t$, non-anticipativity gives that $\tau^S(u)=\tau^S(\tilde u)$, $\Prob$-a.s.~where $\tilde u$ equals $u$ up to and at the stopping time $\tau^S(u)$, but then makes no further interventions. We thus have that
\begin{align*}
\essinf_{u\in\mcU^W_t}J_t(u,\tau^S(u))=\essinf_{u\in\mcU^{\tau^S}_t}J_t(u,\tau^S(u))
\end{align*}
where $\mcU^{\tau^S}_t$ is the subset of $\mcU^W_t$ with all impulse controls that do not have any interventions after $\tau^S(u)$. For any $u\in\mcU^{\tau^S}_t$, we have
\begin{align*}
J_t(u,\tau^S(u))&\geq \E\Big[ -\|\Psi(\cdot,X^{t,u})\|_T-\int_t^{T} |f(r,X^{t,u})|dr+N\delta\Big|\mcF_t\Big]
\\
&\geq -C(1+\|X\|^q_t)+\delta\E\big[N\big|\mcF_t\big]
\end{align*}
There is thus a $C>0$ (that does not depend on $t$, $\tau^S$ and $u$) such that the control $u\in\mcU^{\tau^S}_t$ is dominated by the control $\emptyset\in\mcU^{\tau^S}_t$ (\ie the impulse control that makes no interventions in $[0,T]$) whenever
\begin{align*}
\E\big[N\big|\mcF_t\big] > C(1+\|X\|^q_t).
\end{align*}
Now, for $\eps>0$ let $u^\eps\in\mcU^{\tau^S}_t$ be such that
\begin{align*}
J_t(u^\eps,\tau^S(u^\eps))\leq \essinf_{u\in\mcU^W_t}J_t(u,\tau^S(u))+\eps.
\end{align*}
If $u^\eps$ is dominated by $\emptyset$, then
\begin{align*}
  \essinf_{u\in\mcU^{W,k}_t}J_t(u,\tau^S(u))\leq \essinf_{u\in\mcU^W_t}J_t(u,\tau^S(u))+\eps
\end{align*}
and we can choose a smaller $\eps$ until either $u^\eps$ dominates $\emptyset$ or if this never happens we have
\begin{align*}
  \essinf_{u\in\mcU^{W,k}_t}J_t(u,\tau^S(u))= \essinf_{u\in\mcU^W_t}J_t(u,\tau^S(u))
\end{align*}
We thus assume that $u^\eps$ dominates $\emptyset$ so that
\begin{align}\label{ekv:Neps-L1-bound}
\E\big[N^\eps\big|\mcF_t\big] \leq C(1+\|X\|^q_t),
\end{align}
where $N^\eps:=N^{u^\eps}_T$ is the number of interventions that $u^\eps$ makes in $[t,\tau^s(u^\eps)]$. Since $[u^\eps]_k\in\mcU^{W,k}_t$ and as non-anticipativity of $\tau^S$ implies that $\ett_{[N^\eps\leq k]}\tau^S(u^\eps)=\ett_{[N^\eps\leq k]}\tau^S([u^\eps]_k)$, we get that
\begin{align*}
&\essinf_{u\in\mcU^{W,k}_t}J_t(u,\tau^S(u))- \essinf_{u\in\mcU^W_t}J_t(u,\tau^S(u))\leq J_t([u^\eps]_k,\tau^S([u^\eps]_k))-J_t(u^\eps,\tau^S(u^\eps))+\eps
\\
&\leq \E\Big[ \ett_{[N^\eps>k]}\big(\Psi(\tau^S([u^\eps]_k),X^{t,[u^\eps]_k})-\Psi(\tau^S(u^\eps),X^{t,u^\eps})+\int_t^{T} (f(r,X^{t,[u^\eps]_k})-f(r,X^{t,u^\eps}))dr\big)\Big|\mcF_t\Big]+\eps
\\
&\leq C\E\Big[\ett_{[N^\eps>k]} (1+\|X^{t,[u^\eps]_k}\|_T^q+\|X^{t,u^\eps_k}\|_T^q)\Big|\mcF_t\Big]+\eps
\\
&\leq C\E\big[\ett_{[N^\eps>k]}\big|\mcF_t\big]^{1/2} \E\Big[1+\|X^{t,[u^\eps]_k}\|_T^{2q}+\|X^{t,u^\eps_k}\|_T^{2q}\Big|\mcF_t\Big]^{1/2}+\eps.
\end{align*}
Using \eqref{ekv:Neps-L1-bound} together with Remark~\ref{rem:SFDEmoment}, the second inequality follows as $\tau^S\in\mcT^{S,W}_t$ and $\eps >0$ were arbitrary. The first inequality follows analogously by first taking $\tau^S(u)=T$.\qed\\

\subsection{Discretization of the control sets}
Our approach to prove Theorem~\ref{thm:sdg} requires a discretization of the control sets and we introduce the following sets:
\begin{defn}\label{def:discretization}
For $\eps>0$:
\begin{itemize}
  \item We let $\iota\geq 0$ be the smallest integer such that $2^{-\iota}T\leq\eps$, set $n_{\bbT}^\eps:=2^{\iota}+1$ and introduce the discrete set $\bbT^\eps:=\{t^\eps_j:t^\eps_j=(j-1)2^{-\iota}T,j=1,\ldots,n_{\bbT}^\eps\}$, a discritization of $[0,T]$ with step-size $\Delta t^\eps:=2^{-\iota}T$.
  \item We let $(U^\eps_{l})_{l=1}^{n^\eps_U}$ be a Borel-partition of $U$ such that each $U^\eps_l$ has a diameter that does not exceed $\eps$ and there is a sequence $(b^\eps_l)_{l=1}^{n^\eps_U}$ with $b^\eps_l\in U^\eps_l$ and $b^\eps_l\notin \bar U^\eps_{l'}$ for all $l'\neq l$. We then denote by $\bar U^\eps:=\{b^\eps_1,\ldots,b^\eps_{n^\eps_U}\}$ the corresponding dicretization of $U$.
  \item We let $\mathbf V^{\eps,k}:=\{(t_j,b_j)_{j=1}^{k}:t_1\leq t_2\leq\cdots\leq t_k:(t_j,b_j)\in \bbT^\eps\times \bar U^\eps\}$ with restriction $\mathbf V^{\eps,k}_{[t,s]}:=\{(t_j,b_j)_{j=1}^{k}\in \mathbf V^{\eps,k}:t_j\in [t,s]\}$ to interventions in the interval $[t,s]\subset [0,T]$ and set $\mathbf V^{\eps,k}_{t}:=\mathbf V^{\eps,k}_{[t,T]}$.
\end{itemize}
\end{defn}
For $\eps>0$ and $k\in\bbN$, we then define
\begin{align*}
\mcU^{W,k,\eps}_t:=\{u\in\mcU^{W,k}_t: (\eta_j,\beta_j)\in \bbT^\eps\times \bar U^\eps\text{ for }j=1,\ldots,N,\:\Prob\text{-a.s.}\}
\end{align*}
For $u\in\mcU^{W,k}_t$, we let $\Xi^\eps(u):=\big(\xi_1^\eps(\eta_j), \xi^\eps_2(\beta_j)\big)_{j=1}^k$ with $\xi_1^\eps(s):= \min\{r\in \bbT^\eps:r\geq s\}$ and $\xi^\eps_2(b):=\sum_{l=1}^{n^\eps_U}b^\eps_l\ett_{U^\eps_l}(b)$ and note that $\Xi^\eps(u)\in\mcU^{W,k,\eps}_t$.

Similarly, we let $\mcT^{W,\eps}_t$ be the subset of $\mcT^W_t$ with all stopping times $\tau$ for which $\tau\in\bbT^\eps$, $\Prob$-a.s.
\begin{defn}
For $\eps>0$, let $\mcT^{S,W,\eps}_t$ be the subset of $\mcT^{S,W}_t$ containing all strategies $\tau^S:\mcU^W_t\to\mcT^{W,\eps}_t$ such that on the set $\{\Xi^\eps(u)_{t^\eps_i}=\Xi^\eps(\tilde u)_{t^\eps_i}\}$ the sets $\{\tau^S(u)=t^\eps_i\}$ and $\{\tau^S(\tilde u)=t^\eps_i\}$ differ by at most a $\Prob$-null set for $i=1,\ldots,n_\bbT^\eps$. Moreover, we let $\mcU^{S,W,\eps}_t$ (resp.  $\mcU^{S,W,k,\eps}_t$) be the subset of $\mcU^{S,W}_t$ (resp.  $\mcU^{S,W,k}_t$) with all $u^S:\mcT^W_t\to\mcU^{W,\eps}_t$.
\end{defn}

We then introduce the truncated lower value process with impulse control and stopping discretization as follows
\begin{align}\label{ekv:lvf-trunk-disc}
\underline Y^{k,\eps}_t:=\essinf_{u^S\in\mcU^{S,W,k,\eps}_t}\,\esssup_{\tau\in\mcT^{W,\eps}_t}J_t(u^S(\tau),\tau)
\end{align}
and the corresponding upper value process
\begin{align}\label{ekv:uvf-trunk-disc}
\bar Y^{k,\eps}_t:=\esssup_{\tau^S\in\mcT^{S,W,\eps}_t}\,\essinf_{u\in\mcU^{W,k,\eps}_t}J_t(u,\tau^S(u)).
\end{align}

\begin{lem}\label{lem:underY-eps}
For each $k\in\bbN$ and $t\in[0,T]$, we have $\underline Y^{k,\eps}_t\to \underline Y^{k}_t$ in $L^1(\Omega,\mcF_t,\Prob)$ as $\eps\to 0$.
\end{lem}

\noindent\emph{Proof.} By non-anticipativity, any strategy $u^S\in\mcU^S$ can be written $u^S(\tau):=u_{\tau-}\otimes_\tau \tilde u^S(\tau)$, with $u\in\mcU$ and $\tilde u^S(\tau):=(\tilde \eta^S_j(\tau),\tilde \beta^S_j(\tau))_{j=1}^{\tilde N^S(\tau)}\in\mcU_{\tau}$. On the other hand, we have
\begin{align*}
J_t(u^S(\tau),\tau)&=J_t(u_{\tau-}\otimes_\tau \tilde u^S(\tau),\tau)\geq J_t(u_{\tau-},\tau)
\end{align*}
by Assumption~\ref{ass:oncoeff}.(\ref{ass:@stopping}) and we conclude that we can disregard strategies with $\tilde N^S(\tau)>0$. In particular, this gives that
\begin{align*}
\underline Y^{k}_t=\essinf_{u\in\mcU^{W,k}_t}\,\esssup_{\tau\in\mcT^W_t}J_t(u_{\tau-},\tau)
\end{align*}
and, similarly, we have
\begin{align}
\underline Y^{k,\eps}_t=\essinf_{u\in\mcU^{W,k,\eps}_t}\esssup_{\tau\in\mcT^{W,\eps}_t}J_t(u_{\tau-},\tau).
\end{align}
Hence,
\begin{align*}
\underline Y^{k}_t-\underline Y^{k,\eps}_t&\leq\essinf_{u\in\mcU^{W,k,\eps}_t}\,\esssup_{\tau\in\mcT^W_t}J_t(u_{\tau-},\tau) - \essinf_{u\in\mcU^{W,k,\eps}_t}\esssup_{\tau\in\mcT^{W,\eps}_t}J_t(u_{\tau-},\tau)
\\
&\leq \esssup_{u\in\mcU^{W,k,\eps}_t}\,\esssup_{\tau\in\mcT^W_t}\{J_t(u_{\tau-},\tau)-J_t(u_{\xi^\eps_1(\tau)-},\xi^\eps_1(\tau))\}
\end{align*}
On the other hand, for each $u\in\mcU^{W,k,\eps}_t$ and $\tau\in\mcT^W_t$, we have $u_{\tau-}=u_{\xi_1^\eps(\tau)-}$ and get that
\begin{align*}
J_t(u_{\tau-},\tau)-J_t(u_{\xi^\eps_1(\tau)-},\xi^\eps_1(\tau))&= \E\Big[\Psi(\tau,X^{t,u_{\tau-}})-\Psi(\xi_1^\eps(\tau),X^{t,u_{\xi_1^\eps(\tau)-}})-\int_{\tau}^{\xi_1^\eps(\tau)} f(r,X^{t,u})dr\Big|\mcF_t\Big].
\end{align*}
Since $\esssup_{u\in\mcU^{W,k,\eps}_t}\,\esssup_{\tau\in\mcT^W_t}\{J_t(u_{\tau-},\tau)-J_t(u_{\xi^\eps_1(\tau)-},\xi^\eps_1(\tau))\}\geq 0$, $\Prob$-a.s., it follows by applying Corollary~\ref{cor:rew-cont} that $\E\big[(\underline Y^{k}_t-\underline Y^{k,\eps}_t)^+\big]\to 0$ as $\eps\to 0$.

We approach the opposite inequality by noting that
\begin{align*}
\underline Y^{k,\eps}_t-\underline Y^{k}_t&\leq\essinf_{u\in\mcU^{W,k,\eps}_t}\,\esssup_{\tau\in\mcT^{W,\eps}_t}J_t(u_{\tau-},\tau) - \essinf_{u\in\mcU^{W,k}_t}\esssup_{\tau\in\mcT^{W,\eps}_t}J_t(u_{\tau-},\tau)
\\
&\leq \esssup_{u\in\mcU^{W,k}_t}\,\esssup_{\tau\in\mcT^{W,\eps}_t}\{J_t(\Xi^\eps(u)_{\tau-},\tau) - J_t(u_{\tau-},\tau)\}.
\end{align*}
Now, for each $u\in\mcU^{W,k}_t$ and $\tau\in\mcT^{W,\eps}_t$, we have
\begin{align*}
&\E\big[J_t(\Xi^\eps(u)_{\tau-},\tau)-J_t(u_{\tau-},\tau)\big]
\\
&= \E\Big[\Psi(\tau,X^{t,\Xi^\eps(u)_{\tau-}})-\Psi(\tau,X^{t,u_{\tau-}})+\int_{t}^{\tau} (f(r,X^{t,\Xi^\eps(u)})-f(r,X^{t,u}))dr
\\
&\quad + \sum_{j=1}^{k}\ett_{[\xi^\eps_1(\eta_j) < \tau]}\big(\chi(\xi^\eps_1(\eta_j) ,X^{t,[\Xi^\eps(u)]_{j-1}},\xi^\eps_2(\beta_j))-\chi(\eta_j ,X^{t,[u]_{j-1}},\beta_j)\big)
\\
&\quad-\sum_{j=1}^{k}\ett_{[\eta_j< \tau \leq \xi^\eps_1(\eta_j)]}\chi(\eta_j ,X^{t,[u]_{j-1}},\beta_j)\Big].
\end{align*}
On the other hand, with
\begin{align*}
\tilde u:=\big(\xi^\eps_1(\eta_j) \wedge \tau,\xi^\eps_2(\beta_j)\big)_{j=1}^{N^u_{\tau-}}
\end{align*}
we get by repeatedly appealing to \eqref{ekv:ass@stopping} of Assumption~\ref{ass:oncoeff}.(\ref{ass:@stopping}) that
\begin{align*}
&\Psi(\tau,X^{t,u_{\tau-}})+\sum_{j=1}^{k}\ett_{[\eta_j < \tau \leq \xi^\eps_1(\eta_j)]}\chi(\eta_j ,X^{t,[u]_{j-1}},\beta_j)
\\
&= \Psi(\tau,X^{t,\tilde u})+\sum_{j=1}^{k}\ett_{[\eta_j < \tau \leq \xi^\eps_1(\eta_j)]}\chi(\tau ,X^{t,[\tilde u]_{j-1}},\xi^\eps_2(\beta_j))+\Psi(\tau,X^{t,u_{\tau-}})-\Psi(\tau,X^{t,\tilde u})
\\
&\quad+\sum_{j=1}^{k}\ett_{[\eta_j < \tau \leq \xi^\eps_1(\eta_j)]}\big(\chi(\eta_j ,X^{t,[u]_{j-1}},\beta_j)-\chi(\tau ,X^{t,[\tilde u]_{j-1}},\xi^\eps_2(\beta_j))\big)
\\
&\geq \Psi(\tau,X^{t,\Xi^\eps(u)_{\tau-}})+\Psi(\tau,X^{t,u_{\tau-}})-\Psi(\tau,X^{t,\tilde u})
\\
&\quad+\sum_{j=1}^{k}\ett_{[\eta_j < \tau \leq \xi^\eps_1(\eta_j)]}\big(\chi(\eta_j ,X^{t,[u]_{j-1}},\beta_j)-\chi(\tau ,X^{t,[\tilde u]_{j-1}},\xi^\eps_2(\beta_j))\big).
\end{align*}
Combined, this gives that
\begin{align*}
&\E\big[J_t(\Xi^\eps(u)_{\tau-},\tau)-J_t(u_{\tau-},\tau)\big]
\\
&\leq \E\Big[\Psi(\tau,X^{t,\tilde u})-\Psi(\tau,X^{t,u_{\tau-}})+\int_{t}^{\tau} |f(r,X^{t,\tilde u})-f(r,X^{t,u})|dr
\\
&\quad + \sum_{j=1}^{k}\ett_{[\eta_j < \tau]}\big(\chi(\xi^\eps_1(\eta_j)\wedge \tau ,X^{t,[\tilde u]_{j-1}},\xi^\eps_2(\beta_j))-\chi(\eta_j ,X^{t,[u]_{j-1}},\beta_j)\big)\Big].
\end{align*}
We can thus use Corollary~\ref{cor:rew-cont} to conclude that
\begin{align*}
  \E\big[(\underline Y^{k}_t-\underline Y^{k,\eps}_t)^-\big]\leq \sup_{u\in\mcU^{W,k}_t}\,\sup_{\tau\in\mcT^{W,\eps}_t}\E\big[J_t(\Xi^\eps(u)_{\tau-},\tau)-J_t(u_{\tau-},\tau)\big]\to 0
\end{align*}
as $\eps\to 0$ and the assertion follows.\qed\\

Similarly, we have the following result.
\begin{lem}\label{lem:barY-eps-lower}
For each $k\in\bbN$ and $t\in[0,T]$, we have $\bar Y^{k,\eps}_t\to \bar Y^{k}_t$ in $L^1(\Omega,\mcF_t,\Prob)$ as $\eps\to 0$.
\end{lem}

\noindent\emph{Proof.} We have
\begin{align*}
\bar Y^k_t-\bar Y^{k,\eps}_t&\leq\esssup_{\tau^S\in\mcT^{S,W}_t}\,\essinf_{u\in\mcU_t^{W,k,\eps}}J_t(u,\tau^S(u)) - \esssup_{\tau^S\in\mcT^{S,W,\eps}_t}\,\essinf_{u\in\mcU_t^{W,k,\eps}}J_t(u,\tau^S(u))
\\
&\leq \esssup_{\tau\in\mcT^W_t}\,\esssup_{u\in\mcU_t^{W,k,\eps}}\,(J_t(u,\tau) - J_t(u,\xi_1^\eps(\tau)))
\end{align*}
and we can repeat the argument in the first part of the proof of Lemma~\ref{lem:underY-eps} to conclude that $\limsup_{\eps\to 0}\E[(\bar Y^k_t-\bar Y^{k,\eps}_t)^+]= 0$. Moreover,
\begin{align*}
\bar Y^{k,\eps}_t-\bar Y^k_t&\leq\esssup_{\tau^S\in\mcT^{S,W,\eps}_t}\,\essinf_{u\in\mcU_t^{W,k,\eps}}J_t(u,\tau^S(u)) -\esssup_{\tau\in\mcT^{S,W,\eps}_t}\,\essinf_{u\in\mcU_t^{W,k}}J_t(u,\tau^S(u))
\end{align*}
and using that $\tau^S(u)=\tau^S(\Xi^\eps(u))$, $\Prob$-a.s., whenever $\tau^S\in \mcT^{S,W,\eps}_t$ we conclude that
\begin{align*}
\bar Y^{k,\eps}_t-\bar Y^k_t&\leq \esssup_{\tau\in\mcT^{W,\eps}_t}\,\esssup_{u\in\mcU_t^{W,k}}\,(J_t(\Xi^\eps(u),\tau) - J_t(u,\tau)).
\end{align*}
Repeating the argument in the second part of the proof of Lemma~\ref{lem:underY-eps} then completes the proof.\qed\\

\begin{lem}\label{lem:ost-for-uvf}
For each $k\in\bbN$, $\eps>0$ and $t\in[0,T]$ we have
\begin{align*}
  \underline Y^{k,\eps}_t\leq \bar Y^{k,\eps}_t,
\end{align*}
$\Prob$-a.s.
\end{lem}

\noindent\emph{Proof.} For $0\leq t\leq s\leq T$, $u\in\mcU_{t}$, $\tilde u:=(\tilde\eta_j,\tilde\beta_j)_{j=1}^{\tilde N}\in\mcU_s$ and $\tau\in\mcT_s$, we introduce the random variable
\begin{align}
  J^{t,u}_s(\tilde u,\tau):=\E\Big[\Psi(\tau,X^{t,u\otimes_s \tilde u})+\int_s^{\tau} f(r,X^{t,u\otimes_s \tilde u})dr+\sum_{j=1}^{\tilde N}\ett_{[\tilde \eta_ j\leq\tau]}\chi(\tilde\eta_j,X^{t,u\otimes_s [\tilde u]_{j-1}},\tilde\beta_j)\Big|\mcF^{t,W}_s\Big],\label{ekv:obj-fun-2}
\end{align}
which gives us the cost/reward of using the pair $(\tilde u,\tau)$ given that the minimizer has applied the control $u$ on the interval $[t,s]$. We let $\tilde \mcU^{S,W,k,\eps}_s:=\{u^S\in\mcU^{S,W,k,\eps}_t:\eta_1(\tau)\geq s,\,\forall\tau\in\mcT\}$ and $\tilde \mcT^{W,\eps}_s:=\{\tau\in\mcT^{W,\eps}_t:\tau\geq s\}$ and for $\vecv\in\mathbf V_{[t,s]}^{\eps}$, we define
\begin{align*}
  R^{\vecv,k}_s:=\essinf_{u^S\in\tilde\mcU^{S,W,k,\eps}_s}\,\esssup_{\tau\in\tilde\mcT^{W,\eps}_s}J^{t,\vecv}_s(u^S(\tau),\tau)
\end{align*}
and immediately get that $\underline Y^{k,\eps}_t=R^{\emptyset,k}_t$. Moreover, to simplify notation we let
\begin{align*}
  R^{\vecv,k,+}_{t^\eps_i}:=&\essinf_{u^S\in\tilde\mcU^{S,W,k,\eps}_{t^\eps_{i+1}}}\,\esssup_{\tau\in\tilde\mcT^{W,\eps}_{t^\eps_{i+1}}} J^{t,\vecv}_{t^\eps_{i}}(u^S(\tau),\tau)
  \\
  =&\,\E\Big[\int_{t^\eps_{i}}^{t^\eps_{i+1}}f(r,X^{t,\vecv})dr +R^{\vecv,k}_{t^\eps_{i+1}}\Big|\mcF^{t,W}_{t^\eps_i}\Big],\quad \text{for }i=1,\ldots,n^\eps_{\bbT}
\end{align*}
The family $(R^{\vecv}_s:(\vecv,s)\in\mathbf V_{t}^{\eps}\times \bbT^\eps\cap[t,T])$ satisfies the dynamic programming principle
\begin{align}
  \begin{cases}
  R^{\vecv,k}_{T}=\Psi(T,X^{t,\vecv})
  \\
  R^{\vecv,k}_{t^\eps_i}=\Psi(t^\eps_i,X^{t,\vecv}) \vee \inf_{\vecb\in\bar U^{\eps,k}}\{R^{\vecv\otimes_{t^\eps_i} (t^\eps_i,\vecb),k-|\vecb|,+}_{t^\eps_i}+\chi({t^\eps_i},X^{t,\vecv},\vecb)\},\quad\forall t^\eps_i\in \bbT^\eps\cap[t,T),
  \end{cases}\label{ekv:dpp-for-R}
\end{align}
where $\bar U^{\eps,k}:=\cup_{k'=0}^k (\bar U^{\eps})^{k'}$, $(t^\eps_i,(b_1,\ldots,b_{k'})):=(t^\eps_i,b_j)_{j=1}^{k'}$ and
\begin{align*}
  \chi({t^\eps_i},X^{t,\vecv},(b_1,\ldots,b_{k'})):=\sum_{j=1}^{k'}\chi({t^\eps_i},X^{t,\vecv\otimes_{t^\eps_i}(t^\eps_i,b_1,\ldots,b_{j-1})},b_j)
\end{align*}
To see this note that
\begin{align*}
  R^{\vecv,k}_{t^\eps_i}&=\essinf_{u^S\in\mcU^{S,W,k,\eps}_{t^\eps_i}}\esssup_{\tau\in\mcT^{W,\eps}_{t^\eps_i}}\:(\ett_{[\tau={t^\eps_i}]}+ \ett_{[\tau>{t^\eps_i}]})J^{t,\vecv}_{t^\eps_i}(u^S(\tau),\tau)
  \\
  &=\essinf_{u^S\in\mcU^{S,W,k,\eps}_{t^\eps_i}}\big\{J^{t,\vecv}_{t^\eps_i}(u^S(t^\eps_i),t^\eps_i)\vee \esssup_{\tau\in\mcT^{W,\eps}_{t^\eps_{i+1}}}J^{t,\vecv}_{t^\eps_i}(u^S(\tau),\tau)\big\}
  \\
  &=\Psi(t^\eps_i,X^{t,\vecv}_{t^\eps_i}) \vee \essinf_{u^S\in\mcU^{S,W,k,\eps}_{t^\eps_i}} \esssup_{\tau\in\mcT^{W,\eps}_{t^\eps_{i+1}}} J^{t,\vecv}_{t^\eps_i}(u^S(\tau),\tau)
  \\
  &=\Psi(t^\eps_i,X^{t,\vecv}_{t^\eps_i}) \vee \inf_{\vecb\in\bar U^{\eps,k}}\{R^{\vecv\otimes_{t^\eps_i} (t^\eps_i,\vecb),k-|\vecb|,+}_{t^\eps_i}+\chi({t^\eps_i},X^{t,\vecv},\vecb)\}.
\end{align*}
Letting
\begin{align*}
  \tau^{S,\eps}_k(u):=\inf\{s\in \bbT^\eps\cap[t,T]:R^{\Xi^\eps(u),k-N^u_s}_{s}=\Psi(s,X^{t,\Xi^\eps(u)})\}
\end{align*}
and exploiting \eqref{ekv:dpp-for-R} in a standard fashion gives that
\begin{equation*}
  \underline Y^{k,\eps}_t=\essinf_{u\in\mcU^{W,k,\eps}_t}J_t(u_{\tau^\eps_k(u)-},\tau^\eps_k(u))\leq \bar Y^{k,\eps}_t,
\end{equation*}
where the last inequality holds since $\tau^{S,\eps}_k\in\mcT^{S,W,\eps}_t$.\qed\\

\begin{cor}
For each $k\in\bbN$ and $t\in[0,T]$, we have $\underline Y^k_t\leq \bar Y^k_t$, $\Prob$-a.s.
\end{cor}

\section{Game value by control randomization\label{sec:game-value}}

A successful approach to represent the solution to various types of control problems (including those with path-dependencies) has been to consider a weak formulation where the auxiliary probability space is endowed with an independent Poisson random measure that is used to represent the control. Optimization is then carried out by altering the probability measure to modify the compensator of the random measure, so that (in the limit and on an intuitive level) the path of the corresponding Poisson jump processes has probabilistic characteristics that mimic those of an optimal control.

This approach to stochastic optimal control is termed \emph{control randomization}~\cite{Fuhrman15} and as explained in the introduction, it is intimately connected to BSDEs with a constrained jumps. Despite its efficacy in addressing various types of optimal control problems, the approach pioneered in~\cite{Fuhrman15} has yet to be extended to encompass stochastic differential games. In the present section we bridge this void by establishing a connection between the previously defined lower and upper value functions and a non-linear Snell envelope.

In particular, we first introduce a randomized version of the game where we represent the impulse control by the sequence $(\sigma_j,\zeta_j)$ that appears in the Dirac sum formulation of the random measure, $\mu=\sum_{j\geq 1}\delta_{(\sigma_j,\zeta_j)}$, and then control the integrand in the Dol\'eans-Dade exponential appearing in a Girsanov transformation applied to the probability measure $\Prob$, effectively changing the probability distribution of $\mu$. Applying the same penalization routine as in Section~\ref{sec:Snell}, we show that the value function of this game corresponds to a non-linear Snell envelope defined over solutions to \eqref{ekv:bsde-c-jmp-lin}. We then proceed to show that the value of the randomized game coincides with both the upper and the lower value function of the original game posed in Section~\ref{sec:zsg}, thus proving Theorem~\ref{thm:sdg}.

\subsection{Randomized game and related non-linear Snell envelope}
We introduce the concept of a randomized game and establish its connection to the non-linear Snell envelope. We recall the definition of the set $\mcV$ (resp.~$\mcV^n$) as consisting of all $\Pred(\bbF)\otimes\mcB(U)$-measurable bounded maps $\nu=\nu_t(\omega,e):[0,T]\times\Omega\times U\to [0,\infty)$ (resp. $\to [0,n]$) and let $\mcV_{\inf>0}:=\{\nu\in\mcV:\inf\nu>0\}$ (resp.~$\mcV^n_{\inf>0}:=\{\nu\in\mcV^n:\inf\nu>0\}$). For each $t\in[0,T]$, $\tau\in\mcT_t$ and $\nu\in\mcV$ we define the $\mcF_t$-measurable random variable
\begin{align}
   J^{\mcR}_t(\nu,\tau)=\E^\nu\Big[\Psi(\tau,X)+\int_t^{\tau} f(r,X)dr+\int_t^{\tau}\!\!\!\int_U\chi(s-,X,e) \mu(ds,de)\Big|\mcF_t\Big],\label{ekv:obj-fun-dual}
\end{align}
where $\E^\nu$ is expectation with respect to the probability measure $\Prob^\nu$ on $(\Omega,\mcF)$ defined by $d\Prob^\nu:=\kappa^\nu_T d\Prob$ with
\begin{align*}
\kappa^{\nu}_s&:=\mcE_s\Big(\int_{0}^\cdot\!\!\int_U(\nu_r(e)-1)(\mu(dr,de)-\lambda(de))dr\Big)
\\
&:=\exp\Big(\int_{0}^s\!\!\!\int_U(1-\nu_r(e))\lambda(de)dr\Big)\prod_{\sigma_j\leq s}\nu_{\sigma_j}(\zeta_j).
\end{align*}

\begin{rem}\label{rem:XisEnu-bnd}
As in Lemma~\ref{lem:SFDEmoment}, it can be easily deduced that for each $p\geq 0$, there is a constant $C>0$ such that
\begin{align*}
  \E^\nu\big[\|X^u\|^p_T\big|\mcF_t\big]\leq C(1+\|X^u\|_t^p)
\end{align*}
whenever $t\in[0,T]$, $u\in\mcU$ and $\nu\in\mcV$.
\end{rem}

To establish a relation between our non-linear Snell envelope and the value of a game formulated over randomized controls, we adopt the approximation routine described in Section~\ref{sec:Snell}. Specifically, we consider the unique solution $(Y^n,Z^n,V^n,K^{+,n})\in \mcS^2\times\mcH^2(W)\times\mcH^2(\mu)\times\mcA^2$ to the penalized reflected BSDE
\begin{align}\label{ekv:rbsde-pen}
  \begin{cases}
    Y^n_t=\Psi(T,X)+\int_t^T f(s,X)ds-\int_t^T Z^n_s dW_s-\int_t^T\!\!\!\int_U V^n_{s}(e)\mu(ds,de)+K^{+,n}_T-K^{+,n}_t
    \\
    \quad-n\int_t^T\!\!\!\int_U(V^n_s(e)+\chi(s,X,e))^-\lambda(de)ds,\quad\forall t\in[0,T]\\
    Y^n_t\geq \Psi(t,X),\, \forall t\in [0,T]\mbox{ and } \int_0^T \left(Y^n_t-\Psi(t,X)\right)dK^{n}_t=0
  \end{cases}
\end{align}
and define $K^{-,n}_t:=n\int_0^t\!\int_U(V^n_s(e)+\chi(s,X,e))^-\lambda(de)ds$.

The following representation holds:
\begin{prop}\label{prop:pen-is-imp}
For each $n\in\bbN$ and $t\in[0,T]$, we have the representations
\begin{align}
Y^n_t&=\essinf_{\nu\in\mcV^n_{\inf>0}}\esssup_{\tau\in\mcT_t}J^{\mcR}_t(\nu,\tau)=\esssup_{\tau\in\mcT_t}\essinf_{\nu\in\mcV^n_{\inf>0}} J^{\mcR}_t(\nu,\tau)\label{ekv:Yn-repr-lin}
\end{align}
and
\begin{align}
Y^n_t=\essinf_{\nu\in\mcV^n_{\inf>0}}J^{\mcR}_t(\nu,\tau_n(t)),\label{ekv:Yn-repr-2}
\end{align}
where $\tau_n(t):=\inf\{s\geq t:Y^n_s=\Psi(s,X)\}\wedge T$.
\end{prop}

\noindent\emph{Proof.} We remark that $\Psi(\cdot,X)$ is \cadlag and upper left semi-continuous at predictable stopping times, whereby Remark~\ref{rem:rbsde-pen-n-has-comparison} implies that the condition for Theorem~\ref{thm:Quenez-Sulem-rbsde} holds, ensuring that $\tau_n(t)$ is an optimal stopping time for $Y^n_t$. We let
\begin{align*}
f^\nu(t,x,v):=f(t,x)+\int_U(v(e)+\chi(t,x,e))\nu_t(e)\lambda(de)
\end{align*}
and suppose that for each $(\tau,\nu)\in\mcT\times\mcV^n$, the triple $(Y^{\tau,\nu},Z^{\tau,\nu},V^{\tau,\nu})\in \mcS^2\times\mcH^2(W)\times\mcH^2(\mu)$ is the unique solution to
\begin{align*}
    Y^{\tau,\nu}_t=\Psi(\tau,X)+\int_t^\tau f^{\tau,\nu}(s,X,V^{\tau,\nu}_s)ds-\int_t^\tau Z^{\tau,\nu}_s dW_s-\int_t^\tau\!\!\!\int_U V^{\tau,\nu}_{s}(e)\mu(ds,de),\quad\forall t\in[0,\tau]
\end{align*}
then
\begin{align*}
    Y^{\tau,\nu}_t&=\Psi(\tau,X)+\int_t^\tau f(s,X)ds-\int_t^\tau Z^{\tau,\nu}_s dW_s-\int_t^\tau\!\!\!\int_U V^{\tau,\nu}_{s}(e)\mu(ds,de)
    \\
    &\quad+\int_t^\tau\!\!\!\int_U(V^{\tau,\nu}_s(e)+\chi(s,X,e))\nu_s(e)\lambda(de)ds
    \\
    &=\Psi(\tau,X)+\int_t^\tau f(s,X)ds-\int_t^\tau Z^{\tau,\nu}_s dW_s+\int_t^\tau\!\!\!\int_U \chi(s-,X,e)\mu(ds,de)
    \\
    &\quad-\int_t^\tau\!\!\!\int_U(V^{\tau,\nu}_{s}(e)+\chi(s-,X,e))\tilde\mu^{\nu}(ds,de),
\end{align*}
where $\tilde\mu^\nu(ds,de):=\mu(ds,de)-\nu_s(e)\lambda(de)ds$. Now, under the measure $\Prob^\nu$ defined above,  the compensator of $\mu$ is $\nu_\cdot(e)\lambda(de)$. Hence,
\begin{align*}
Y^{\tau,\nu}_t&=\E^\nu\Big[\Psi(\tau,X)+\int_t^{\tau} f(r,X)dr + \int_t^{\tau}\!\!\!\int_U\chi(s-,X,e) \mu(ds,de)\Big|\mcF_t\Big]
\\
&=\essinf_{\tau\in\mcT_t}J^{\mcR}_t(\nu,\tau).
\end{align*}
Finally, letting $v^*_t(e):=n\ett_{[V^n_t(e)<-\chi(t,X,e)]}$ and arguing as in the proof of Lemma~\ref{lem:Yn-repr} gives that $(\nu^*,\tau_n(t))$ is a saddle-point for the game, \ie for any $(\nu,\tau)\in\mcV^n\times\mcT_t$ we have
\begin{align*}
  J^{\mcR}_t(\nu^*,\tau)\leq J^{\mcR}_t(\nu^*,\tau_n(t))\leq J^{\mcR}_t(\nu,\tau_n(t)).
\end{align*}
In particular, as the restriction to positive densities becomes irrelevant when taking the infimum, the representations \eqref{ekv:Yn-repr-lin} and \eqref{ekv:Yn-repr-2} hold.\qed\\

By combining the above proposition with Theorem~\ref{thm:nl-Snell}, we obtain the following corollary:
\begin{cor}
There is a $Y\in\mcS^2$ such that for any $\eta\in\mcT$, we have $Y_\eta=\esssup_{\tau\in\mcT_\eta} Y^{\tau}_\eta$, where for each $\tau\in\mcT$, the quadruple $(Y^\tau,Z^\tau,V^\tau,K^{-,\tau})\in \mcS^2\times\mcH^2(W)\times\mcH^2(\mu)\times\mcA^2$ is the unique maximal solution to \eqref{ekv:bsde-c-jmp-lin}. Moreover, we have the representation
\begin{align*}
  Y_t=\essinf_{\nu\in\mcV_{\inf>0}}\esssup_{\tau\in\mcT_t}J^{\mcR}_t(\nu,\tau)
\end{align*}
that holds for all $t\in [0,T]$.
\end{cor}

To complete the proof of Theorem~\ref{thm:sdg}, we need to relate the value of the randomized game to that of the original game. This is accomplished in the subsequent subsections, where we begin by demonstrating that the expected value of the upper value function is dominated by that of $Y$.

\subsection{Proving that $\bar Y_t\leq Y_t$\label{subsec:barY-leq-Y}}
We begin by examining the inequality $\E\big[\bar Y_t-Y_t\big]\leq 0$, which can be readily deduced from the findings reported in Section 4.1 of \cite{Bandini18}. This relation will suffice since the principal result of the next subsection below implies that $\bar Y_t\geq Y_t$, $\Prob$-a.s., leading us to conclude that $\bar Y_t = Y_t$, $\Prob$-a.s. To maintain formality, we state the result in the following proposition:
\begin{prop}\label{prop:rand-geq}
For each $t\in [0,T]$, we have $\E[\bar Y_t-Y_t]\leq 0$.
\end{prop}

The novel work in \cite{Fuhrman15} considered a weak formulation of the control problem where, in addition to a supremum over controls, the value function was obtained by taking the supremum over all conceivable probability spaces. In particular, this made it straightforward to prove that the value in the original problem dominates that of the randomized version. An essential contribution made in \cite{Bandini18} was to consider a strong version of the control problem, where the probability space is fixed. Considering the type of zero-sum games that we analyze does not lead to a significant increase of the complexity compared to the analysis in \cite{Bandini18}. To see this, note that for each $\eps>0$, there is a stopping strategy $\tau^{S,\eps}\in\mcT^{S}_t$ (the set of non-anticipative maps $\tau^S:\mcU_t\to\mcT_t$), such that
\begin{align*}
  \E[\bar Y_t]&\leq \inf_{u\in \mcU^{W}_t}\E[J_t(u,\tau^{S,\eps}(u))]+\eps
  \\
  &=\inf_{u\in \mcU^{W}_t}\E\Big[\Psi(\tau^{S,\eps}(u),X^{t,u})+\int_t^T\ett_{[s\leq \tau^{S,\eps}(u)]}f(s,X^{t,u})ds+\sum_{j\geq 1}\ett_{[\eta_j\leq \tau^{S,\eps}(u)]}\gamma(\eta_j,X^{t,[u]_{j-1}},\beta_j)\Big]+\eps.
\end{align*}
For each $\tau^{S,\eps}\in\mcT^{S}_t$, the expression on the right-hand side represents an impulse control problem commencing at time $t$, built upon the historical trajectory $(X_s:0\leq t\leq s)$. Notably, this problem deviates from the standard archetype, as the pertinent information $(\Psi(\tau^{S,\eps}(u),\cdot),\ett{[s\leq \tau^{S,\eps}(u)]}f,\ett_{[\eta_j\leq \tau^{S,\eps}(u)]}\gamma)$ does not conform to conventional regularity assumptions.

Conversely, the pivotal outcomes leading up to Proposition 4.2 in \cite{Bandini18} hinge exclusively upon the properties of measurability, without necessitating any further regularity constraints on this dataset. Consequently, we are able to replicate the reasoning delineated in the corresponding proofs, culminating in the deduction that for each $n\in\bbN$, we have
\begin{align*}
  \inf_{u\in \mcU^{W}_t}\E[J_t(u,\tau^{S,\eps}(u))] =  \essinf_{\nu\in\mcV^n_{t,\inf>0}}\inf_{u\in \mcU_t}\E^\nu\Big[\Psi(\tau^{S,\eps}(u),X^{t,u})+\int_t^T\ett_{[s\leq \tau^{S,\eps}(u)]}f(s,X^{t,u})ds
  \\
  \quad+\sum_{j\geq 1}\ett_{[\eta_j\leq \tau^{S,\eps}(u)]}\gamma(\eta_j,X^{t,[u]_{j-1}},\beta_j)\Big],
\end{align*}
where $\mcV^n_{t,\inf>0}:=\{\nu\in\mcV^n_{\inf>0}:\nu_s\equiv 1,\forall s\in[0,t]\}$. Now, as $u^{t,\mu}:=(\sigma_{j+N^\mu_{t}},\zeta_{j+N^\mu_{t}})_{j\geq 1}\in\mcU_t$, where $N^\mu_t:=\mu((0,t],U)$, and $\tau:=\tau^{S,\eps}(u^{t,\mu})\in\mcT_t$, it follows by the second representation of $Y^n$ in \eqref{ekv:Yn-repr-lin} that the right-hand side is dominated by $\E[Y^n_t]$. Taking the limit as $n\to\infty$ and using dominated convergence, we conclude that Proposition~\ref{prop:rand-geq} holds.

\subsection{Proving that $Y_t\leq \underline Y_t$}
We finish the proof of Theorem~\ref{thm:sdg} by showing the following:

\begin{prop}\label{prop:unif-conv}
$Y_t\leq \underline Y_t$, $\Prob$-a.s.~for each $t\in [0,T]$.
\end{prop}

Effectively this proposition together with the preceding one implies that $Y_t\leq \underline Y_t\leq \bar Y_t\leq Y^n_t$, $\Prob$-a.s., for each $t\in [0,T]$ and $n\in\bbN$, enabling us to show that the game has a value.

The proof of this proposition is more involved than that of Proposition~\ref{prop:rand-geq} and is distributed over several lemmata. Leveraging Lemma~\ref{lem:underY-eps} and (H.\ref{hyp:unif-conv}), it suffices to prove that $\E\big[(Y_t-\underline Y^{k,\eps}_t)^+\big]\to 0$ as $\eps\to 0$ for each $k\in\bbN$. To achieve this, we use the fact that each $u\in \mcU^{k,\eps}$ can be approximated to arbitrary precision by a random measure in the spirit of Section 4.1.2 in \cite{Fuhrman15} or Section 4.2 in \cite{Fuhrman2020}. However, it should be noted that the game framework that we examine requires a different approach to the aforementioned works. Mainly, this is due to the fact that the optimal stopping times $\tau_n(t)$ depend on $n$. To resolve this issue we resort to a discretization of the set of stopping times $\mcT_t$ in the game representation of $Y^n$, both by restricting all stopping times to take values in $\bbT^\eps$ and by restricting the information that is used through only considering stopping times in a smaller filtration. It is worth noting that our approach has been specifically tailored to address both the game setting and the conditional framework, distinguishing it from the methods employed in \cite{Fuhrman15,Fuhrman2020}.

To begin, we restrict the stopping set for the randomized version of the game and impose a similar restriction on the number of interventions in the randomized control as we did in the original version. Specifically, we define $N^\mu_{s,t}:=\mu((s,t],U)$, so that $N^\mu_{s,t}$ represents the number of interventions in the control corresponding to $\mu$ within the interval $(s,t]$. Since we are only considering bounded-from-below $\nu$, we cannot place a $\Prob^\nu$-a.s. upper bound on $N^\mu_{t,T}$ in the optimization. Instead, we introduce the set:
\begin{align*}
  \mcV^{n,k,t}_{\inf>0}:=\{\nu\in\mcV^{n}_{\inf>0}:\nu_s(\omega,e)\leq 1,\forall (s,\omega,e)\in \{(r,\omega')\in[t,T]\times\Omega:\mu(\omega',(t,r),U)\geq k\}\times U\},
\end{align*}
so that for any $\nu\in\mcV^{n,k,t}_{\inf>0}$, we have $\E^\nu\big[N^\mu_{t,T}\big]\leq k+T\lambda(U)$ and more importantly there is a $C>0$ such that
\begin{align*}
  \E^\nu\big[(N^\mu_{t,T})^4\big]\leq C(1+k^4) 
\end{align*}
for all $(n,k,t)\in\bbN^2\times [0,T]$ and $\nu\in\mcV^{n,k,t}_{\inf>0}$. We now let
\begin{align}
Y^{n,k,\eps}_t&:=\esssup_{\tau\in\mcT^{\eps}_t}\essinf_{\nu\in\mcV^{n,k,t}_{\inf>0}}J^{\mcR}_t(\nu,\tau),\label{ekv:Yn-eps}
\end{align}
where $\mcT^{\eps}_t$ is the set of stopping times with respect to the filtration $\bbF^{t,W}\vee\bbF^{\Xi^\eps(\mu)}$ valued in $\bbT^\eps\cap[t,T]$. Here, $\Xi^\eps(\mu)=\Xi^\eps((\sigma_j,\zeta_j)^{N^\mu_T}_{j= 1})$ and for each $u\in\mcU$, the filtration $\bbF^{u}=(\mcF^{u}_s)_{s\in[0,T]}$ is defined as $\mcF^{u}_s:=\sigma([u]_s)$. To support this approximation we introduce the following objects. We let
\begin{align*}
\mu^{t,\eps}=(\sigma^{t,\eps}_j,\zeta^{t,\eps}_j)^{N^\mu_{T}}_{j= 1}:=(\sigma_j,\zeta_j)^{N^\mu_{t}}_{j= 1}\otimes_t\Xi^\eps\big((\sigma_j,\zeta_j)_{j=N^\mu_{t}+1}^{N^\mu_{T}}\big)
\end{align*}
and introduce the process $X^{t,\tau,\eps}:=\lim_{j\to\infty}X^{t,\tau,\eps,j}_{\cdot\wedge\tau}$, where the sequence $(X^{t,\tau,\eps,j})_{j\in\bbN}$ is defined recursively by letting $X^{t,\tau,\eps,j}$ be the unique \cadlag process that satisfies
\begin{align*}
X^{t,\tau,\eps,j}_s&=x_0+\int_0^s a(r,X^{t,\tau,\eps,j})dr+\int_0^s \sigma(r,X^{t,\tau,\eps,j})dW_r+\sum_{i= 1}^j\ett_{[\sigma_i\leq s]}\gamma(\sigma^{t,\eps}_i\wedge \tau,X^{t,\tau,\eps,i-1},\zeta^{t,\eps}_i),\quad\forall s\in [0,\tau]
\end{align*}
for each $\tau\in\mcT_t$. We adopt the notation $X^{t,\eps}:=X^{t,T,\eps}$ and let
\begin{align*}
   J^{\mcR,t,\eps}_s(\nu,\tau):=\E^\nu\Big[\Psi(\tau,X^{t,\eps})+\int_s^{\tau} f(r,X^{t,\eps})dr+\sum_{j\geq 1}\ett_{[s< \sigma_j\leq\tau]}\chi(\sigma_j,X^{t,\sigma_j,\eps,j-1},\zeta^{t,\eps}_j)\Big|\mcF_s\Big].
\end{align*}
With this definition, $J^{\mcR,t,\eps}$ corresponds to delaying the jumps in the state $X^{t,\eps}$ so that for each $s\in\bbT^\eps$, the new state $X^{t,\eps}_{\cdot\wedge s}$ is $\mcF^{W}_s\vee\mcF^{\mu^{t,\eps}}_s$-measurable while additionally discretizing the interventions to take values in $\bar U^\eps$. This discretization of the randomized impulse control allows us to solve the corresponding optimal stopping problem in a straightforward manner as shown in the next lemma.

\begin{lem}\label{lem:stop-proc-disc}
For each $\eps>0$ and $t\in[0,T]$, there is a non-increasing sequence of stopping times $(\tau_n^{\eps})_{n\in\bbN}\in\mcT^{\eps}_{t}$, such that
\begin{align}\label{ekv:stop-proc-disc}
\esssup_{\tau\in\mcT_t}\essinf_{\nu\in\mcV^{n,k,t}_{\inf>0}}J^{\mcR,t,\eps}_t(\nu,\xi^\eps_1(\tau)) = \essinf_{\nu\in\mcV^{n,k,t}_{\inf>0}}J^{\mcR,t,\eps}_t(\nu,\tau^\eps_n),
\end{align}
$\Prob$-a.s., for all $(k,n)\in\bbN^2$.
\end{lem}

\noindent\emph{Proof.} To show the existence of an optimal stopping time we use dynamic programming and introduce the processes
\begin{align*}
\tilde Y^{n,k,t,\eps}_s&:=\esssup_{\tau\in\mcT_s}\essinf_{\nu\in\mcV^{n,k,t}_{\inf>0}}J^{\mcR,t,\eps}_s(\nu,\xi_1^\eps(\tau)),
\end{align*}
from which we extract the $\bbF$-stopping times
\begin{align*}
\tau^\eps_n&:=\inf\{s\in \bbT^\eps\cap [t,T]:\tilde Y^{n,k,t,\eps}_s=\Psi(s,X^{t,\eps})\}.
\end{align*}
Since $\tilde Y^{n,k,t,\eps}_s$ is non-increasing in $n$, the sequence of stopping times $(\tau_n^{\eps})_{n\in\bbN}\in\mcT^{\eps}_{t}$ is non-increasing in $n$. Now, in the left-hand side of \eqref{ekv:stop-proc-disc}, stopping outside of the set $\bbT^\eps$ is suboptimal from the point of view of the maximizer. On $\bbT^\eps$, we thus have
\begin{align*}
\tilde Y^{n,k,t,\eps}_{t^\eps_i}&=\esssup_{\tau\in\mcT_{t^\eps_i}}\essinf_{\nu\in\mcV^{n,k,t}_{\inf>0}} \{\ett_{[\tau={t^\eps_i}]}J^{\mcR,t,\eps}_{t^\eps_i}(\nu,{t^\eps_i}) + \ett_{[\tau>{t^\eps_i}]}J^{\mcR,t,\eps}_{t^\eps_i}(\nu,\xi^\eps_1(\tau))\}
\\
&=\Psi({t^\eps_i},X^{t,\eps})\vee \esssup_{\tau\in\mcT_{t^\eps_{i+1}}}\essinf_{\nu\in\mcV^{n,k,t}_{\inf>0}}\E^\nu\Big[\int_{t^\eps_{i}}^{t^\eps_{i+1}} f(r,X^{t,\eps})dr
\\
&\quad+\sum_{j\geq 1}\ett_{[\sigma^{t,\eps}_j=t^\eps_{i+1}]}\chi(\sigma_j,X^{t,\sigma_j,\eps,j-1},\zeta^{t,\eps}_j) +J^{\mcR,t,\eps}_{t^\eps_{i+1}}(\nu,\xi^\eps_1(\tau))\Big|\mcF_{t^\eps_i}\Big].
\end{align*}
By a regular BSDE argument, the non-linear expectation $\essinf_{\nu\in\mcV^{n,k,t}_{\inf>0}}\E^\nu$ satisfies a tower property and we get that for arbitrary $\tau\in\mcT_{t^\eps_i}$,
\begin{align*}
&\essinf_{\nu\in\mcV^{n,k,t}_{\inf>0}}\E^\nu\Big[\int_{t^\eps_{i}}^{t^\eps_{i+1}} f(r,X^{t,T,\eps})dr
+\sum_{j\geq 1}\ett_{[\sigma^{t,\eps}_j=t^\eps_{i+1}]}\chi(\sigma_j,X^{t,\sigma_j,\eps,j-1},\zeta^{t,\eps}_j) + J^{\mcR,t,\eps}_{t^\eps_{i+1}}(\nu,\xi^\eps_1(\tau))\Big|\mcF_{t^\eps_i}\Big]
\\
&=\essinf_{\nu\in\mcV^{n,k,t}_{\inf>0}}\E^\nu\Big[\int_{t^\eps_{i}}^{t^\eps_{i+1}} f(r,X^{t,\eps})dr
+\sum_{j\geq 1}\ett_{[\sigma^{t,\eps}_j=t^\eps_{i+1}]}\chi(\sigma_j,X^{t,\sigma_j,\eps,j-1},\zeta^{t,\eps}_j) + \essinf_{\nu'\in\mcV^{n,k,t}_{\inf>0}}J^{\mcR,t,\eps}_{t^\eps_{i+1}}(\nu',\xi^\eps_1(\tau))\Big|\mcF_{t^\eps_i}\Big].
\end{align*}
This leads us to the conclusion that $\tilde Y^{n,k,t,\eps}$ satisfies the weak dynamic programming principle,
\begin{align*}
\tilde Y^{n,k,t,\eps}_{t^\eps_i}&\leq \Psi({t^\eps_i},X^{t,\eps})\vee \essinf_{\nu\in\mcV^{n,k,t}_{\inf>0}}\E^\nu\Big[\int_{t^\eps_{i}}^{t^\eps_{i+1}} f(r,X^{t,\eps})dr
\\
&\quad+\sum_{j\geq 1}\ett_{[\sigma^{t,\eps}_j=t^\eps_{i+1}]}\chi(\sigma_j,X^{t,\sigma_j,\eps,j-1},\zeta^{t,\eps}_j)+\tilde Y^{n,k,t,\eps}_{t^\eps_{i+1}}\Big|\mcF_{t^\eps_i}\Big].
\end{align*}
On the other hand, by iteration and again using the tower property, we find that the right-hand side is bounded from above by $\essinf_{\nu\in\mcV^{n,k,t}_{\inf>0}} J^{\mcR,t,\eps}(\nu,\tau^\eps_n)$. Standard arguments then give that $\tilde Y^{n,k,t,\eps}_s$ is $\mcF^{W}_s\vee\mcF^{\mu^{t,\eps}}_s$-measurable whenever $s\in \bbT^\eps$, implying that $\tau^\eps_n$ is an $\bbF^{W}\vee\bbF^{\mu^{t,\eps}}$-stopping time and thus belongs to $\mcT^{\eps}_{t}$ for each $n\in\bbN$. This proves the assertion whenever $t\in\bbT^\eps$. The generalization to arbitrary $t\in [0,T]$ is straightforward.\qed\\

In addition, we can compare the approximation to the original value process as follows.

\begin{lem}\label{lem:stop-disc}
For each $t\in[0,T]$, we have
\begin{align}\label{ekv:stop-disc}
 \lim_{\eps\to 0}\sup_{n\in\bbN}\E\big[(Y^n_t-\essinf_{\nu\in\mcV^{n,k,t}_{\inf>0}}J^{\mcR}_t(\nu,\tau^{\eps}_n))^+\big] = 0
\end{align}
for all $k\in\bbN$.
\end{lem}

\noindent\emph{Proof.} First, the representation \eqref{ekv:Yn-repr-lin} gives
\begin{align*}
  Y^n_t-\essinf_{\nu\in\mcV^{n,k,t}_{\inf>0}}J^{\mcR}_t(\nu,\tau^{\eps}_n)=\esssup_{\tau\in\mcT_t}\essinf_{\nu\in\mcV^n_{\inf>0}} J^{\mcR}_t(\nu,\tau)-\essinf_{\nu\in\mcV^{n,k,t}_{\inf>0}}J^{\mcR}_t(\nu,\tau^{\eps}_n).
\end{align*}
On the other hand, Lemma~\ref{lem:stop-proc-disc} implies that
\begin{align*}
  \essinf_{\nu\in\mcV^{n,k,t}_{\inf>0}}J^{\mcR}_t(\nu,\tau^{\eps}_n)&=  \essinf_{\nu\in\mcV^{n,k,t}_{\inf>0}}\big(J^{\mcR,t,\eps}_t(\nu,\tau^{\eps}_n) + (J^{\mcR}_t(\nu,\tau^{\eps}_n)-J^{\mcR,t,\eps}_t(\nu,\tau^{\eps}_n))\big)
  \\
  &\geq \esssup_{\tau\in\mcT_t}\essinf_{\nu\in\mcV^{n,k,t}_{\inf>0}}J^{\mcR,t,\eps}_t(\nu,\xi_1^\eps(\tau)) + \essinf_{\nu\in\mcV^{n,k,t}_{\inf>0}}\big(J^{\mcR}_t(\nu,\tau^{\eps}_n)-J^{\mcR,t,\eps}_t(\nu,\tau^{\eps}_n)\big)
  \\
  &\geq \esssup_{\tau\in\mcT_t}\essinf_{\nu\in\mcV^{n,k,t}_{\inf>0}}J^{\mcR}_t(\nu,\xi_1^\eps(\tau)) - 2\esssup_{\tau\in\mcT_t}\esssup_{\nu\in\mcV^{n,k,t}_{\inf>0}}|J^{\mcR,t,\eps}_t(\nu,\xi_1^\eps(\tau))-J^{\mcR}_t(\nu,\xi^\eps_1(\tau))|.
\end{align*}
Combined, this gives that
\begin{align}\nonumber
   Y^n_t-\essinf_{\nu\in\mcV^{n,k,t}_{\inf>0}}J^{\mcR}_t(\nu,\tau^{\eps}_n) &\leq \esssup_{\tau\in\mcT_t}\esssup_{\nu\in\mcV^{n,k,t}_{\inf>0}} \big(J^{\mcR}_t(\nu,\tau)- J^{\mcR}_t(\nu,\xi_1^\eps(\tau))\big)
   \\
   &\quad + 2\esssup_{\tau\in\mcT_t}\esssup_{\nu\in\mcV^{n,k,t}_{\inf>0}} | J^{\mcR,t,\eps}_t(\nu,\xi^\eps_1(\tau)) - J^{\mcR}_t(\nu,\xi^\eps_1(\tau))|.\label{ekv:2-times}
\end{align}
Concerning the first term, we have for fixed $k\in\bbN$ and any $\nu$ that belongs to $\mcV^{n,k,t}_{\inf>0}$ for some $n\in\bbN$ and any $\tau\in\mcT_t$, that
\begin{align*}
  J^{\mcR}_t(\nu,\tau)- J^{\mcR}_t(\nu,\xi_1^\eps(\tau))&=\E^\nu\Big[\Psi(\tau,X)-\Psi(\xi_1^\eps(\tau),X)-\int_{\tau}^{\xi_1^\eps(\tau)} f(r,X)dr
  \\
  &\quad-\int_{\tau}^{\xi_1^\eps(\tau)}\!\!\!\int_U\chi(s-,X,e) \mu(ds,de)\Big|\mcF_t\Big].
\end{align*}
For $j\in\bbN$, we let $X^j$ solve the SDE
\begin{align*}
X^{j}_s&=x_0+\int_0^s a(r,X^{j})dr+\int_0^s \sigma(r,X^{j})dW_r+\sum_{i= 1}^j\ett_{[\sigma_i\leq s]}\gamma(\sigma_i,X^{i-1},\zeta_i)
\end{align*}
and find that the process $(\chi(s,X^{j-1},b):s\in[\sigma_{j-1},T])$ only jumps at predictable stopping times. Since the times $(\sigma_j)_{j=1}^\infty$ are totally inaccessible under $\Prob^\nu$ whenever $\nu$ is bounded, Equation \eqref{ekv:ass@stopping} in Assumption~\ref{ass:oncoeff}.(\ref{ass:@stopping}) gives that
\begin{align*}
&\Psi(\xi_1^\eps(\tau),X)+\int_{\tau}^{\xi_1^\eps(\tau)}\!\!\!\int_U\chi(s-,X,e) \mu(ds,de)
\\
&=\Psi(\xi_1^\eps(\tau),X)+\sum_{j=N^{\mu}_{\tau}+1}^{N^{\mu}_{\xi_1^\eps(\tau)}}\big(\Psi(\sigma_j,X)+\chi(\sigma_j ,X^{j-1},\zeta_j)-\Psi(\sigma_j,X)\big)
\\
&\geq \Psi(\xi_1^\eps(\tau),X)+\sum_{j=N^{\mu}_{\tau}+1}^{N^{\mu}_{\xi_1^\eps(\tau)}}\big(\Psi(\sigma_j,X^{j-1})-\Psi(\sigma_j,X)\big)
\\
&=\Psi(\tau,X)+\sum_{j=N^{\mu}_{\tau}}^{N^{\mu}_{\xi_1^\eps(\tau)}} \big(\Psi(\sigma_{j+1}\wedge\xi_1^\eps(\tau) ,X^{j})-\Psi(\sigma_j\vee\tau,X)\big),
\end{align*}
with $\sigma_0:=0$. Assuming w.l.o.g.~that $\nu\equiv 1$ on $[0,t]$ we get that for each $K>0$,
\begin{align*}
\E\Big[\big(J^{\mcR}_t(\nu,\tau)- J^{\mcR}_t(\nu,\xi_1^\eps(\tau))\big)^+\Big]&\leq \E^\nu\Big[\int_{\tau}^{\xi_1^\eps(\tau)} |f(r,X)|dr+\sum_{j=N^{\mu}_{\tau}}^{N^{\mu}_{\xi_1^\eps(\tau)}} \big(\Psi(\sigma_j\vee\tau,X)-\Psi(\sigma_{j+1}\wedge\xi_1^\eps(\tau) ,X^{j})\big)^+\Big]
\\
&\leq \E^\nu\Big[C(\eps+\ett_{[\|X\|_T+N^\mu_{t,T}\geq K]}N^\mu_{t,T})(1+\|X\|_T^q)
\\
&\quad+\sum_{j=N^{\mu}_{t}}^{N^{\mu}_{T}\wedge (N^{\mu}_{t}+K)} \varpi_K(\mathbf d[(\sigma_j\vee\tau,X^j),(\sigma_{j+1}\wedge\xi_1^\eps(\tau) ,X^{j})])\big)^+\Big].
\end{align*}
Now, $N^\mu_{t,T}$ has fourth order moment under $\Prob^\nu$ for $\nu\in\mcV^{n,k,t}_{\inf>0}$ that is uniformly bounded in $n$. Remark~\ref{rem:XisEnu-bnd} then gives that
\begin{align*}
   \E^\nu\big[C(\eps+\ett_{[\|X\|_T+N^\mu_{t,T}\geq K]}N^\mu_{t,T})(1+\|X\|_T^q)\big]&\leq \E^\nu\big[C(\eps+\frac{N^\mu_{t,T}(N^\mu_{t,T}+\|X\|_T)}{K})(1+\|X\|_T^q)\big]
   \\
   &\leq C(\eps+\frac{1}{K})\E^\nu\big[\|X\|_T^{2(q+1)}\big]^{1/2}
   \\
   &\leq C(\eps+\frac{1}{K}),
\end{align*}
while repeating the argument from Section~\ref{subsec:barY-leq-Y} gives
\begin{align}\nonumber
&\E^\nu\Big[\sum_{j=N^{\mu}_{t}}^{N^{\mu}_{T}\wedge (N^{\mu}_{t}+K)} \varpi_K(\mathbf d[(\sigma_j\vee\tau,X^j),(\sigma_{j+1}\wedge\xi_1^\eps(\tau) ,X^{j})])\Big]
\\
&\leq \sup_{u\in\mcU^K_t}\E\Big[\sum_{j= 0}^{N}\varpi_K(\mathbf d[(\eta_{j}\vee\tau,X^{t,[u]_j}),(\eta_{j+1}\wedge \xi^\eps_1(\tau),X^{t, [u]_j})])\Big]\label{ekv:nr-2}
\end{align}
Since the right-hand side of \eqref{ekv:nr-2} tends to 0 as $\eps\to 0$ uniformly in $(\nu,\tau)$ by Corollary~\ref{cor:rew-cont}, we conclude that
\begin{align*}
  \limsup_{\eps\to 0}\sup_{n\in\bbN}\sup_{\tau\in\mcT_t}\sup_{\nu\in\mcV^{n,k,t}_{\inf>0}}\E\big[\big(J^{\mcR}_t(\nu,\tau)- J^{\mcR}_t(\nu,\xi_1^\eps(\tau))\big)^+\big]\leq \frac{C}{K}
\end{align*}
implying that the left-hand side must equal zero as $K>0$ was arbitrary. We move on to the second term on the right-hand side of \eqref{ekv:2-times} and have for any $\nu\in\mcV^{n,k,t}_{\inf>0}$ and $\tau\in\mcT_t$, valued in $\bbT^\eps$, that
\begin{align*}
  &|J^{\mcR}_t(\nu,\tau)-J^{\mcR,t,\eps}_t(\nu,\tau)|\leq\E^\nu\Big[|\Psi(\tau,X)-\Psi(\tau,X^{t,\eps})|+\int_t^{\tau} |f(r,X)-f(r,X^{t,\eps})|dr
  \\
  &\quad +\sum_{j\geq 1}\ett_{[t\leq \sigma_j \leq \tau]}|\chi(\sigma_j,X^{j-1},\zeta_j)-\chi(\sigma_j,X^{t,\sigma_j,\eps,j-1},\zeta^{t,\eps}_j)|\Big|\mcF_t\Big].
\end{align*}
Repeating the latter part of the above argument now gives that for each $K\geq 0$, it holds that
\begin{align*}
  &\E\big[|J^{\mcR}_t(\nu,\tau)-J^{\mcR,t,\eps}_t(\nu,\tau)|\big]
  \\
&\leq C(\eps+\frac{1}{K}) + K\sup_{u,\tilde u\in\mcU^K_t,\,|u-\tilde u|\leq 2K\eps}\E\Big[\sum_{j= 0}^{N\wedge \tilde N}\varpi_K(\mathbf d[(\eta_{j}\vee\tilde\eta_{j},X^{t,[u]_{j-1}}),(\eta_j\vee\tilde\eta_j,X^{t,[\tilde u]_{j-1}})]+\eps)\Big]
\end{align*}
implying that
\begin{align*}
\limsup_{\eps\to 0}\sup_{n\in\bbN}\sup_{\nu\in\mcV^{n,k,t}_{\inf>0}}\sup_{\tau\in\mcT_t}\E\big[|J^{\mcR}_t(\nu,\xi^\eps_1(\tau))-J^{\mcR,t,\eps}_t(\nu,\xi^\eps_1(\tau))|\big]=0
\end{align*}
which proves the assertion.\qed\\

Inspired by Section 4.3 of \cite{Fuhrman2020}, we introduce an auxiliary probability space $( \Omega''',\mcF''',\Prob''')$ on which lives real-valued random variables $(U^m_j,S^m_j)_{m,j\in\bbN}$ and random measures $(\pi^l)_{l\geq 1}$ such that
\begin{enumerate}
  \item the $U^m_j$ are all uniformly distributed on $(0,1)$,
  \item the probability distribution of $S^m_j$ admits a density $f^m_j$ with respect to the Lebesgue measure, that has support on the interval $((1-2^{j})/m,(1-2^{j-1})/m)$, so that $0<S^m_1<S^m_2<\cdots<1/m$ for every $m\in\bbN$,
  \item every $\pi^l$ is a Poisson random measure on $(0,\infty)\times U$, with compensator $l^{-1}\lambda(da)dt$, with respect to its natural filtration;
  \item the random elements $U^m_j,S^{m'}_{j'},\pi^l$ are all independent.
\end{enumerate}
Now, we define $\hat\Omega:=\Omega\times\Omega'''$, let $\hat \mcF$ be the $\Prob\otimes\Prob'''$ completion of $\mcF\otimes\mcF'''$ and let $\hat\Prob$ denote the extension of $\Prob\otimes\Prob'''$ to $\hat\mcF$. Further, we let $\hat W,\hat\mu,\hat U^m_j,\hat S^{m'}_{j'}$ and $\hat \pi^l$ denote the canonical extensions of $W,\mu,U^m_j,S^{m'}_{j'}$ and $\pi^l$ to $\hat\Omega$. For $t\in [0,T]$, $u\in\check\mcU^{\hat W}_t$ and $\tau\in\check\mcT^{\hat W}_t$ (which are extensions of $\mcU^W_t$ and $\mcT^W_t$ to $\hat\Omega$, that are more carefully defined below), we define
\begin{align*}
  \hat J_t(u,\tau)=\hat\E\Big[\Psi(\tau,\hat X^{t,u})+\int_t^{\tau} f(r,\hat X^{t,u})dr+\sum_{j=1}^{N}\ett_{[\eta_ j\leq\tau]}\chi(\eta_j,\hat X^{t,[u]_{j-1}},\beta_j)\Big|\hat\mcF_t\Big],
\end{align*}
where $\hat\E$ is expectation with respect to $\hat\Prob$, the filtration $\hat\bbF:=(\hat \mcF_t)_{{t\geq 0}}$ is the $\hat\Prob$-augmented natural filtration on $(\hat\Omega,\hat\mcF)$ generated by $\hat W$ and $\hat\mu$ and $\hat X^{t,u}$ solves
\begin{align*}
\hat X^{t,u}_s&=x_0+\int_0^sa(r,\hat X^{t,u})dr+\int_0^s\sigma(r,\hat X^{t,u})d\hat W_r+\int_0^{s\wedge t}\!\!\!\int_U\gamma(r-,\hat X^{t,u},e)\hat \mu(dr,de)
\\
&\quad+\sum_{j=1}^N \ett_{[\eta_j\leq s]}\gamma(\eta_j,\hat X^{t,[u]_{j-1}},\beta_j).
\end{align*}
When extending the basic notations to the space $(\hat\Omega,\hat\mcF,\hat\Prob)$, we introduce two versions of most objects depending on whether they utilize the information in the $\sigma$-algebra $\mcF'''$ or not. We denote objects that incorporate the information in $\mcF'''$ with a check symbol, while objects that do not use this information are denoted with a hat symbol. Specifically, we make the following definitions:
\begin{itemize}
\item For each $t\in[0,T]$, we let $\hat\bbF^{t,\hat W}$ (resp.~$\check\bbF^{t,\hat W}$) be the $\hat\Prob$-completion of the filtration $(\mcF^{t,W}_s \times \Omega''')_{s\geq 0}$ (resp.~$(\mcF^{t,W}_s\otimes \mcF''')_{s\geq 0}$).
\item We let $\hat\mcT^{\hat W}_t$ be the set of all $\hat\bbF^{t,\hat W}$-stopping times $\tau$ with $\tau\in [t,T]$, $\hat\Prob$-a.s.
\item For each $t\in[0,T]$, $\hat\mcU^{\hat W}_t$ (resp. $\check\mcU^{\hat W}_t$) is the set of impulse controls $ u:=(\eta_j,\beta_j)_{j=1}^{ N}$, where $(\eta_j)_{j=1}^\infty$ is a non-decreasing sequence of $\hat\bbF^{t,\hat W}$-stopping times (resp.~$\check\bbF^{t,\hat W}$-stopping times) with $\eta_1\geq t$, $\beta_j$ a $U$-valued, $\hat\mcF^{t,\hat W}_{\eta_j}$-measurable (resp.~$\check\mcF^{t,\hat W}_{\eta_j}$-measurable) random variable and $N:=\max\{j:\eta_j\leq T\}$.
\item For each $t\in[0,T]$ and $k\in\bbN$, we denote by $\hat\mcU^{\hat W,k}_t$ (resp. $\check\mcU^{\hat W,k}_t$) the subset of $\hat\mcU^{\hat W}_t$ (resp. $\check\mcU^{\hat W}_t$) containing all impulse controls with $N\leq k$, $\hat\Prob$-a.s.
\item For each $\eps>0$, $\hat\mcU^{\hat W,k,\eps}_t$ is the subset of $\hat\mcU^{\hat W,k}_t$, with all impulse controls $\hat u:=(\hat\eta_j,\hat\beta_j)_{j=1}^{\hat N}$ for which $(\hat\eta_j,\hat\beta_j)$ is $\bbT^\eps\times\bar U^\eps$-valued for $j=1,\ldots,\hat N$.
\item For each $\check u=(\check\eta_j,\check\beta_j)_{j=1}^{\check N}\in\check\mcU^{\hat W}_t$, we denote by $\hat\mcT^{\hat W,\check u}_t$ the set of $\hat\bbF^{t,\hat W}\vee \hat\bbF^{\check u}$-stopping times, where $\hat\bbF^{\check u}:=(\hat\mcF^{\check u}_s)_{s\geq 0}$ is defined as $\hat\mcF^{\check u}_s:=\sigma\big((\check\eta_j,\check\beta_j)_{j=1}^{\check N(s)}\big)$, with $\check N(s):=\max\{j\geq 0:\check\eta_j\leq s\}$ for all $s\geq 0$, and we let $\hat\mcT^{\hat W,\check u,\eps}_t$ be the restriction to stopping times taking values in $\bbT^\eps$.
\end{itemize}

The idea is to use the sequences $\hat U^m_j$ and $\hat S^{m'}_{j'}$ to ``randomize'' an impulse control $\hat u\in\hat\mcU^{\hat W,k,\eps}_t$ and then add $\hat \pi^l$ to get a new sequence $\check u:=(\check\eta_j,\check\beta_j)_{j=1}^\infty$ of random variables such that the $\hat\Prob$-compensator of the corresponding random measure\footnote{As above, we abuse notation slightly by letting $\check u$ denote both the sequence of random variables and the corresponding random measure.} $\check u:=\sum_{j=1}^\infty\delta_{(\check\eta_j,\check\beta_j)}$ has a density $\check\nu$ with respect to $\lambda(da)dt$ which is strictly positive and such that $\check u$ is sufficiently ``close'' to $\hat u$.

Following the above procedure, we define $\hat{\underline Y}^{k,\eps}_t$ as the canonical extension of $\underline Y^{k,\eps}_t$ to $\hat\Omega$. Before we proceed to prove Proposition~\ref{prop:unif-conv}, we present the following lemma.

\begin{lem}\label{lem:will-conv}
For each $t\in [0,T]$, $k\in\bbN$ and $\varrho>0$, there is an $\eps\in (0,\varrho]$ and a $\check u\in \check\mcU^{\hat W}_{t}$ such that
\begin{align*}
 \hat \E\big[(\hat J_t(\check u,\tau) - \hat{\underline Y}_t^{k})^+\big]\leq \varrho
\end{align*}
for each $\tau\in\hat\mcT^{\hat W,\Xi^{\eps}(\check u),\eps}_t$ and the random measure on $(t,T]\times U$ corresponding to $\check u$ has a $\hat\Prob$-compensator with respect to the filtration $\hat\bbF^{t,\hat W}\vee \hat\bbF^{\check u}$ that is absolutely continuous with respect to $\lambda$ and takes the form
\begin{align*}
\check\nu_t(\hat\omega,a)\lambda(da)dt
\end{align*}
where $\check\nu:(t,T]\times \hat\Omega\times U\to [0,\infty)$ is $\Pred(\hat\bbF^{t,\hat W}\vee \hat\bbF^{\check u}) \otimes\mcB(U)$-measurable and bounded away from zero. Moreover, the number of interventions of the impulse control $\check u$ (\ie the $\check u$-measure of the set $(t,T]\times U$), denoted $\hat N^{\check u}_{t,T}$, has moments of all orders under $\hat\Prob$.
\end{lem}

\noindent\emph{Proof.} As noted in the proof of Lemma~\ref{lem:underY-eps}, non-anticipativity implies that any strategy $u^S\in\mcU^{S,W}_t$ can be written $u^S(\tau):=u_{\tau-}\otimes_\tau \tilde u^S(\tau)$, with $u\in\mcU^W_t$ and $\tilde u^S(\tau):=(\tilde \eta^S_j(\tau),\tilde \beta^S_j(\tau))_{j=1}^{\tilde N^S(\tau)}\in\mcU^S_{\tau}\cap \mcU^{S,W}_t$ and by Assumption~\ref{ass:oncoeff}.(\ref{ass:@stopping}) it is optimal to have $\tilde N^S\equiv 0$. Lemma~\ref{lem:underY-eps} now implies that in $L^1(\hat\Omega,\hat\mcF_t,\hat\Prob)$ we have
\begin{align*}
  \hat{\underline Y}_t^{k}&=\lim_{\eps\to 0}\essinf_{u\in\hat\mcU^{\hat W,k,\eps}_t}\esssup_{\tau\in\hat\mcT^{\hat W,\eps}_t}\hat J_t(u_{\tau-},\tau)=\lim_{\eps\to 0}\essinf_{u\in\hat\mcU^{\hat W,k,\eps}_t}\esssup_{\tau\in\hat\mcT^{\hat W,u,\eps}_t}\hat J_t(u_{\tau-},\tau).
\end{align*}
In particular, there is an $\eps\in (0,\varrho]$ such that
\begin{align}
  \hat\E\big[|\hat{\underline Y}_t^{k}-\essinf_{u\in\hat\mcU^{\hat W,k,\eps}_t}\esssup_{\tau\in\hat\mcT^{\hat W,u,\eps}_t}\hat J_t(u_{\tau-},\tau)|\big]\leq \varrho/4.\label{ekv:interm}
\end{align}
Using \eqref{ekv:interm}, we prove the lemma in two steps:\\

\textbf{Step 1:} We first show existence of a non-negative (not necessarily bounded away from zero) map $\nu$ satisfying the first part of the lemma. By the definition of the essential supremum and stability under pasting of the set $\hat\mcU^{\hat W,k,\eps}_{t}$, there is a $\hat u^{\eps,\varrho}:=(\hat \eta^{\eps,\varrho}_j,\hat \beta^{\eps,\varrho}_j)_{j=1}^{\hat N^{\eps,\varrho}}\in\hat\mcU^{\hat W,k,\eps}_t$ (with $\hat \eta^{\eps,\varrho}_{\hat N^{\eps,\varrho}}<T$) such that
\begin{align}\label{ekv:hat-u-eps-opt}
  \essinf_{u\in\hat\mcU^{\hat W,k,\eps}_{t}}\,\esssup_{\tau\in\hat\mcT^{\hat W,u,\eps}}\hat J_t(u_{\tau-},\tau)\geq \esssup_{\tau\in\hat\mcT^{\hat W,\hat u^{\eps,\varrho},\eps}}\hat J_t(\hat u^{\eps,\varrho}_{\tau-},\tau)-\varrho/4,
\end{align}
$\hat\Prob$-a.s. Define for each $m\in\bbN$, the transition kernel $q^m(b,da)$ on $U$ as in the proof of Lemma 4.4 of \cite{Fuhrman2020}, let
\begin{align*}
  \check \eta^m_j:=\hat\eta^{\eps,\varrho}_j+\hat S^m_j,\qquad \check\beta^m_j:=q^m(\hat\beta^{\eps,\varrho}_j,\hat U^m_j),\qquad \check N^m:=\inf\{j\geq 0:\check \eta^m_j<T\}
\end{align*}
and introduce the impulse control $\check u^m:=(\check \eta^m_j,\check\beta^m_j)_{j=1}^{\check N^m}\in\check\mcU^{\hat W,k}_t$. According to Lemma A.11 in~\cite{Fuhrman15} the corresponding $\hat\Prob$-compensator with respect to $\hat\bbF^{t,\hat W}\vee \hat\bbF^{\check u^m}$ is given by the explicit formula
\begin{align*}
  \sum_{j=1}^k\ett_{(\hat\eta^{\eps,\varrho}_j\vee \check \eta^m_{j-1},\check \eta^m_j]}(s)q^m(\hat\beta^{\eps,\varrho}_j,da)\frac{f^m_j(s-\hat\eta^{\eps,\varrho}_j)}{1-F^m_j(s-\hat\eta^{\eps,\varrho}_j)}ds,
\end{align*}
with $F^m_j(s):=\int_{-\infty}^s f^m_j(r)dr$. For each $m\in\bbN$, this compensator is clearly $\Pred(\hat\bbF^{t,\hat W}\vee \hat\bbF^{\check u^m}) \otimes\mcB(U)$-measurable.

Moreover, there is an $m'\in\bbN$ such that the densities of the random variables in the sequence $(S^{m}_j)_{j\in\bbN}$ all have support in $(0,\Delta t^\eps)$ and $\xi_2^\eps(\check\beta^m_j)=\hat\beta^{\eps,\varrho}_j$ for $j=1,\ldots,k$, whenever $m\geq m'$. We thus conclude that for any such $m$ and each $\tau\in\hat\mcT^{\hat W,\hat u^{\eps,\varrho},\eps}_t$, the impulse controls $\check u^m$ and $\hat u^{\eps,\varrho}_{\tau-}$ have the same number of interventions in the interval $[t,\tau]$, which gives
\begin{align*}
  & \hat J_t(\check u^m,\tau)- \hat J_t(\hat u^{\eps,\varrho}_{\tau-},\tau)\leq \hat\E\Big[\Psi(\tau,\hat X^{t,\check u^m})-\Psi(\tau,\hat X^{t,\hat u^{\eps,\varrho}_{\tau-}})+\int_t^{\tau} |f(r,X^{t,\check u^m})-f(r,X^{t,\hat u^{\eps,\varrho}})|dr
\\
&\quad+\sum_{j=1}^{k}\ett_{[\hat\eta^{\eps,\varrho}_ j<\tau]}\big(\chi(\check \eta^m_j ,\hat X^{t,[\check u^m]_{j-1}},\check\beta^m_j)-\chi(\hat\eta^{\eps,\varrho}_j ,\hat X^{t,[\hat u^{\eps,\varrho}]_{j-1}},\hat\beta^{\eps,\varrho}_j)\big)\Big|\mcF_t\Big].
\end{align*}
Arguing as in the proof of Lemma~\ref{lem:underY-eps} while appealing to a slightly adjusted version of Corollary~\ref{cor:rew-cont}, where we allow impulse controls in $\check \mcU^{\hat W}_t$, we now find that there is a $m''\geq m'$ such that
\begin{align}\label{ekv:m-approx}
  \hat\E\big[(\hat J_t(\check u^m,\tau)-\hat J_t(\hat u^{\eps,\varrho}_{\tau-},\tau))^+\big]\leq \varrho/4
\end{align}
for all $\tau\in\hat\mcT^{\hat W,\hat u^{\eps,\varrho},\eps}_t$, whenever $m\geq m''$. Now, whenever $m\geq m''$ we have $\Xi^\eps(\check u^m)=(\hat \eta^{\eps,\varrho}_{j}+\Delta t^\eps,\hat\beta^{\eps,\varrho}_j)_{j=1}^{\hat N^{\eps,\varrho}}$. In particular, this gives that $\hat\mcT^{\hat W,\Xi^\eps(\check u^m),\eps}_t\subset \hat\mcT^{\hat W,\hat u^{\eps,\varrho},\eps}_t$ and by \eqref{ekv:interm}-\eqref{ekv:m-approx} we get
\begin{align*}
  \hat\E\big[(\hat J_t(\check u^m,\tau) - \hat {\underline Y}_t^{k})^+\big]&\leq \hat\E\big[(\hat J_t(\check u^m,\tau)-\esssup_{\tau'\in\hat\mcT^{\hat W,\hat u^{\eps,\varrho},\eps}_t}\hat J_t(\hat u^{\eps,\varrho}_{\tau'-},\tau'))^+\big]+\varrho/2
  \\
  &\leq \sup_{\tau'\in\hat\mcT^{\hat W,\hat u^{\eps,\varrho},\eps}_t}\hat\E\big[(\hat J_t(\check u^m,\tau')-\hat J_t(\hat u^{\eps,\varrho}_{\tau'-},\tau'))^+\big]+\varrho/2
  \\
  &\leq 3\varrho/4
\end{align*}
for any $\tau\in \hat\mcT^{\hat W,\Xi^\eps(\check u^m),\eps}_t$.

\bigskip

\textbf{Step 2:} To establish the claim, we need to modify $\check u^m$ so that the corresponding density with respect to $\lambda$ is bounded away from $0$ on $[t,T]$. We therefore consider the control $\check u^{m,l}:=(\check\eta_j^{m,l},\check\beta^{m,l}_j)_{j=1}^{\check N^{m,l}}$ corresponding to the random measure $\check u^m+\hat \pi^l(\cdot\cap [t,T],\cdot)$ and note that the number of interventions of $\check u^{m,l}$ on $[t,T]$ is bounded by $k+\hat N^{\hat \pi^l}$, where $\hat N^{\hat \pi^l}_{t,T}:=\hat \pi^l([t,T],U)$ is Poisson distributed with parameter $\lambda(U)(T-t)/l$ under $\hat\Prob$.  In particular, this gives that $\hat N_{t,T}^{\hat \pi^l}$ and then also $k+\hat N_{t,T}^{\hat \pi^l}$ has moments of all orders under $\hat\Prob$.

We let $A:=\{\omega:\hat N_{t,T}^{\hat \pi^l}=0\}$ and get that for any $\tau\in\hat\mcT^{\hat W,\Xi^\eps(\check u^{m}),\eps}_t$,
\begin{align*}
  &\hat J_t(\check u^{m,l},\tau)-\hat J_t(\check u^{m},\tau)
  \\
  &\leq \hat\E\Big[\ett_{A^c}\big(\Psi(\tau,\hat X^{t,\check u^{m,l}})-\Psi(\tau,\hat X^{t,\check u^m})+\int_t^{\tau} (|f(r,X^{t,\check u^{m,l}}|+|f(r,X^{t,\check u^m}|)dr
  \\
  &\quad+\sum_{j=1}^{\infty}\ett_{[\check \eta^{m,l}_ j\leq\tau]}\chi(\check \eta^{m,l}_j ,\hat X^{t,[\check u^{m,l}]_{j-1}},\check\beta^{m,l}_j)-\sum_{j=1}^{k}\ett_{[\check\eta^m_ j\leq\tau]}\chi(\check \eta^{m}_j ,\hat X^{t,[\check u^m]_{j-1}},\check\beta^m_j)\big)\Big|\mcF_t\Big]
  \\
  &\leq C\hat\E\Big[\ett_{A^c}\check N^{m,l}\big(1+\|\hat X^{t,\check u^{m,l}}\|^q_T+\|\hat X^{t,\check u^{m}}\|^q_T\big)\Big|\mcF_t\Big]
  \\
  &\leq C\hat\E\Big[\hat N_{t,T}^{\hat \pi^l}\big(1+\|\hat X^{t,\check u^{m,l}}\|^q_T+\|\hat X^{t,\check u^{m}}\|^q_T\big)\Big|\mcF_t\Big],
\end{align*}
where the right-hand side tends to 0, $\hat\Prob$-a.s., when $l\to\infty$, since $\hat\E\big[(\hat N_{t,T}^{\hat \pi^l})^2\big|\mcF_t\big]\to 0$ as $l\to\infty$. In particular, there is a $l'\in\bbN$ such that
\begin{align*}
  \hat\E\big[(\hat J_t(\check u^{m,l},\tau)-\hat J_t(\check u^{m},\tau))^+\big]\leq \varrho/4
\end{align*}
for all $\tau\in\hat\mcT^{\hat W,\Xi^\eps(\check u^{m}),\eps}_t$, whenever $l\geq l'$. Since, for any $\tau\in\hat\mcT^{\Xi^\eps(\check u^{m,l}),\eps}_t$, there is a $\tau'\in\hat\mcT^{\Xi^\eps(\check u^{m}),\eps}_t$ such that $\ett_A\tau=\ett_A\tau'$, $\hat\Prob$-a.s., there is thus a $l'\in\bbN$ such that
\begin{align*}
  \hat\E\big[(\hat J_t(\check u^{m,l},\tau)- \hat{\underline Y}^{k,\eps}_t)^+\big]\leq \varrho
\end{align*}
for all $\tau\in\hat\mcT^{\Xi^\eps(\check u^{m,l}),\eps}_t$, whenever $m\geq m''$ and $l\geq l'$. The assertion now follows by noting that $\check u^{m,l}$ has a $\hat\Prob$-compensator with respect to the filtration $\hat\bbF^{t,\hat W}\vee \hat\bbF^{\check u^{m,l}}$ that has a density with respect to $\lambda$ which is bounded from below by $1/l$ on $[t,T]$.\qed\\

\begin{rem}\label{rem:check-nu-n}
We may extend $\check\nu$ in the statement of Lemma~\ref{lem:will-conv} to $[0,T]\times \hat\Omega\times U$ by setting $\check\nu_s\equiv 1$, whenever $s\in[0,t]$. Letting $\check \nu^n:=\check\nu\wedge n$ we then find that $\check \nu^n\in\check\mcV^{n,k,t}_{\inf>0}$. Here, analogously to above, $\mcV^{n,k,t}_{\inf>0}$ is the set of all $\Pred(\hat\bbF^{t,\hat W}\vee \hat\bbF^{\check u}) \otimes\mcB(U)$-measurable maps $\nu:[t,T]\times\hat\Omega\times U\to [0,n]$ with $\inf\nu>0$ and $\nu_s(\omega,e)\leq 1$ on the set $\{(s,\omega,e)\in [0,T]\times\hat\Omega\times U:N^{\check u}_{t,s}(\omega)\geq k\}$.
\end{rem}

\bigskip

\noindent\emph{Proof of Proposition~\ref{prop:unif-conv}.} We fix $t\in[0,T]$ and $k\in\bbN$ and note that for each $\varrho>0$, there is by Lemma~\ref{lem:stop-disc} a $\varrho'>0$ such that for all $n\in\bbN$, it holds that
\begin{align}\label{Yn-rel-varrho}
 \E\big[(Y^n_t- \essinf_{\nu\in\mcV^{n,k,t}_{\inf>0}}J^{\mcR}_t(\nu,\tau^\eps_n))^+\big]\leq \varrho
\end{align}
whenever $\eps\in (0,\varrho']$. Moreover, using Lemma~\ref{lem:will-conv}, with $\varrho\wedge\varrho'$ in place of $\varrho$, gives that there is an $\eps\in (0,\varrho\wedge\varrho']$ and a random measure $\check u$ such that
\begin{align}\label{ekv:relate-lem}
  \hat\E\big[(\hat J_t(\check u,\tau) - \hat {\underline Y}^{k}_t)^+\big]\leq \varrho
\end{align}
for each $\tau\in\hat\mcT^{\hat W,\Xi^{\eps}(\check u),\eps}_t$.

To establish a correspondence between the original game and its randomized version, we combine $\mu$ and $\check u$ to obtain the random measure $\check\mu:=\mu(\cdot\cup (0,t],\cdot)+\check u(\cdot\cup (t,T],\cdot)$. By Lemma~\ref{lem:will-conv}, $\check\mu$ has a $\hat\Prob$-compensator with respect to the filtration $\hat\bbF^{t,\hat W}\vee \hat\bbF^{\check u}$, and this compensator has a density $\ett_{[0,t]}+\ett_{(t,T]}\check\nu$ with respect to $\lambda$. We will continue our abuse of notation and use $\check \nu$ to denote this density, the infimum of which is strictly positive whereas the supremum may be unbounded.

The proof, which is divided into three steps, uses an auxiliary formulation of the randomized version of the game with a state process $\check X$ that is driven by $\hat W$ and $\check\mu$.\\

\textbf{Step 1:} We begin by defining the auxiliary randomized version of the game. Letting $(\check\sigma_j,\check\zeta_j)_{j\geq 1}$ be the impulse control corresponding to $\check\mu$ we find, since $\check\nu$ is bounded from below, that
\begin{align*}
\hat M_s:=\exp\Big(\int_{0}^s\!\int_U(1-(\check\nu_r(a))^{-1})\lambda(da)dr\Big)\prod_{\check\sigma_j\leq s}(\check\nu_{\check\sigma_j}(\check\zeta_j))^{-1}
\end{align*}
is a strictly positive martingale with respect to the filtration $\hat\bbF^{\hat W}\vee\hat\bbF^{\check\mu}$ under $\hat\Prob$. Furthermore, as $\check\nu_s\equiv 1$ for all $s\in [0,t]$, we have $\hat M\equiv 1$ on $[0,t]$. We define the equivalent probability measure $\check\Prob$ on $(\hat\Omega,\hat\mcF)$ as $d\check\Prob=\hat M_Td\hat\Prob$. By the Girsanov theorem, $\check\mu$ has $\check\Prob$-compensator $\lambda(da)ds$ with respect to the filtration $\hat\bbF^{\hat W}\vee \hat\bbF^{\check\mu}$. Moreover, despite the fact that $\check\nu$ is generally not bounded we still have a Dol{\'e}ans-Dade exponential
\begin{align*}
\hat\kappa^{\check\nu}_s:=\exp\Big(\int_{0}^s\!\int_U(1-\check\nu_r(e))\lambda(de)dr\Big)\prod_{\hat\sigma_j\leq s}\check\nu_{\hat\sigma_j}(\hat\zeta_j)
\end{align*}
for which $\check \E[\hat\kappa^{\check\nu}_T]=\hat \E[\hat M_T\hat\kappa^{\check\nu}_T]=1$, proving that $\hat\kappa^{\check\nu}$ is a $\check\Prob$-martingale. We can thus define a corresponding probability measure, $\check\Prob^{\check\nu}$, on $(\hat\Omega,\hat\mcF)$ as $d\check\Prob^{\check\nu}:=\hat\kappa^{\check\nu}_Td\check\Prob$, and since $\hat M_T\hat\kappa^{\check\nu}_T\equiv 1$, we conclude that $\check\Prob^{\check\nu}=\hat\Prob$ on $(\hat\Omega,\hat\mcF)$. We further extend this definition by letting $d\check\Prob^{\nu}:=\hat\kappa^{\nu}_Td\check\Prob$ whenever $\nu\in\check\mcV_{\inf>0}$. Here, $\check\mcV_{\inf>0}$ is the set of all $\hat\bbF^{\hat W}\vee \hat\bbF^{\check\mu}$-predictably measurable bounded maps $\nu=\nu_t(\hat\omega,e):[0,T]\times\hat\Omega\times U\to [0,\infty)$ with $\inf \nu_t>0$. In particular, with
\begin{align*}
\check X_s&=x_0+\int_0^s a(r,\check X)dr+\int_0^s\sigma(r,\check X)d\hat W_r+\int_0^s\!\!\!\int_U\gamma(r-,\check X,e)\check\mu(dr,de)
\end{align*}
we get that $\hat J_t(\check u,\tau)=\check J^{\mcR}_t(\check\nu,\tau)$, $\check\Prob$-a.s., for each $\tau\in\hat\mcT^{\hat W,\check u}$, where
\begin{align*}
\check J^{\mcR}_t(\nu,\tau) := \check\E^{\nu}\Big[\Psi(\tau,\check X)+\int_t^{\tau} f(r,\check X)dr+\int_t^{\tau} \!\!\!\int_U\chi(s-,\check X,e) \check\mu(ds,de)\Big|\hat\mcF_t\Big]
\end{align*}
and $\check\E^{\nu}$ is expectation with respect to $\check\Prob^v$. Combined with \eqref{ekv:relate-lem}, this gives us that
\begin{align}\label{ekv:below-underline-Y}
  \hat\E\big[(\check J^{\mcR}_t(\check\nu,\tau)- \hat{\underline Y}^k_t)^+\big]\leq \varrho
\end{align}
for each $\tau\in\hat\mcT^{\hat W,\Xi^{\eps}(\check u),\eps}_t$.

Since the probability space $(\hat\Omega,\hat\mcF,\check\Prob,\hat W,\check\mu)$ is a setting for our penalized BSDEs \eqref{ekv:rbsde-pen}, there is a unique quadruple $(\check Y^n,\check Z^n,\check V^n,\check K^{-,n})\in \check\mcS^2\times\check\mcH^2(\hat W)\times\check\mcH^2(\check\mu)\times \check\mcA^2$, where $\check\mcS^2$, $\check\mcH^2(\hat W)$, $\check\mcH^2(\check \mu)$ and $\check\mcA^2$ are defined as $\mcS^2$, $\mcH^2(W)$, $\mcH^2(\mu)$ and $\mcA^2$ but on the probability space $(\hat\Omega,\hat\mcF,\check\Prob,\hat W,\check\mu)$, such that
\begin{align}\label{ekv:rbsde-pen-hat}
  \begin{cases}
    \check Y^n_s=\Psi(T,\check X)+\int_s^T f(r,\check X)dr-\int_s^T \check Z^n_r d\hat W_r-\int_s^T\!\!\!\int_U \check V^n_r(e)\check\mu(dr,de)+\check K^{-,n}_T-\check K^{-,n}_s
    \\
    \quad-n\int_s^T\!\!\!\int_U(\check V^n_r(e)+\chi(r,\check X,e))^-\lambda(de)dr,\quad \forall s\in [0,T],\\
    \check Y^n_s\geq \Psi(s,\check X),\, \forall s\in [0,T]\mbox{ and } \int_0^T \big(\check Y^n_s-\Psi(s,\check X)\big)d\check K^{-,n}_s=0.
  \end{cases}
\end{align}
By Proposition~\ref{prop:pen-is-imp}, we also have the representation
\begin{align}
\check Y^n_t&=\esssup_{\tau\in\hat\mcT^{\hat W,\check\mu}_t}\essinf_{\nu\in\check\mcV^n_{\inf>0}}\check J^{\mcR}_t( \nu,\tau),\label{ekv:Yn-repr-hat}
\end{align}
employing the obvious notation $\check\mcV^n_{\inf>0}:=\{\nu\in\check\mcV_{\inf>0}:\nu\leq n\}$. Now, for each $\omega\in\Omega$ and $\omega'''\in\Omega'''$, the random measures $\check\mu(\omega,\omega''',\cdot,\cdot)$ and $\mu(\omega,\cdot,\cdot)$ agree on any measurable subset of $[0,t]\times U$ and since $\msL_{\Prob}(X,\mu,W)=\msL_{\check\Prob}(\check X,\check \mu,\hat W)$ we conclude that $\check Y^n_t=\hat Y^n_t$, $\check\Prob$-a.s.

\bigskip

\textbf{Step 2:} Using once more that $(\hat\Omega,\hat\mcF,\check\Prob,\hat W,\check\mu)$ is a setting for our penalized BSDEs \eqref{ekv:rbsde-pen}, Lemma~\ref{lem:stop-proc-disc} guarantees the existence of a non-increasing sequence $(\check\tau^\eps_n)_{n\in\bbN}$ of stopping times in\footnote{Actually, by Lemma~\ref{lem:stop-proc-disc} these are $\bbT^\eps\cap[t,T]$-valued stopping times with respect to the $\check\Prob$-completion of the filtration $(\mcF^{t,W}_s\times \Omega''')_{s\geq 0}\vee \hat\bbF^{\Xi^\eps(\check u)}$. However, since $\check\Prob$ and $\hat\Prob$ are equivalent, this set of stopping times coincides with $\hat\mcT^{\hat W,\Xi^{\eps}(\check u),\eps}_t$.} $\hat\mcT^{\hat W,\Xi^{\eps}(\check u),\eps}_t$ such that
\begin{align}\label{Yn-rel-varrho-2}
 \check\E\big[(\hat Y^n_t- \essinf_{\nu\in\check\mcV^{n,k,t}_{\inf>0}}\check J^{\mcR}_t(\nu,\check\tau^\eps_n))^+\big]\leq \varrho
\end{align}
for each $n\in\bbN$. As the sequence is non-increasing (outside of a $\check\Prob$-null set), there is a stopping time $\check\tau^{*,\eps}\in\hat\mcT^{\hat W,\Xi^{\eps}(\check u),\eps}_t$, such that $\check\tau^\eps_n\searrow \check\tau^{*,\eps}$. Moreover, since the stopping times $(\check\tau^\eps_n)_{n\in\bbN}$ take values in a finite set, the sequence is eventually constant. Hence, the sets $\check A^\eps_n:=\{\omega\in\Omega:\check\tau^{*,\eps}<\check\tau^\eps_n\}$ form a monotone sequence, the limit of which is a $\check\Prob$-null set. Now, there is a $C>0$ such that for each $\eps>0$ and $n\in\bbN$, we have
\begin{align*}
\check J^{\mcR}_t(\nu,\check \tau^\eps_n)-\check J^{\mcR}_t(\nu,\check \tau^{*,\eps}) &= \check \E^\nu\Big[\ett_{\check A^\eps_n}\big(\Psi(\check \tau^\eps_n, \check X)-\Psi(\check \tau^{*,\eps},\check X)+\int_{\check \tau^{*,\eps}}^{\check \tau^\eps_n} f(r,\check X)dr
  \\
  &\quad+\int_{\check \tau^{*,\eps}}^{\check \tau^\eps_n} \!\!\int_U\chi(s-,\check X,e) \check\mu(ds,de)\big)\Big|\hat\mcF_t\Big]
  \\
  &\leq C\check \E^\nu\big[\ett_{\check A^\eps_n}(1+\hat N^{\check \mu}_{t,T})(1+\|\check X\|_T^q)\big|\hat\mcF_t\big]
  \\
  &\leq C\check \E^\nu\big[\ett_{\check A^\eps_n}\big|\hat\mcF_t\big]^{1/2}(1+\|\check X\|^{q}_t),
\end{align*}
$\check\Prob$-a.s., for all $\nu\in\check\mcV^{n,k,t}_{\inf>0}$. It thus follows from \eqref{Yn-rel-varrho-2} that
\begin{align}
 &\check\E\big[(\hat Y^n_t- \essinf_{\nu\in\check\mcV^{n,k,t}_{\inf>0}}\{\check J^{\mcR}_t(\nu,\check \tau^{*,\eps}) + C\check \E^\nu\big[\ett_{\check A^\eps_n}\big|\hat\mcF_t\big]^{1/2}(1+\|\check X\|^{q}_t)\})^+\big]\leq\varrho.\label{ekv:Yn-bnd-2}
\end{align}

\bigskip

\textbf{Step 3:} For each $n\in\bbN$, we define the truncation of $\check\nu$ as $\check\nu^n:=\check\nu\wedge n$ and thus have that $\check\nu^n\in\check\mcV^{n,k,t}_{\inf>0}$ (see Remark~\ref{rem:check-nu-n}). Since $\hat \Prob$ and $\check\Prob$ agree on $\mcF_t$, this effectively gives us through \eqref{ekv:below-underline-Y} and \eqref{ekv:Yn-bnd-2} that
\begin{align}
\hat \E\big[(\hat Y^n_t - \hat{\underline Y}^k_t)^+\big]&\leq \varrho +\check\E\big[(\essinf_{\nu\in\check\mcV^{n,k,t}_{\inf>0}}\{\check J^{\mcR}_t(\nu,\check \tau^{*,\eps}) + C\check\E^\nu\big[\ett_{\check A^\eps_n}\big|\hat\mcF_t\big]^{1/2}(1+\|\hat X\|^{q}_t)\} - \hat{\underline Y}^k_t)^+\big]\nonumber
\\
&\leq \varrho +\check\E\big[(\check J^{\mcR}_t(\check \nu^n,\check \tau^{*,\eps}) + C\check\E^{\check \nu^n}\big[\ett_{\check A^\eps_n}\big|\hat\mcF_t\big]^{1/2}(1+\|\hat X\|^{q}_t)-\hat{\underline Y}^k_t)^+\big]\nonumber
\\
&\leq\check \E\big[(\check J^{\mcR}_t(\check\nu^n ,\check \tau^{*,\eps})-\check J^{\mcR}_t(\check\nu ,\check \tau^{*,\eps}))^+ + C\check\E^{\check\nu^n}\big[\ett_{\hat A^\eps_n}\big|\hat\mcF_t\big]^{1/2}(1+\|\hat X\|^q_t)\big]+2\varrho
\label{ekv:Yn-repr-hat-2}.
\end{align}
It remains to show that the expected value on the right-hand side tends to 0 as $n\to\infty$ to conclude that $\hat \E\big[(\hat Y^n_t- \hat{\underline Y}^k_t)^+\big]\to 0$ as $n\to \infty$, since $\varrho>0$ was arbitrary. Letting
\begin{align*}
  \Phi_t(\tau):=\Psi(\tau,\check X)+\int_t^{\tau} f(r,\check X)dr+\int_t^{\tau} \!\!\!\int_U\chi(s-,\check X,e) \check\mu(ds,de),
\end{align*}
we have
\begin{align*}
  |\Phi_t(\tau)|\leq C(1+\hat N_{t,T}^{\check u})(1+\|\check X\|^q_T)=:\bar\Phi.
\end{align*}
and get that
\begin{align*}
  \check J^{\mcR}_t(\check\nu^n,\check \tau^{*,\eps})-\check J^{\mcR}_t(\check\nu,\check \tau^{*,\eps})&\leq \esssup_{\tau\in\hat\mcT^{\hat W,\check\mu,\eps}_t}(\check J^{\mcR}_t(\check\nu^n,\tau)-\check J^{\mcR}_t(\check\nu,\tau))
  \\
  &= \esssup_{\tau\in\hat\mcT^{\hat W,\check\mu,\eps}_t}\check\E\big[(\hat\kappa^{\check\nu^n}_T-\hat\kappa^{\check\nu}_T)\Phi_t(\tau)\big|\hat\mcF_t\big]
  \\
  &\leq \check\E\big[|\hat\kappa^{\check\nu^n}_T-\hat\kappa^{\check\nu}_T|\bar\Phi\big|\hat\mcF_t\big].
\end{align*}
We will show that the right-hand side tends to 0, $\check\Prob$-a.s., then dominated convergence implies that the first term on the right-hand side of \eqref{ekv:Yn-repr-hat-2} tends to 0. Letting $E_K:=\{\omega:\|\check X\|_T+\hat N^{\check u}_{t,T}\leq K\}$, we get
\begin{align*}
  \check\E\big[|\hat\kappa^{\check\nu^n}_T-\hat\kappa^{\check\nu}_T|\bar\Phi\big|\hat\mcF_{t}\big]\leq \check\E\big[\ett_{E_K}|\hat\kappa^{\check\nu^n}_T-\hat\kappa^{\check\nu}_T|\bar\Phi\big|\hat\mcF_{t}\big] + \check\E^{\check\nu^n}\big[\ett_{E_K^c}\bar\Phi\big|\hat\mcF_{t}\big]+ \check\E^{\check\nu}\big[\ett_{E_K^c}\bar\Phi\big|\hat\mcF_{t}\big].
\end{align*}
Concerning the last two terms, we have
\begin{align*}
  \check\E^{\check\nu^n}\big[\ett_{E_K^c}\bar\Phi\big|\hat\mcF_{t}\big]&= C\check\E^{\check\nu^n}\big[\ett_{E_K^c}(1+\hat N_{t,T}^{\check u})(1+\|\check X\|^q_T)\big|\hat\mcF_{t}\big]
  \\
  &\leq \frac{C}{K}\check\E^{\check\nu^n}\big[(\|\check X\|_T+\hat N_{t,T}^{\check u})(1+\hat N_{t,T}^{\check u})(1+\|\check X\|^q_T)\big|\hat\mcF_{t}\big]
  \\
  &\leq \frac{C}{K}(1+\|\check X\|^{q+1}_{t})
\end{align*}
and similarly for the last term, where $C>0$ does not depend on $K$ or $n$ since $\check\E^{\check\nu^n}[(\hat N_{t,T}^{\check u})^4]\leq \check\E^{\check\nu}[(\hat N_{t,T}^{\check u})^4]<\infty$. Concerning the first term, we have
\begin{align*}
\check\E\big[\ett_{E_K}|\hat\kappa^{\check\nu^n}_{T}-\hat\kappa^{\check\nu}_T|\bar\Phi\big|\hat\mcF_{t}\big] &\leq C\check\E\big[|\hat\kappa^{\check\nu^n}_T-\hat\kappa^{\check\nu}_T|(1+K^{q+1})\big|\hat\mcF_{t}\big],
\end{align*}
On the other hand, by dominated convergence we have
\begin{align*}
\int_{0}^T\int_U(1-\check\nu^n_r(e))\lambda(de)dr \to \int_{0}^T\int_U(1-\check\nu_r(e))\lambda(de)dr,
\end{align*}
$\check\Prob$-a.s.~and
\begin{align*}
  \prod_{\hat\sigma_j\leq s}\check\nu^n_{\hat\sigma_j}(\hat\zeta_j)\to \prod_{\hat\sigma_j\leq s}\check\nu_{\hat\sigma_j}(\hat\zeta_j),
\end{align*}
$\check\Prob$-a.s.~as the number of terms in the product is $\check\Prob$-a.s.~finite effectively implying that $\hat\kappa^{\check\nu^n}_T\to \hat\kappa^{\check\nu}_T$, $\check\Prob$-a.s. For each $K>0$, there is thus a $\check\Prob$-null set $E\subset\hat\mcF$ such that
\begin{align*}
  \lim_{n\to\infty}\check\E\big[|\hat\kappa^{\check\nu^n}_T-\hat\kappa^{\check\nu}_T|\bar\Phi\big|\hat\mcF_{t}\big]&\leq \frac{C}{K}(1+\|\check X\|^{q+1}_{t})
\end{align*}
on $\hat\Omega\setminus E$. Since $K>0$ was arbitrary, we conclude that the left-hand side equals 0, $\check\Prob$-a.s. By dominated convergence, it follows that
\begin{align*}
  \limsup_{n\to \infty}\check \E\big[(\check J^{\mcR}_t(\check\nu^n,\check \tau^{*,\eps})-\check J^{\mcR}_t(\check\nu,\check \tau^{*,\eps}))^+\big]=0.
\end{align*}
Moreover, since $\check A^\eps_n$ is a monotone sequence of sets the limit of which is a $\check\Prob$-null set we have by dominated convergence that
\begin{align*}
  \lim_{n\to\infty}\check\E^{\check\nu^n}\big[\ett_{\check A^\eps_n}\big|\hat\mcF_{t}\big]=\lim_{n\to\infty}\check\E\big[\hat\kappa^{\check\nu^n}_T\ett_{\check A^\eps_n}\big|\hat\mcF_{t}\big]=0,
\end{align*}
$\check\Prob$-a.s. Combined, it follows that
\begin{align*}
  \lim_{n\to \infty}\check \E\big[(\hat Y^n_t-\hat{\underline Y}^{k}_t)^+\big]=0
\end{align*}
and from Fatou's lemma we conclude that $\hat \E\big[(\hat Y_t-\hat{\underline Y}^k_t)^+\big]=0$, where $\hat Y:=\lim_{n\to\infty}\hat Y^n$. The statement in the proposition now follows by projecting back on the original probability space $(\Omega,\mcF,\Prob)$ and leveraging (H.\ref{hyp:unif-conv}).\qed\\

\bibliographystyle{plain}
\bibliography{imp-stop-game_ref}
\end{document}